\pdfoutput=1
\documentclass[psamsfonts]{amsart}
\usepackage[utf8]{inputenc}
\usepackage{amsfonts, comment, enumerate}
\usepackage{mathtools}
\usepackage[mathcal]{eucal}
\usepackage[hidelinks]{hyperref}
\hypersetup{
pdftitle={Complexity of deep computations via topology of function spaces},
pdfsubject={Mathematics, Set Theory},
pdfauthor={Eduardo Dueñez, José Iovino, Tonatiuh Matos-Wiederhold, Luciano Salvetti, Franklin D. Tall},
pdfkeywords={}
}
\usepackage{amsmath,amssymb}
\usepackage{amsthm}
\usepackage{pdflscape}
\usepackage{pgfplots}
\usepackage{mathrsfs}
\usepackage{tikz}
\usetikzlibrary{trees}
\usepackage{graphicx}
\graphicspath{ {images/} }
\usepackage{subcaption}
\usepackage{stmaryrd}
\usepackage{cancel}
\usepackage{array}

\newtheorem{thm}{Theorem}[section]
\newtheorem*{AxExt}{Extendibility Axiom}
\newtheorem{cor}[thm]{Corollary}
\newtheorem{prop}[thm]{Proposition}
\newtheorem{lem}[thm]{Lemma}

\newtheorem{fact}[thm]{Fact}

\theoremstyle{definition}
\newtheorem{defn}[thm]{Definition}
\newtheorem{question}[thm]{Question}

\newtheorem{notn}[thm]{Notation}

\theoremstyle{remark}
\newtheorem{rem}[thm]{Remark}
\newtheorem{rems}[thm]{Remarks}

\DeclareMathOperator{\tp}{tp}

\makeatletter
\DeclareFontFamily{U}{tipa}{}
\DeclareFontShape{U}{tipa}{m}{n}{<->tipa10}{}
\newcommand{\arc@char}{{\usefont{U}{tipa}{m}{n}\symbol{62}}}%

\newcommand{\arc}[1]{\mathpalette\arc@arc{#1}}

\newcommand{\arc@arc}[2]{%
  \sbox0{$\m@th#1#2$}%
  \vbox{
    \hbox{\resizebox{\wd0}{\height}{\arc@char}}
    \nointerlineskip
    \box0
  }%
}
\makeatother

\newcommand{\hatC}{\hat{\mathbb{C}}}
\newcommand{\Cp}{\mathrm{C}_{\mathrm{p}}}
\newcommand{\Lshard}{\mathcal{L}_{\mathrm{sh}}}
\newcommand{\NIP}{\text{NIP}}
\newcommand{\sh}{\text{sh}}
\newcommand{\rb}{r_{\bullet}}
\newcommand{\restr}{{\restriction}}
\newcommand{\RP}{\mathbb{R}\mathbb{P}}

\newcommand{\Brack}[1]{\left\llbracket{#1}\right\rrbracket}

\makeatletter
\let\c@equation\c@thm
\makeatother
\numberwithin{equation}{section}

\title[Complexity of deep computations via topology of function spaces]{Complexity of deep computations\\ via topology of function spaces}

\author[Dueñez, Iovino, Matos-Wiederhold, Salvetti, Tall]{
Eduardo Dueñez$^{1}$ \qquad
José Iovino$^{1}$ \qquad
Tonatiuh Matos-Wiederhold$^{2}$ \qquad
Luciano Salvetti$^{2}$ \qquad
Franklin D. Tall$^{2}$
}

\date{\today}

\begin{document}
\setlength{\emergencystretch}{1.5em}
\enlargethispage{3pt}

\maketitle

{\centering\tiny\vspace{-0.6cm}
  $^{1}$Department of Mathematics, University of Texas at San Antonio\\
  $^{2}$Department of Mathematics, University of Toronto\\
}

\begin{abstract}
We bring Shelah's theory of independence and Todorčević's classification of Rosenthal compacta to bear on limits of computations beyond the continuous case.
The first paper of this series characterized computability through continuous extendibility in a framework for floating-point computation.
Limits of continuous computations may nevertheless be discontinuous, as familiar iterative algorithms demonstrate.
We identify the regime in which these limits remain of the first Baire class and apply its structure theory to classify their complexity.
For countably many real-valued features, and assuming continuous extendibility of computations, the absence of the independence property on the prescribed invariant bounded regions is equivalent to first-Baire-class regularity of all limits of the extended computations on the corresponding completed regions.
Every such limit is consequently the pointwise limit of a sequence of extended computations on each completed region.
For countable families acting continuously on a compact metrizable state space, this is the tame regime of topological dynamics, whose hierarchy Glasner and Megrelishvili obtained from Todorčević's trichotomy.
We transport that classification to computation structures and realize every level by a rotation, the power method, the planar and spatial Bresenham algorithms, and the action generated by Euclidean steps on the projective line.
Newton's method for the square roots of unity on the Riemann sphere lies beyond the tame regime.
Under the same extendibility assumption, a measure-theoretic counterpart for arbitrary feature families characterizes the absence of independence by featurewise universal measurability of all limits of extended computations on the corresponding completed regions.
For every Radon input probability measure, each output feature is then uniformly approximable by polynomials in the input features outside a set of arbitrarily small probability.
\end{abstract}

\section*{Introduction}

This paper brings Shelah's theory of independence and Todorčević's classification of Rosenthal compacta to bear on the complexity of deep computations, the limits of computations along ultrafilters; with these tools we determine when a deep computation is the limit of a sequence of computations, and we stratify deep computations into a strict hierarchy of complexity classes.
The distinction between tame and wild behavior is a central concern of \emph{classification theory}, the branch of model theory that studies theories through structural dividing lines.
Shelah's independence property supplies the dividing line used here.
The first paper of this series characterized the continuous case within a framework for floating-point computation; here we proceed to limits that may be discontinuous.
A pointwise limit of continuous computations may be discontinuous, yet remain of the first Baire class: it may still be the pointwise limit of a sequence of continuous maps.
For a countable system of continuous maps on a compact metrizable space, its \emph{enveloping semigroup} is the pointwise closure of its transformations.
The system is \emph{tame} when every map in this closure is of the first Baire class.
We bring that viewpoint to computation, using the Bourgain--Fremlin--Talagrand theorem to identify the tame regime and Todorčević's structure theory to distinguish levels within it.

The illustrations are drawn from familiar algorithms of numerical analysis, computer graphics and arithmetic: Newton's method, the power method, Bresenham's line algorithm and the Euclidean algorithm.
These classical procedures provide concrete instances of the structural distinctions studied here.
Newton's method for square roots on the real projective line and the power method have discontinuous limiting maps that are nevertheless attained sequentially.
Newton's method on the Riemann sphere, by contrast, has limits beyond the first Baire class.
These examples, developed in Sections~\ref{sec:examples} and~\ref{sec:iterative}, place the failure of continuity within a classification of limiting computations.
Limits also enter the recurrent, equilibrium and continuous-depth architectures of deep learning.
A recurrent network repeatedly applies a layer with fixed weights, while the deep equilibrium models of Bai, Kolter and Koltun~\cite{Bai:2019} seek a fixed point by root finding.
The neural ordinary differential equations of Chen, Rubanova, Bettencourt and Duvenaud~\cite{Chen:2018} replace discrete layers by continuous evolution of the hidden state.
Related efforts include the physics-informed neural networks of Raissi, Perdikaris and Karniadakis~\cite{Raissi-Perdikaris-Karniadakis:2019}, which incorporate differential equations into training.
These architectures motivate a framework in which limits, as well as finite compositions, are objects of mathematical analysis.

Alva, Dueñez, Iovino and Walton introduced compositional computation structures in the first paper of this series~\cite{alva2024approximability}, to appear in \emph{Mathematical Structures in Computer Science}.
Computations act on states whose real-valued \emph{features} are observed and manipulated to finite precision.
The \emph{type} of a state records the values of all its features; adjoining the limits of these records gives the \emph{space of types}.
The analysis takes place on \emph{shards}, the sets of states on which each feature obeys a prescribed bound.
The closure of a shard in the space of types is its \emph{type shard}, which includes its limit states.
A \emph{deep computation} is a pointwise limit of computations along an ultrafilter, evaluated on actual states.
An \emph{ultracomputation} is obtained by taking such a limit after extending the computations to the type shards.
Our classification concerns these limits of extended computations; for the compact dynamical structures used in the algorithmic examples, every type is realized and the two notions coincide.

The first paper treats the \emph{stable} sets of computations, whose limits on the type shards remain continuous.
Its principal tool is Grothendieck's double limit criterion~\cite{Grothendieck:1952}.
For a uniformly bounded family of continuous functions on a compact space, the criterion characterizes continuity of every pointwise limit by the exchange of iterated limits over sequences of functions and points, whenever both iterated limits exist.
In the computation framework, this criterion yields polynomial approximability: on each shard, every output feature of a limiting computation is a uniform limit of polynomials in the input features.
This is the notion of computability established in that paper.
The passage from continuous to Baire-class-1 limits retains sequential approximation while admitting the discontinuities that arise in natural algorithms.
Its analysis calls for the structure theory of Rosenthal compacta in place of the continuous-function criterion.

Beyond stability, Shelah's classification theory~\cite{Shelah:1990} singles out \emph{Shelah's independence property}, which he isolated in 1971~\cite[Definition~4.1]{Shelah:1971}.
A formula has the independence property when, for every natural number~$n$, there are $n$ parameters on which every pattern of truth and falsity of the formula is realized by some element of a model.
A theory none of whose formulas has the independence property is called \emph{NIP}, or \emph{dependent}.
Within a few years, the same dividing line appeared independently in two other fields, as a dichotomy between polynomial and exponential growth.
A family of sets shatters a finite set when it cuts out all subsets of that set as traces, that is, as intersections of the set with members of the family.
In these terms, either a family shatters arbitrarily large finite sets, or the number of its traces on every finite set is bounded by a polynomial in the size of the set.
In statistical learning, Vapnik and Chervonenkis announced this dichotomy in 1968~\cite{Vapnik-Chervonenkis:1968} and published its proof in 1971~\cite{Vapnik-Chervonenkis:1971}, in their study of the uniform convergence of relative frequencies to probabilities.
In combinatorics, Sauer~\cite{Sauer:1972} and Shelah~\cite{Shelah:1972} proved the dichotomy in 1972, and it is now known as the Sauer--Shelah lemma.
The largest size of a finite set shattered by a family is called the \emph{Vapnik--Chervonenkis dimension} of the family.

The distinction between tame and wild behavior also animates infinite Ramsey theory, where partition and subsequence theorems extract coherent structure from arbitrary configurations.
Rosenthal's dichotomy gives this distinction an analytic form: a uniformly bounded sequence of real-valued functions on a set has either a pointwise convergent subsequence or an independent subsequence.
Here independence means that there are a lower and an upper threshold such that, for any two disjoint finite sets of indices, some point makes the functions indexed by the first fall below the lower threshold and those indexed by the second exceed the upper one.
This dichotomy lies at the heart of Rosenthal's $\ell_1$-theorem of 1974~\cite{Rosenthal:1974}: every bounded sequence in a Banach space has a subsequence that is either weakly Cauchy or equivalent to the unit vector basis of~$\ell_1$, the space of absolutely summable sequences.
The Bourgain--Fremlin--Talagrand theorem~\cite{BFT_1978_PCompactBaire} carries this analysis from sequences to pointwise closures of families of functions, making the first Baire class the natural domain of the tame alternative.
For pointwise bounded families of continuous functions on a Polish space, Bourgain, Fremlin and Talagrand characterized in 1978 the requirement that every function in the pointwise closure belong to the first Baire class (see Theorem~\ref{thm:BFT}).
Recall that a Polish space is a separable and completely metrizable space.
On a compact metrizable space, their condition is exactly the absence of the independence property for the family of functions (see Theorem~\ref{thm:BFT2}).
Compact sets of functions of the first Baire class on Polish spaces, with the topology of pointwise convergence, were investigated by Rosenthal in 1977~\cite{Rosenthal:1977}.
Such compacta were named \emph{Rosenthal compacta} after him by Godefroy in 1980~\cite{Godefroy:1980}, and have since been studied intensively in functional analysis, especially in the theory of Banach spaces.

The classification of Rosenthal compacta took a decisive form in 1999, when Todorčević~\cite{Todorcevic_1999_CompactSubsetsBaire} proved his trichotomy.
Every separable Rosenthal compactum that is not metrizable contains one of a few canonical compacta.
The alternatives are governed by hereditary separability and first countability (see Theorem~\ref{thm:Todorcevic-trichotomy}).
The theorem supplies a structural classification of the ways in which such compacta can fail to be metrizable.
Its proof combines infinite Ramsey theory with set-theoretic forcing: the compacta are analyzed through their interpretations in forcing extensions, and absoluteness arguments allow that analysis to yield conclusions in the ground universe.
These methods establish the structural alternatives in ZFC.
Todorčević's methods continue to be analyzed and exploited~\cite{Todorcevic:2023,horowitz2019compactsetsbaireclass}.

Glasner and Megrelishvili~\cite{glasner2022tame} brought this structure theory into topological dynamics.
For a system on a compact metrizable space, its \emph{enveloping semigroup} is the pointwise closure of its family of transformations.
The system is \emph{tame} when every map in this closure is of the first Baire class~\cite{Kohler1995,Glasner2006}.
The Glasner--Megrelishvili hierarchy distinguishes tame systems by metrizability, hereditary separability and first countability of their enveloping semigroups.
This dynamical meaning of tameness governs our use of the term for computations.
The absence of the independence property identifies the tame regime, and Todorčević's trichotomy makes possible a finer classification within it.

Model theory provides the connection between these function spaces and computation structures.
The passage from states to types places computations and their limits in spaces of functions with the topology of pointwise convergence, the setting of \emph{$\Cp$-theory}.
Shelah's independence property becomes a condition on the finite patterns of high and low output values realized by computations on suitable inputs.
For countably many features, its absence, denoted by \emph{NIP}, corresponds under the extendibility and confinement hypotheses to first-Baire-class regularity on the type shards.
This correspondence links the classification of Rosenthal compacta with the model-theoretic passage beyond stability.
It also connects the framework with Vapnik--Chervonenkis (VC) dimension and PAC (probably approximately correct) learning~\cite{Vapnik-Chervonenkis:1971,Vapnik-Chervonenkis:1974,Valiant:1984}: finiteness of VC dimension implies NIP, and the two notions coincide for first-order formulas (see Section~\ref{sec:NIP-poly-exp}).

Our first main result, Theorem~\ref{thm:nip-baire1def}, concerns computation structures satisfying the Extendibility Axiom and having countably many features; the axiom requires computations to extend continuously to the type shards.
Fix a collection of shards that covers the actual states and that every computation of the family under consideration preserves.
Applying the Bourgain--Fremlin--Talagrand theorem, we characterize NIP on these shards, feature by feature, by the requirement that every limit of the extended computations be of Baire class~1 on the corresponding type shards.
The pointwise closures of these families of extended computations are then Rosenthal compacta, and every member is the pointwise limit of a sequence of extended computations on the chosen type shard.
For countable families of computations these compacta are separable, so Todorčević's trichotomy applies.
We distinguish the levels $\NIP_3$, $\NIP_2$ and $\NIP_1$ by metrizability, hereditary separability and first countability, respectively (see Definition~\ref{def:NIP-levels}).
Together with stability and NIP, they form the hierarchy
\[
\text{stable}\Rightarrow\NIP_3\Rightarrow\NIP_2\Rightarrow\NIP_1\Rightarrow\NIP.
\]
For the compact metrizable dynamical structures considered below, this is the Glasner--Megrelishvili hierarchy; the type shards allow the same analysis in the more general computation framework.

The topological distinctions admit a precise interpretation in terms of observations (see Proposition~\ref{prop:levels-tests}).
A \emph{test} consists of an input state type and a feature evaluated at the output.
At the metrizable level, one countable collection of tests separates all the limiting computations on a type shard.
First countability allows the determining collection to depend on the individual limiting computation.
For these separable Rosenthal compacta, hereditary separability excludes a family of the cardinality of the continuum in which each member can be isolated from the other members by finitely many tests and a positive separation margin.
The hierarchy therefore distinguishes the organization of the observations needed to separate limiting computations, as well as their regularity as maps.

The familiar algorithms introduced above realize the distinctions throughout the hierarchy.
In Section~\ref{sec:iterative} we analyze their associated dynamics and show that every implication in the hierarchy is strict.
An irrational rotation of the circle is stable, whereas the power method is $\NIP_3$ but not stable.
Bresenham's planar line algorithm, represented by its Sturmian dynamics, is $\NIP_2$ but not $\NIP_3$.
The corresponding construction in three dimensions is $\NIP_1$ but not $\NIP_2$.
The action on the projective line of the group generated by Euclidean steps is NIP but not $\NIP_1$; its limiting computations include exact equality tests.
The logistic family supplies further examples within and beyond the hierarchy, while Newton's method on the Riemann sphere lies beyond NIP for the polynomials considered here.
These classifications concern the spaces of limiting maps generated by the procedures, and need not follow from the apparent simplicity of a single computational step.

{\emergencystretch=2.5em
The examples also distinguish the canonical obstructions supplied by Todorčević's theorem.
At the level $\NIP_1$ but not $\NIP_2$, the three-dimensional Bresenham construction realizes the extended duplicate of the split Cantor set and contains no Alexandroff duplicate of the Cantor set.
A selected countable family of compositions of Euclidean steps and their inverses realizes the other alternative, containing the Alexandroff duplicate (see Proposition~\ref{prop:euclid-duplicate}).
Our examples at the level NIP but not $\NIP_1$ arise from actions of groups with two generators.
Whether a single iterated continuous map can reach that level remains open, as does the corresponding realization of the Alexandroff-duplicate alternative at the preceding level (see Questions~\ref{q:cascade-NIP1} and~\ref{q:cascade-duplicate}).\par}

For families indexed by the binary tree, we also use the refinement of the structure theory due to Argyros, Dodos and Kanellopoulos~\cite{argyros2008rosenthal}.
Their seven canonical minimal families retain the tree indexing in the comparison of pointwise closures.
In Section~\ref{sec:heptachotomy} we describe these families and interpret their constituent prefix and threshold tests as primitives of binary search.
This refinement gives a more detailed analysis of the canonical families within the tame regime.

Finally, in Section~\ref{sec:measure-theoretic-NIP} we pass from topology to measure, dispensing with countability of the feature family.
The measure-theoretic form of the Bourgain--Fremlin--Talagrand theorem characterizes NIP by universal measurability of all limits of extended computations on the prescribed invariant type shards (see Theorem~\ref{thm:NIP-Monte-Carlo}).
Universal measurability is taken feature by feature with respect to every Radon probability measure on the type shard.
For each such measure and each output feature, Luzin's theorem gives continuity on a closed subset whose complement has arbitrarily small measure.
The Stone--Weierstrass theorem then gives uniform polynomial approximation in the input features on that subset.
This approximation property is the consequence of universal Monte Carlo computability established here; the definition requires a featurewise universally measurable extension.
We obtain the same approximation conclusion from universal Talagrand stability and show that, for countable families of computations, this condition coincides with NIP (see Corollary~\ref{cor:NIP-Talagrand-computations}).

The paper presupposes no acquaintance with advanced topology, model theory or the theory of computation.
The necessary background is assembled in Section~\ref{sec:prelim}, beginning from undergraduate topology and measure theory.
The purpose of that development is to make the structural results available in their natural generality while keeping their computational applications in view.

\subsection*{Organization of the paper}
In Section~\ref{sec:prelim} we present the topological and combinatorial preliminaries,
including five canonical Rosenthal compacta (Section~\ref{sec:five-compacta})
and the role of $\NIP$ in model theory, combinatorics and learning (Section~\ref{sec:NIP-poly-exp}).
In Section~\ref{sec:CCS} we recall compositional computation structures, the source of the structural and model-theoretic viewpoint of the paper, and their stable sets of computations.
That section ends with two motivating examples, Newton’s method for square roots and Newton’s method on the Riemann sphere, whose different behaviors the classification is designed to separate (Section~\ref{sec:preview}).
Section~\ref{sec:Classification} contains the classification of deep computations and the main theorems.
In Section~\ref{sec:examples} we work out Newton’s method on the Riemann sphere and on the real projective line.
Next, in Section~\ref{sec:iterative}, we place classical algorithms at every level of the hierarchy and beyond it.
In Section~\ref{sec:measure-theoretic-NIP} we replace topology by measure and study deep computations that are universally measurable on each type shard.
The output features of such deep computations are therefore, outside sets of arbitrarily small probability, uniform limits of polynomials in the input features.
We intend to develop this measure-theoretic viewpoint further in future research.
Section~\ref{sec:open-questions} contains the questions that remain open.
In Appendix~\ref{app:thue-morse} we prove that the shift on the Thue--Morse subshift generates an action of~$\mathbb{Z}$ that is not tame, a fact used in Section~\ref{sec:beyond-NIP}.

\tableofcontents

\section{Topological preliminaries:\\
  From continuity to Baire class~1}
\label{sec:prelim}

In this section we assemble the topology of function spaces on which the rest of the paper rests.
The two themes of this paper are the approximation of a deep computation by a sequence of computations and the complexity of the set of all deep computations.
Since deep computations are pointwise limits of continuous maps along ultrafilters, both themes belong to the theory of spaces of functions with the topology of pointwise convergence.
We first recall Polish spaces, Baire category and the first Baire class;
we then turn to Rosenthal's dichotomy and the Bourgain--Fremlin--Talagrand dichotomy, from which the NIP emerges.
The Rosenthal compacta come next, and five of them serve as yardsticks for the classification.
The section closes with the role of the NIP in model theory, in combinatorics and in statistical learning.
All spaces are assumed to be Hausdorff unless explicitly stated otherwise.

We write $\mathbb{N}$ for the set $\{0,1,2,\dots\}$ of natural numbers, so that the semigroups of iterates considered below contain the identity.
A subset of a topological space is $F_\sigma$ if it is a countable union of closed sets, and $G_\delta$ if it is a countable intersection of open sets.
In a metric space every open set is~$F_\sigma$, and, equivalently, every closed set is~$G_\delta$.

A \emph{Polish space} is a separable topological space whose topology is induced by a complete metric.
Countable products of Polish spaces, with the product topology, are Polish and furnish classical examples besides the real line.
They include the Cantor space~$2^\mathbb{N}$ of all infinite binary sequences, the Baire space $\mathbb{N}^\mathbb{N}$ of all infinite sequences of natural numbers, and the space $\mathbb{R}^\mathbb{N}$ of all sequences of real numbers.

\begin{fact}\cite[4.3.24]{Engelking:1989}
  A subset of a Polish space is itself Polish in the subspace topology if and only if it is a $G_{\delta }$ set.
  In particular, closed subsets and open subsets of Polish spaces are also Polish spaces.
\end{fact}

A subset $A$ of a topological space $X$ is \emph{nowhere dense} if its closure has empty interior.
The set $A$ is of the \emph{first category}, or \emph{meager}, if it is contained in a countable union of closed nowhere dense sets, and of the \emph{second category} otherwise.
The complements of meager sets are called \emph{comeager}.
Since the closure of a nowhere dense set is nowhere dense, the meager subsets of $X$ form the $\sigma$-ideal generated by the nowhere dense subsets of $X$.
A topological space $X$ is \emph{Baire} if every comeager subset of $X$ is dense, or, equivalently, if every countable intersection of dense open sets is dense.
The Baire category theorem states that compact Hausdorff spaces and completely metrizable spaces, Polish spaces in particular, are Baire.%
\footnote{More generally, any Čech-complete space is Baire~\cite[Theorem 3.9.3]{Engelking:1989}.}

For topological spaces $X$ and $Y$, we denote by $\Cp(X,Y)$ the set of all continuous functions $f\colon X\to Y$ with the topology of pointwise convergence, and we abbreviate $\Cp(X,\mathbb{R})$ to $\Cp(X)$.
The topology of pointwise convergence is the one that $\Cp(X,Y)$ inherits as a subspace of the product $Y^X$, in which $X$ serves merely as an index set.

A function $f\colon X\to Y$ between topological spaces is of \emph{Baire class~1} if it is the pointwise limit of a sequence $(f_n)_{n\in\mathbb{N}}$ of continuous functions, that is, $f(x)=\lim_{n\to\infty}f_n(x)$ for all $x\in X$.
When $Y$ is metrizable, an elementary argument shows that a function of Baire class~1 is continuous at every point outside a set of the first category.
We denote by $\mathrm{B}_1(X,Y)$ the set of functions $X\to Y$ of Baire class~1 with the topology of pointwise convergence, and we write $\mathrm{B}_1(X)$ for $\mathrm{B}_1(X,\mathbb{R})$;
the inclusions $\Cp(X,Y)\subseteq \mathrm{B}_1(X,Y)\subseteq Y^X$ are topological.
Relatively compact subsets of $\Cp(X)$ and of $\mathrm{B}_1(X)$, especially for Polish $X$, are central objects of study in functional analysis and in topological dynamics.

A topological space $X$ is \emph{perfectly normal} if it is normal and every closed subset of $X$ is a $G_\delta$ set (equivalently, every open subset of $X$ is an $F_\sigma$ set).
Every metrizable space (hence, every Polish space) is perfectly normal.

The following characterization of Baire class~1 functions is classical; for a proof, see~\cite[Section 10]{Todorcevic_1997_TopicsTop}.

\begin{fact}[Lebesgue, Hausdorff, Baire\protect\footnotemark]
  \label{baire}
  If $X$ is perfectly normal, then the following conditions are equivalent for a function $f:X\to \mathbb{R}$:
  \begin{itemize}
  \item[(1)]
    The function $f$ is of Baire class~1, that is, it is a pointwise limit of continuous functions.
  \item[(2)]
    The preimage $f^{-1}[U]$ is an $F_\sigma$ subset of $X$ whenever $U\subseteq\mathbb{R}$ is open.
  \end{itemize}
  If $X$ is Polish, then (1) and (2) are equivalent to:
  \begin{itemize}
  \item[(3)]
    For every nonempty closed $F\subseteq X$, the restriction $f|_F$ has a point of continuity.
  \end{itemize}
  Moreover, if $X$ is Polish and $f\notin \mathrm{B}_1(X)$, then there exist countable $D_0,D_1\subseteq X$ and reals $a<b$ such that
  \[
    D_0\subseteq f^{-1}(-\infty,a],
    \quad
    D_1\subseteq f^{-1}[b,\infty),
    \quad
    \overline{D_0}=\overline{D_1}.
  \]
\end{fact}
\footnotetext{The implication (1)$\Rightarrow$(2) holds in every topological space.
Indeed, if $f$ is the pointwise limit of a sequence $(f_n)_{n\in\mathbb{N}}$ of continuous functions and $a$ is a real number, then
\[
  f^{-1}(a,\infty)=\bigcup_{k\ge1}\bigcup_{m\in\mathbb{N}}\bigcap_{n\ge m}f_n^{-1}[a+1/k,\infty)
\]
is a countable union of closed sets.
The same is true of $f^{-1}(-\infty,a)$;
since every open subset of~$\mathbb{R}$ is a countable union of intervals $(a,b)=(a,\infty)\cap(-\infty,b)$, condition~(2) follows.
The reverse implication, (2)$\Rightarrow$(1), was established by Lebesgue for real-valued functions of finitely many real variables~\cite{Lebesgue:1904} and extended by Hausdorff to real-valued functions on metric spaces~\cite[Theorem~24.10]{kechris1995cdst}.
The proof of (2)$\Rightarrow$(1) requires of~$X$ only that it be normal and that each of its open subsets be an $F_\sigma$ set.
These two requirements are precisely the defining properties of a perfectly normal space.
The implication (1)$\Rightarrow$(3) holds whenever $X$ is Polish.
If $f$ satisfies~(1) and $F$ is a nonempty closed subset of~$X$, then the restriction $f|_F$ is of Baire class~1.
As noted before Fact~\ref{baire}, the points of discontinuity of $f|_F$ therefore form a meager subset of~$F$.
Since $F$ is closed in a Polish space, it is itself Polish, hence a Baire space that this meager set cannot exhaust.
Consequently, the restriction $f|_F$ has a point of continuity.
The reverse implication (3)$\Rightarrow$(1) is due to Baire, who proved it in his thesis for real-valued functions of real variables~\cite{Baire:1899}; for Polish spaces, see~\cite[Theorem~24.15]{kechris1995cdst}.
Finally, the last assertion of Fact~\ref{baire} follows from the implication (3)$\Rightarrow$(1) together with the Baire category theorem.}

\subsection{From Rosenthal’s dichotomy to the Bourgain--Fremlin--Talagrand dichotomy to Shelah’s NIP}

The question when a deep computation is accessible to approximation, that is, when it is the limit of a \emph{sequence} of computations, is answered by Rosenthal's dichotomy and the Bourgain--Fremlin--Talagrand dichotomy.
We recall both dichotomies in this subsection;
from them the NIP emerges in the theory of function spaces.
In a metrizable space every point of the closure of a set is the limit of a sequence of points of the set.
In spaces of functions with the topology of pointwise convergence, by contrast, the closure of a set need not admit this sequential description, and the two dichotomies determine when it does not.

On the Cantor space $X=2^\mathbb{N}$, the $n$-th coordinate map $\pi_n\colon X\to\{0,1\}$ is continuous for each $n\in \mathbb{N}$;
however, the sequence $(\pi_n)_{n\in\mathbb{N}}$ has \emph{no} convergent subsequence in $\mathbb{R}^X$ (see~\cite[Chapter~1.1]{Todorcevic_1997_TopicsTop}).
In a sense, this sequence exhibits the most extreme failure of sequential convergence.
Indeed, the closure of $\{\pi_n:n\in\mathbb{N}\}$ in $\{0, 1\}^X$, or in $\mathbb{R}^X$, is homeomorphic to the Stone--Čech compactification $\beta\mathbb{N}$ of the discrete space of natural numbers.
The space $\beta \mathbb{N}$ has cardinality $2^{2^{\aleph_0}}$, the largest possible for a separable Hausdorff space, and its points correspond to the ultrafilters on~$\mathbb{N}$.%
\footnote{\label{fn:id-ufilters}Addition on~$\mathbb{N}$ extends to an associative operation on~$\beta\mathbb{N}$, whose idempotent elements, the idempotent ultrafilters, exist by the Ellis--Numakura lemma.
Alva, Dueñez, Iovino and Walton~\cite{alva2024approximability} use idempotent ultrafilters to prove the existence of \emph{deep equilibria} in compositional computation structures.
Deep equilibria are limits of the iterates of a single computation, taken along an ultrafilter on~$\mathbb{N}$, that are idempotent under composition.}

In 1974 Haskell Rosenthal proved the following dichotomy for pointwise bounded sequences of continuous functions on a Polish space.
His theorem shows that for these sequences a copy of $\beta\mathbb{N}$, such as the closure of the coordinate maps, is the only obstruction to the existence of a convergent subsequence.

\begin{thm}[Rosenthal’s dichotomy~\cite{Rosenthal:1974}]
  \label{thm:Rosenthal}
  If $X$ is a Polish space and $(f_n)_{n\in\mathbb{N}}$ is a pointwise bounded sequence in $\Cp(X)$, then the sequence $(f_n)$ has either a convergent subsequence or a subsequence whose closure in $\mathbb{R}^X$ is homeomorphic to~$\beta\mathbb{N}$.
\end{thm}

Rosenthal's dichotomy grew out of \emph{Rosenthal’s $\ell_1$-theorem}, in whose statement $\ell_1$ denotes the space of absolutely summable sequences.
This theorem states that every bounded sequence in a Banach space has a subsequence that is either weakly Cauchy or equivalent to the unit vector basis of~$\ell_1$.
In particular, a Banach space contains no isomorphic copy of $\ell_1$ if and only if every bounded sequence in it has a weakly Cauchy subsequence.

The Bourgain--Fremlin--Talagrand theorem, which concerns sets rather than sequences, characterizes the relative compactness in the larger space $\mathrm{B}_1(X)$ of a pointwise bounded set of continuous functions on a Polish space~$X$.
For such a set, relative compactness in $\mathrm{B}_1(X)$ is the property that every function in its pointwise closure be of Baire class~1.
This theorem is a counterpart of Rosenthal’s dichotomy, although it is usually stated as an equivalence rather than as a dichotomy (see Theorem~\ref{thm:BFT} below).

The statement of the theorem requires the following notation.
For $f\in \mathbb{R}^X$ and a real number $a$, define
\begin{align*}
X^f_{\le a} & \coloneqq f^{-1}(-\infty,a],  &
X^f_{\ge a} & \coloneqq f^{-1}[a, +\infty).
\end{align*}
For a set $A\subseteq \mathbb{R}^X$ and a real number $a$, define, by a slight abuse of notation,
\begin{align*}
X^A_{\le a} & \coloneqq \bigcap_{f\in A}X^f_{\le a} = \{x\in X: \text{$f(x)\le a$ for all $f\in A$}\}, \\
X^A_{\ge a} & \coloneqq \bigcap_{f\in A}X^f_{\ge a} = \{x\in X: \text{$f(x)\ge a$ for all $f\in A$}\};
\end{align*}
in particular, we have $X^\emptyset_{\ge a} = X = X^\emptyset_{\le a}$.
For a sequence $(f_n)_{n\in\mathbb{N}}$ in $\mathbb{R}^X$ and a set $I\subseteq \mathbb{N}$, define $I^{\complement} \coloneqq \mathbb{N}\setminus I$ and $f_I \coloneqq \{f_i: i\in I\}$.

\begin{thm}[The Bourgain--Fremlin--Talagrand dichotomy~\cite{BFT_1978_PCompactBaire}]
  \label{thm:BFT}
  Let $X$ be a Polish space, and let $A\subseteq \Cp(X)$ be pointwise bounded.
  The following are equivalent:
  \begin{enumerate}[(i)]
  \item
    The set $A$ is relatively compact in $\mathrm{B}_1(X)$.
  \item
    For every $(f_n)_{n\in\mathbb{N}}\subseteq A$ and every $a<b$ there is $I\subseteq\mathbb{N}$ such that
    \begin{equation*}
    X_{\le a}^{f_I} \cap X_{\ge b}^{f_{I^{\complement}}} = \emptyset.
    \end{equation*}
  \end{enumerate}
\end{thm}
(In this form, the Bourgain--Fremlin--Talagrand dichotomy is a special case of the equivalence (ii)$\Leftrightarrow$(v) of Corollary~4G of Bourgain, Fremlin and Talagrand~\cite{BFT_1978_PCompactBaire}.)

In condition~(ii) of Theorem~\ref{thm:BFT}, the set $X_{\le a}^{f_I}$ is defined by the inequalities $f_i(x)\le a$ with $i\in I$, and the set $X_{\ge b}^{f_{I^{\complement}}}$ by the inequalities $f_j(x)\ge b$ with $j\in I^{\complement}$.
Since $I$ or $I^{\complement}$ is infinite, at least one of the two sets is defined by countably infinitely many inequalities.
A computation, however, can check only finitely many inequalities at a time, and the following definitions accordingly impose only finitely many.

\begin{defn}\label{def:NIP-fns}
  A set $A\subseteq \mathbb{R}^X$ of functions has the \emph{finitary No-Independence Property (NIP)} if, for every sequence $(f_n)_{n\in\mathbb{N}}$ in $A$ and all reals $a<b$, there exist finite disjoint sets $I,J\subseteq\mathbb{N}$ such that $X_{\le a}^{f_I} \cap X_{\ge b}^{f_J} = \emptyset$.
  We say that finite disjoint sets $I, J$ with this property \emph{witness finitary NIP for $(f_n)$ and $a, b$}.

  A set $A\subseteq \mathbb{R}^X$ has the \emph{finitary Independence Property (IP)} if it does not have finitary NIP.
  Explicitly, finitary IP means that there exist a sequence $(f_n)_{n\in\mathbb{N}}$ in $A$ and reals $a<b$ such that $X_{\le a}^{f_I} \cap X_{\ge b}^{f_J} \ne \emptyset$ for every pair of finite disjoint sets $I,J\subseteq\mathbb{N}$.
\end{defn}

If the word “finite” is omitted in these definitions, we obtain \emph{countable} NIP and \emph{countable} IP.
Countable NIP is weaker than finitary NIP, and countable IP is accordingly stronger than finitary IP.
If, in addition, the witnesses are required to satisfy $J = I^{\complement}$, we obtain the properties that we call \emph{BFT-NIP} and \emph{BFT-IP}.
Theorem~\ref{thm:BFT} thus states that, on a Polish space~$X$, a pointwise bounded set $A\subseteq \Cp(X)$ is relatively compact in $\mathrm{B}_1(X)$ if and only if it has BFT-NIP.

From now on, IP and NIP mean \emph{finitary} IP and NIP unless otherwise specified.

\begin{thm}
  \label{thm:BFT2}
  Let $X$ be a compact metrizable (hence Polish) space.
  For every pointwise bounded set $A\subseteq \Cp(X)$, the following are equivalent:
  \begin{enumerate}[(i)]
  \item
    the set $A$ is relatively compact in~$\mathrm{B}_1(X)$;
  \item
    the set $A$ has BFT-NIP;
  \item
    the set $A$ has countable NIP;
  \item
    the set $A$ has finitary NIP.
  \end{enumerate}
\end{thm}
(The equivalence \emph{(ii)}$\Leftrightarrow$\emph{(iv)} holds for every compact space~$X$, even without the assumption that $A$ is pointwise bounded, and the equivalence \emph{(ii)}$\Leftrightarrow$\emph{(iii)} holds for every set $A\subseteq\mathbb{R}^X$ on an arbitrary space~$X$.)
\begin{proof}
  The equivalence (i)$\Leftrightarrow$(ii) is Theorem~\ref{thm:BFT}.
  The implications (iv)$\Rightarrow$(iii) and (ii)$\Rightarrow$(iii) are trivial.
  The converse (iii)$\Rightarrow$(ii) holds because witnesses $I$ and $J$ of countable NIP satisfy $J\subseteq I^{\complement}$, whence $I$ and $I^{\complement}$ witness BFT-NIP.
  To prove (ii)$\Rightarrow$(iv), suppose that $(f_n)$ and $a<b$ witness finitary IP, and let $I\subseteq\mathbb{N}$ be given.
  The closed sets $X^{f_n}_{\le a}$ ($n\in I$) and $X^{f_m}_{\ge b}$ ($m\in I^{\complement}$) have the finite intersection property; by compactness, therefore, we obtain $X^{f_I}_{\le a}\cap X^{f_{I^{\complement}}}_{\ge b}\neq\emptyset$.
  Since $I$ is arbitrary, the set $A$ does not have BFT-NIP, and (ii) fails.
\end{proof}
Theorem~\ref{thm:BFT2} can be cast as a dichotomy: on a compact metrizable space~$X$, every pointwise bounded set $A\subseteq\Cp(X)$ either is relatively compact in~$\mathrm{B}_1(X)$ or has IP in the finitary and, equivalently, in the countable sense.

The independence property was isolated by Saharon Shelah~\cite[Definition 4.1]{Shelah:1971} as one of the dividing lines of his classification of first-order theories~\cite{Shelah:1990}.
We return to it in Section~\ref{sec:NIP-poly-exp}, where it appears together with its counterparts in combinatorics and in statistical learning.

\subsection{Rosenthal compacta}

The Bourgain--Fremlin--Talagrand dichotomy leads to the Rosenthal compacta, the compact spaces within which the deep computations of this paper are classified, and we recall them here.
By Theorem~\ref{thm:BFT}, if a pointwise bounded set of continuous functions on a Polish space has BFT-NIP, then its pointwise closure is a compact set of functions of the first Baire class.
Compact sets of this kind had appeared in Rosenthal's work on the $\ell_1$-theorem, and in 1977 he made them the subject of a study of their own~\cite{Rosenthal:1977}.
Godefroy~\cite{Godefroy:1980} named them after Rosenthal in the following definition.

\begin{defn}
  \label{def:Rosenthal-compacta}
  A \emph{Rosenthal compactum} is a topological space homeomorphic to a compact subset of $\mathrm{B}_1(X)$, with the topology of pointwise convergence, for some Polish space~$X$.
\end{defn}

Rosenthal compacta have strong tameness properties;
for instance, every point in the closure of a subset of a Rosenthal compactum is the limit of a sequence of points of the subset~\cite[Theorem~3F]{BFT_1978_PCompactBaire}.
In this paper Rosenthal compacta serve as spaces of deep computations.
Theorem~\ref{thm:nip-baire1def} below shows that for countably many features the deep computations of an NIP set of computations, extended to the limit states, form a Rosenthal compactum on each type shard.
Here a type shard is a compact set of states and limit states on which every feature obeys a prescribed bound (see Section~\ref{sec:shards}).

\subsection{The special case~\texorpdfstring{$\mathrm{B}_1(X,\mathbb{R}^\mathcal{P})$}{B\_1(X, R\textasciicircum P)} with~\texorpdfstring{$\mathcal{P}$}{P} countable}

This subsection is devoted to reducing the study of maps with values in a product of copies of the real line to that of real-valued functions.
The reduction is needed because a state of a computation is described by the values of its features.
A computation is therefore a map with values in the product of copies of the real line indexed by the features.
Fix a countable set~$\mathcal{P}$ (in this paper, “countable” means finite or countably infinite), whose elements we call \emph{predicate symbols} or \emph{formal predicates}.
Our goal is to characterize the relatively compact subsets of $\mathrm{B}_1(X,\mathbb{R}^\mathcal{P})$, where $X$ is a perfectly normal space, typically a Polish space.

Whenever $\mathcal{P}$ is regarded as a topological space, it carries the discrete topology, and the product $X\times \mathcal{P}$ is regarded either as a set or as a topological product, depending on the context.
The spaces $\Cp(X, \mathbb{R}^{\mathcal{P}})\subseteq \mathrm{B}_1(X, \mathbb{R}^{\mathcal{P}})$ are subspaces of the ambient space $(\mathbb{R}^{\mathcal{P}})^X$, which can be described in two other ways, through the natural homeomorphisms
\begin{equation*}
(\mathbb{R}^{\mathcal{P}})^X \cong \mathbb{R}^{X\times \mathcal{P}} \cong (\mathbb{R}^X)^{\mathcal{P}},
\end{equation*}
given by
\begin{align*}
\mathbb{R}^{X\times \mathcal{P}} \to (\mathbb{R}^{\mathcal{P}})^X: \Phi & \mapsto {\hat{}}\Phi    \\
\mathbb{R}^{X\times \mathcal{P}} \to (\mathbb{R}^X)^{\mathcal{P}}: \Phi & \mapsto \Phi{\,\hat{}},
\end{align*}
where
\begin{align*}
\hat{}\Phi(x)\coloneqq \Phi(x, \cdot)\in \mathbb{R}^{\mathcal{P}}, &  &
\Phi\,\hat{}(P)\coloneqq \Phi(\cdot, P) \in \mathbb{R}^X.
\end{align*}
These identifications treat $X$ and $\mathcal{P}$ as mere sets;
in particular, the topology of~$X$ plays no role in them.

The following proposition reduces the study of continuous functions, and of functions of Baire class~1, from $X$ to $\mathbb{R}^{\mathcal{P}}$ to the case of real-valued functions.

\begin{prop}\label{prop:Baire-identifications}
  The identification $(\mathbb{R}^{\mathcal{P}})^X \cong \mathbb{R}^{X\times \mathcal{P}} \cong (\mathbb{R}^X)^{\mathcal{P}}$ induces identifications
  \begin{align*}
  \Cp(X, \mathbb{R}^{\mathcal{P}}) \cong \Cp(X\times \mathcal{P}) \cong \Cp(X)^{\mathcal{P}}, &  &
  \mathrm{B}_1(X, \mathbb{R}^{\mathcal{P}}) \cong \mathrm{B}_1(X\times \mathcal{P}) \cong \mathrm{B}_1(X)^{\mathcal{P}}.
  \end{align*}
\end{prop}
(The cardinality of $\mathcal{P}$ plays no role here.)
\begin{proof}
  The identification of the $\Cp$-spaces follows from the definition of the product topology and the discreteness of~$\mathcal{P}$.
  By the former, a continuous map $X\to \mathbb{R}^{\mathcal{P}}$ is precisely a $\mathcal{P}$-indexed family of continuous functions $X\to \mathbb{R}$;
  by the latter, such families correspond to continuous functions $X\times \mathcal{P}\to \mathbb{R}$.
  For the spaces of functions of Baire class~1, recall that membership in the first Baire class is defined through pointwise limits of sequences of continuous functions.
  Since a sequence of $\mathcal{P}$-indexed families of continuous functions is the same as a $\mathcal{P}$-indexed family of sequences of continuous functions, the identification of these spaces follows as well.
\end{proof}

For $A\subseteq Y^X$ and $K\subseteq X$, we write $A|_K\coloneqq\{f|_K:f\in A\}$ for the set of restrictions to $K$ of the functions in $A$.
Similarly, for $A\subseteq(\mathbb{R}^{\mathcal{P}})^X$ and $P\in\mathcal{P}$, we write $\pi_P\circ A\coloneqq\{\pi_P\circ f:f\in A\}$, where $\pi_P\colon\mathbb{R}^{\mathcal{P}}\to\mathbb{R}$ is the projection onto the coordinate~$P$.
The following theorem extends Theorem~\ref{thm:BFT2} to $\mathbb{R}^{\mathcal{P}}$-valued functions.

\begin{thm}\label{thm:Generalized-BFT}
  Assume that $\mathcal{P}$ is countable, $X$ is a Polish space, and $A\subseteq\allowbreak \Cp(X,\mathbb{R}^\mathcal{P})$ is pointwise bounded in the sense that the set $\pi_P\circ A\subseteq\Cp(X)$ is pointwise bounded for every $P\in\mathcal{P}$.
  The following are equivalent for every compact $K\subseteq X$:

  \begin{enumerate}[(i)]
  \item
    the set $A|_K$ is relatively compact in $\mathrm{B}_1(K,\mathbb{R}^\mathcal{P})$;
  \item
    the set $\pi_P\circ A|_K$ has NIP for every $P\in\mathcal{P}$.
  \end{enumerate}
\end{thm}
\begin{proof}
  By Proposition~\ref{prop:Baire-identifications} and Tychonov’s theorem, condition~(i) holds if and only if $\pi_P\circ A|_K$ is relatively compact in $\mathrm{B}_1(K)$ for every $P\in\mathcal{P}$.
  Since $K$ is compact metrizable, Theorem~\ref{thm:BFT2} shows that this is equivalent to~(ii).
\end{proof}

In the applications, the NIP is verified on a set of states, whereas compactness is available only on the closure of that set.
The following lemma shows that the NIP passes from a set to its closure.

\begin{lem}\label{lem:NIP-and-closure}
  Assume that $X$ is Hausdorff and that $A\subseteq \Cp(X)$.
  The following are equivalent for every $L\subseteq X$:
  \begin{enumerate}[(i)]
  \item
    the set $A|_L$ satisfies the NIP;
  \item
    the set $A|_{\overline{L}}$ satisfies the NIP.
  \end{enumerate}
\end{lem}
\begin{proof}
  The implication (ii)$\Rightarrow$(i) is trivial.
  Conversely, suppose that a sequence $(f_n)$ in $A$ and reals $a<b$ witness the IP on $\overline{L}$, and choose $a'$ and $b'$ with $a<a'<b'<b$.
  For finite disjoint $I,J\subseteq\mathbb{N}$, the open set $\bigcap_{n\in I}f_n^{-1}(-\infty,a')\cap\bigcap_{n\in J}f_n^{-1}(b',\infty)$ contains a point of $\overline{L}$, namely a point at which $f_n\le a$ for $n\in I$ and $f_n\ge b$ for $n\in J$;
  hence it meets~$L$.
  Therefore the sequence $(f_n)$ and the reals $a'<b'$ witness the IP on~$L$.
\end{proof}

\subsection{Five canonical Rosenthal compacta}
\label{sec:five-compacta}
Here we fix the notation for five Rosenthal compacta, all built on the Cantor set, that are the canonical members of the classes of Todorčević’s trichotomy (see Theorem~\ref{thm:Todorcevic-trichotomy} below).
By the trichotomy, the space of deep computations of each non-metrizable $\NIP$ example of Section~\ref{sec:iterative} contains a copy of one of these compacta.
The five compacta are distinguished from one another by three topological properties.
Recall that a topological space is \emph{separable} if it has a countable dense subset and \emph{hereditarily separable} if each of its subspaces is separable.
Finally, a topological space is \emph{first countable} if each of its points has a countable base of neighborhoods.

Let $2^{\mathbb{N}}$ carry the product topology and the lexicographic order~$<$, whose order topology is the product topology.
The set of finite binary words is denoted $2^{<\mathbb{N}}$, and the length of a word $t\in2^{<\mathbb{N}}$ is denoted~$|t|$.
For $a\in2^{\mathbb{N}}$ and $n\in\mathbb{N}$, we write $a\restr n$ for the word formed by the first $n$ bits of~$a$.
Two words are \emph{comparable} if one extends the other, and an \emph{antichain} is a set of pairwise incomparable words.
For a finite word $t\in2^{<\mathbb{N}}$, let $t^\frown0^\infty$ and $t^\frown1^\infty$ be the sequences that begin with $t$ and continue with all $0$’s and with all $1$’s, respectively.
The set of sequences extending~$t$ is then the clopen interval $[t^\frown0^\infty,t^\frown1^\infty]$.
Copies of the binary tree within itself index the dense subset of the last of the five compacta below, and the families of Section~\ref{sec:heptachotomy} are compared along these copies.
These copies are the dyadic subtrees: a \emph{dyadic subtree} of $2^{<\mathbb{N}}$ is a subset isomorphic to $2^{<\mathbb{N}}$ under the order of extension.
We write a dyadic subtree as $\{w_t:t\in2^{<\mathbb{N}}\}$, indexed through the canonical isomorphism, which preserves extension and the lexicographic order.
A dyadic subtree is called \emph{regular} if each of its levels is contained in a level of~$2^{<\mathbb{N}}$.
A sequence $(u_k)$ of words \emph{converges} to $a\in2^{\mathbb{N}}$ if, for every $n\in\mathbb{N}$, all but finitely many of the words $u_k$ extend $a\restr n$.

Double brackets denote truth values: $\Brack{\dots}$ is~$1$ if the statement inside the brackets is true and~$0$ if it is false.
We introduce the following functions $2^{\mathbb{N}}\to\{0,1\}$:
\begin{itemize}
\item
  the \emph{prefix test} $v_t\coloneqq \Brack{\cdot{\restriction}|t| = t}$, the characteristic function of the set of sequences extending a fixed word~$t\in2^{<\mathbb{N}}$;
\item
  the \emph{equality test} $\delta_a\coloneqq\Brack{\cdot = a}$, the characteristic function of the singleton $\{a\}$ for fixed~$a\in2^{\mathbb{N}}$;
\item
  the \emph{threshold tests} $\chi_{[a}\coloneqq\Brack{a\le \cdot}$ and $\chi_{(a}\coloneqq\Brack{a<\cdot}$, the characteristic functions of the final intervals $\{x:a\le x\}$ (closed) and~$\{x:a<x\}$ (open), for fixed~$a\in2^{\mathbb{N}}$;
  the special case $\mathbf{0}\coloneqq \chi_{(1^\infty}$ is the identically zero function.
\end{itemize}
Prefix tests are continuous, and so are the threshold tests $\chi_{[t^\frown0^\infty}$ and $\chi_{(t^\frown1^\infty}$, since each of these functions depends only on the first $|t|$ bits of its argument.
By contrast, equality tests and threshold tests at all other points are discontinuous functions of Baire class~$1$.
All closures below are taken in the topology of pointwise convergence.

\begin{enumerate}
\item
  The \emph{Alexandroff compactification}
  $A(2^{\mathbb{N}})=\{\delta_a:a\in2^{\mathbb{N}}\}\cup\{\mathbf{0}\}$.
  Each $\delta_a$ is isolated, and every neighborhood of $\mathbf{0}$ contains all but finitely many of the functions $\delta_a$.
  Thus the Alexandroff compactification is neither separable nor first countable at~$\mathbf{0}$.
\item
  The \emph{extended Alexandroff compactification}
  $\hat{A}(2^{\mathbb{N}})=\overline{\{v_t:t\in2^{<\mathbb{N}}\}}=\{v_t:t\in2^{<\mathbb{N}}\}\cup A(2^{\mathbb{N}})$.
  As $n\to\infty$, the prefix tests $v_{a\restr n}$ converge to $\delta_a$;
  along an antichain they converge to $\mathbf{0}$.
  Hence the extended Alexandroff compactification is separable but not first countable.
\item
  The \emph{split Cantor set}
  $S(2^{\mathbb{N}})=\{\chi_{[a},\chi_{(a}:a\in2^{\mathbb{N}}\}$.
  The continuous thresholds $\chi_{[t^\frown0^\infty}$, together with the isolated point $\mathbf{0}=\chi_{(1^\infty}$, form a countable dense subset:
  thresholds approaching $a$ from below converge to $\chi_{[a}$, and from above to~$\chi_{(a}$.
  It follows that the split Cantor set is separable, first countable and hereditarily separable, but not metrizable.
\item
  The \emph{Alexandroff duplicate}
  $D(2^{\mathbb{N}})=\{(\mathbf{0},a),(\delta_a,a):a\in2^{\mathbb{N}}\}\subseteq A(2^{\mathbb{N}})\times2^{\mathbb{N}}$.
  The points $(\delta_a,a)$ are isolated, whereas every neighborhood of $(\mathbf{0},a)$ contains, for all $b$ in some neighborhood of~$a$, the point $(\mathbf{0},b)$ and, with finitely many exceptions, the point $(\delta_b,b)$.
  Consequently, the Alexandroff duplicate is first countable but not separable;
  it is contained in the separable compactum $\hat{D}(2^{\mathbb{N}})=\overline{\{(v_t,t^\frown0^\infty):t\in2^{<\mathbb{N}}\}}=D(2^{\mathbb{N}})\cup\{(v_t,t^\frown0^\infty):t\in2^{<\mathbb{N}}\}$.
\item
  The \emph{extended duplicate of the split Cantor set}.
  Following Argyros, Dodos and Kanellopoulos~\cite[4.3.7]{argyros2008rosenthal}, we define it along a regular dyadic subtree $\{s_t:t\in2^{<\mathbb{N}}\}$ of $2^{<\mathbb{N}}$, fixed once and for all.
  This subtree is required to satisfy $s_t^\frown0^\infty\neq s_{t'}^\frown0^\infty$ and $s_t^\frown1^\infty\neq s_{t'}^\frown1^\infty$ for all distinct $t,t'\in2^{<\mathbb{N}}$, and the extended duplicate is then
  \[
    \hat{D}(S(2^{\mathbb{N}}))=\overline{\{(v_{s_t},\chi_{[s_t^\frown0^\infty}):t\in2^{<\mathbb{N}}\}}\subseteq\hat{A}(2^{\mathbb{N}})\times S(2^{\mathbb{N}}).
  \]
  Call a sequence $a\in2^{\mathbb{N}}$ a \emph{branch} of the subtree if each of its initial segments $a\restr n$ is extended by some word~$s_t$;
  the branches form a perfect subset of~$2^{\mathbb{N}}$.
  The compactum $\hat{D}(S(2^{\mathbb{N}}))$ consists of the isolated points $(v_{s_t},\chi_{[s_t^\frown0^\infty})$ and, for each branch~$a$, of the point $(\delta_a,\chi_{[a})$ and of at most two points with first coordinate~$\mathbf{0}$~\cite[4.3.7]{argyros2008rosenthal}.
  The point $(\mathbf{0},\chi_{[a})$ belongs to the compactum exactly when some sequence $(u_k)$ of pairwise incomparable words of the subtree converges to~$a$ with $u_k^\frown0^\infty$ increasing in~$k$.
  Likewise, the point $(\mathbf{0},\chi_{(a})$ belongs to the compactum exactly when some such sequence converges to~$a$ with $u_k^\frown0^\infty$ decreasing in~$k$.
  The points $(\delta_a,\chi_{[a})$ are reached along the branches.
  Indeed, if $a$ is a branch and $u_1,u_2,\dots$ are the words of the subtree that are initial segments of~$a$, listed by increasing length, then the sequences $u_k^\frown0^\infty$ increase strictly to~$a$.
  The points $(v_{u_k},\chi_{[u_k^\frown0^\infty})$ then converge to $(\delta_a,\chi_{[a})$ as $k\to\infty$.
  As the closure of a countable set, the extended duplicate of the split Cantor set is separable.
  This compactum is also first countable, since it embeds into the Helly space of nondecreasing functions from $[0,1]$ into itself, which is first countable~\cite[4.3.7]{argyros2008rosenthal}.
  For each branch~$a$, the point $(\delta_a,\chi_{[a})$ is the only point $(\delta_b,\chi_{[b})$, with $b$ a branch, whose first coordinate takes the value~$1$ at~$a$.
  The points $(\delta_a,\chi_{[a})$, with $a$ a branch, therefore form a discrete subspace, which has the cardinality of the continuum because the branches form a perfect set.
  Hence the extended duplicate of the split Cantor set is not hereditarily separable.
\end{enumerate}

Figure~\ref{fig:five-compacta} shows the five compacta as patterns of convergence of words of the binary tree~$2^{<\mathbb{N}}$.

\begin{figure}[htbp]
  \centering
  \definecolor{cpBerry}{HTML}{9C1C6E}%
  \definecolor{cpPink}{HTML}{F25CA8}%
  \definecolor{cpMain}{HTML}{2BA9C9}%
  \definecolor{cpMainD}{HTML}{1C6F86}%
  \definecolor{cpLav}{HTML}{A58BFF}%
  \definecolor{cpLavD}{HTML}{5B3FC4}%
  \def\cpW{3.5}\def\cpH{0.40}%
  \def\cpN#1#2{({((#2)+0.5)*\cpW/(2^(#1))-\cpW/2},{(#1)*\cpH})}%
  \def\cpTree{%
    \foreach \l in {0,1,2,3}{%
      \pgfmathtruncatemacro{\cpm}{2^\l-1}%
      \pgfmathtruncatemacro{\cpL}{\l+1}%
      \foreach \i in {0,...,\cpm}{%
        \pgfmathtruncatemacro{\cpa}{2*\i}%
        \pgfmathtruncatemacro{\cpb}{2*\i+1}%
        \draw[black!22,line width=0.4pt] \cpN{\l}{\i} -- \cpN{\cpL}{\cpa};
        \draw[black!22,line width=0.4pt] \cpN{\l}{\i} -- \cpN{\cpL}{\cpb};
      }}%
    \foreach \l in {0,...,4}{%
      \pgfmathtruncatemacro{\cpm}{2^\l-1}%
      \foreach \i in {0,...,\cpm}{\fill[black!30] \cpN{\l}{\i} circle (0.75pt);}}}%
  \def\cpBranch#1{%
    \draw[cpBerry,line width=1.1pt] \cpN{0}{0} -- \cpN{1}{0} -- \cpN{2}{1} -- \cpN{3}{3} -- \cpN{4}{6};
    \draw[cpBerry,line width=0.8pt,dashed,-latex] \cpN{4}{6} -- ++(0.06,0.80);
    \foreach \l/\i in {0/0,1/0,2/1,3/3,4/6}{\fill[cpPink,draw=cpBerry,line width=0.4pt] \cpN{\l}{\i} circle (1.7pt);}
    \path \cpN{4}{6} ++(0.06,0.82) node[cpBerry,font=\scriptsize,anchor=south] {#1};}%
  \def\cpBelow#1#2{%
    \draw[cpMainD,line width=0.4pt,densely dotted] \cpN{2}{0} -- \cpN{3}{2};
    \draw[cpMainD,line width=0.6pt,dashed,-latex] \cpN{3}{2} -- (#2);
    \foreach \l/\i in {2/0,3/2}{\fill[cpMain,draw=cpMainD,line width=0.4pt] \cpN{\l}{\i} circle (1.7pt);}
    \node[cpMainD,font=\scriptsize,anchor=south] at (#2) {#1};}%
  \def\cpAbove#1#2{%
    \draw[cpLavD,line width=0.4pt,densely dotted] \cpN{1}{1} -- \cpN{4}{7};
    \draw[cpLavD,line width=0.6pt,dashed,-latex] \cpN{4}{7} -- (#2);
    \foreach \l/\i in {1/1,4/7}{\fill[cpLav,draw=cpLavD,line width=0.4pt] \cpN{\l}{\i} circle (1.7pt);}
    \node[cpLavD,font=\scriptsize,anchor=south] at (#2) {#1};}%
  \def\cpTitle#1{\node[font=\small,anchor=north] at (0,-0.22) {#1};}%
  \begin{tikzpicture}
    \begin{scope}[shift={(0,0)}]
      \cpTree
      \draw[cpBerry,line width=0.9pt] \cpN{0}{0} -- \cpN{1}{0} -- \cpN{2}{1} -- \cpN{3}{3} -- \cpN{4}{6};
      \draw[cpBerry,line width=0.9pt] \cpN{0}{0} -- \cpN{1}{1} -- \cpN{2}{2} -- \cpN{3}{4} -- \cpN{4}{9};
      \draw[cpBerry,line width=0.9pt] \cpN{0}{0} -- \cpN{1}{1} -- \cpN{2}{3} -- \cpN{3}{7} -- \cpN{4}{14};
      \foreach \i/\cpc in {6/cpBerry,9/cpBerry,14/cpBerry}{%
        \draw[\cpc,line width=0.6pt,dashed] \cpN{4}{\i} -- ++(0,0.3);
        \draw[\cpc,fill=white,line width=0.6pt] \cpN{4}{\i} ++(0,0.36) circle (1.8pt);
        \draw[\cpc,line width=0.5pt,-latex] \cpN{4}{\i} ++(0,0.44) -- (0.15,2.62);}
      \node[font=\scriptsize,anchor=south] at (0.15,2.62) {$\delta_{a_1},\delta_{a_2},\delta_{a_3},\ldots\to\mathbf{0}$};
      \cpTitle{$A(2^{\mathbb{N}})$}
    \end{scope}
    \begin{scope}[shift={(4.35,0)}]
      \cpTree
      \cpBelow{$\mathbf{0}$}{-1.45,2.2}
      \cpAbove{$\mathbf{0}$}{1.35,2.2}
      \cpBranch{$\delta_a$}
      \cpTitle{$\hat{A}(2^{\mathbb{N}})$}
    \end{scope}
    \begin{scope}[shift={(8.7,0)}]
      \cpTree
      \cpBelow{$\chi_{[a}$}{-1.45,2.2}
      \cpAbove{$\chi_{(a}$}{1.35,2.2}
      \cpBranch{$\chi_{[a}$}
      \cpTitle{$S(2^{\mathbb{N}})$}
    \end{scope}
    \begin{scope}[shift={(2.175,-3.8)}]
      \cpTree
      \cpBelow{$(\mathbf{0},a)$}{-1.45,2.2}
      \cpAbove{$(\mathbf{0},a)$}{1.35,2.2}
      \cpBranch{$(\delta_a,a)$}
      \cpTitle{$\hat{D}(2^{\mathbb{N}})\supseteq D(2^{\mathbb{N}})$}
    \end{scope}
    \begin{scope}[shift={(6.525,-3.8)}]
      \cpTree
      \cpBelow{$(\mathbf{0},\chi_{[a})$}{-1.6,2.2}
      \cpAbove{$(\mathbf{0},\chi_{(a})$}{1.5,2.2}
      \cpBranch{$(\delta_a,\chi_{[a})$}
      \cpTitle{$\hat{D}(S(2^{\mathbb{N}}))$}
    \end{scope}
  \end{tikzpicture}%
  \caption{The five compacta of Section~\ref{sec:five-compacta}, drawn on the binary tree $2^{<\mathbb{N}}$.
    The tree grows upward, and the extension of a word by~$0$ lies to the left of its extension by~$1$.
    In the first tree, which stands for $A(2^{\mathbb{N}})$, the equality tests $\delta_{a_k}$ of distinct sequences $a_k\in2^{\mathbb{N}}$ converge to the zero function~$\mathbf{0}$.
    In the other trees, the thick path runs through the initial segments of a sequence $a\in2^{\mathbb{N}}$, here $a=0110\cdots$.
    The dots on the left of the path mark words~$t$ with $t^\frown0^\infty<a$, and those on the right mark words with $t^\frown0^\infty>a$.
    Both sequences of words converge to~$a$, and the dashed arrows end at the limits.
    A word~$t$ stands for the prefix test~$v_t$ in $\hat A(2^{\mathbb{N}})$ and for the threshold test $\chi_{[t^\frown0^\infty}$ in $S(2^{\mathbb{N}})$.
    In $\hat D(2^{\mathbb{N}})$, which contains the Alexandroff duplicate $D(2^{\mathbb{N}})$, the word~$t$ stands for the pair $(v_t,t^\frown0^\infty)$.
    In $\hat D(S(2^{\mathbb{N}}))$, the tree is the dyadic subtree of words $s_t$ fixed in Section~\ref{sec:five-compacta}, and the node $s_t$ stands for the pair $(v_{s_t},\chi_{[s_t^\frown0^\infty})$.}
  \label{fig:five-compacta}
\end{figure}

The constructions of the Alexandroff compactification and of the Alexandroff duplicate make sense with $2^{\mathbb{N}}$ replaced by any compact metric space~$K$;
the resulting spaces are denoted $A(K)$ and $D(K)$, respectively.
Since $A(K)$ is the one-point compactification of a discrete space of the cardinality of~$K$, it depends only on that cardinality;
for instance, $A([0,1])$ is homeomorphic to~$A(2^{\mathbb{N}})$.

\subsection{NIP in model theory, combinatorics and learning: polynomial versus exponential growth}
\label{sec:NIP-poly-exp}

We now trace the NIP of Definition~\ref{def:NIP-fns} to its counterparts in model theory, combinatorics and statistical learning, and we make precise how it is related to them.
This relation is used in Section~\ref{sec:iterative}, where a bound on a Vapnik--Chervonenkis dimension yields the NIP of Bresenham’s algorithm.
For families of $\{0,1\}$-valued functions (that is, of sets), a combinatorial dividing line closely related to the NIP appeared independently in three fields between 1968 and 1972, some years before Theorem~\ref{thm:BFT}.
In model theory the line arose in the work of Saharon Shelah~\cite{Shelah:1971,Shelah:1990}, and in combinatorics in the work of Norbert Sauer~\cite{Sauer:1972} and, again, of Shelah~\cite{Shelah:1972,Shelah:1990}.
In statistical learning, finally, the same line appeared in the work of Vladimir Vapnik and Alexey Chervonenkis~\cite{Vapnik-Chervonenkis:1968,Vapnik-Chervonenkis:1971,Vapnik-Chervonenkis:1974}.
In each of the three fields, the line separates the “tame” families of sets from the “wild” ones.
The number of traces of a tame family on an $n$-element set (the intersections of the members of the family with the set) is bounded by a polynomial in~$n$.
In a wild family, by contrast, for every $n\in\mathbb{N}$ some $n$-element set has all $2^n$ possible traces.

\begin{description}

\item[In model theory]
  Among the dividing lines of Shelah’s classification theory~\cite{Shelah:1990}, the independence property is detected by counting types over finite sets of parameters.
  Let $T$ be a complete first-order theory, and let $\varphi(x;y)$ be a formula of its language, in which $x$ and $y$ are finite tuples of variables of lengths $|x|$ and~$|y|$.
  Finally, let $B$ be a finite set of $|y|$-tuples of elements of a model $M$ of~$T$.
  A \emph{$\varphi$-type} over $B$ is a maximal set of formulas $\varphi(x;b)$ and $\neg\varphi(x;b)$, with $b\in B$, that is consistent with the theory of $M$ with constants for the elements of~$M$.
  Since $B$ is finite, every $\varphi$-type over $B$ is realized by an $|x|$-tuple $a$ of elements of~$M$, and the type is determined by the set $\{b\in B: M\models\varphi(a;b)\}$.
  Hence the set $S_\varphi(B)$ of $\varphi$-types over $B$ is in bijection with the traces on $B$ of the family of sets $\{b: M\models\varphi(a;b)\}$, with $a$ ranging over $|x|$-tuples of elements of~$M$.
  This bijection gives the bound $|S_\varphi(B)|\le2^{|B|}$.
  The formula $\varphi$ has the \emph{independence property} (IP) if, for every $n\in\mathbb{N}$, some set $B$ of $n$ tuples satisfies $|S_\varphi(B)|=2^n$.
  Otherwise the formula is \emph{NIP}, or \emph{dependent}, and the theory $T$ is NIP if all its formulas are NIP.
  Since $T$ is complete, these notions do not depend on the model~$M$.

  If $\varphi$ is NIP, the Sauer--Shelah lemma below bounds $|S_\varphi(B)|$ by a polynomial in $|B|$.
  For a set $A$ of parameters, we write $S_\varphi(A)$ for the set of $\varphi$-types over the set of $|y|$-tuples of elements of~$A$.
  If $\varphi$ has IP, then for every $n$ some set $A$ of at most $n|y|$ parameters satisfies $|S_\varphi(A)|\ge2^n$~\cite[Section~2.3.4]{simon2015nip}.
  Consequently, $\varphi$ is NIP if and only if $|S_\varphi(A)|$ is bounded by a polynomial in $|A|$ as $A$ ranges over the finite sets of parameters.
  Over infinite sets of parameters, counting $\varphi$-types does not separate NIP from IP without further set-theoretic hypotheses.
  Indeed, every unstable formula, NIP or not, already has $2^\kappa$ $\varphi$-types over some set of $\kappa$ parameters whenever $\kappa$ is infinite and $2^\kappa=\kappa^+$~\cite[Section~2.3.4]{simon2015nip}.

\item[In combinatorics]
  Sauer~\cite{Sauer:1972} and Shelah~\cite{Shelah:1972} proved independently the following dichotomy, now known as the \emph{Sauer--Shelah lemma}.
  Let $\mathscr{F}=\{C_i:i\in I\}$ be a family of subsets of a set $\Omega$, indexed by a set~$I$.
  Then either
  for every $n\in\mathbb{N}$ there is a set $B\subseteq \Omega$ with $|B|=n$ such that $|\{C_i\cap B: i\in I\}|=2^n$,
  or
  there exists $N\in\mathbb{N}$ such that, for every $B\subseteq \Omega$ with $|B|\ge N$,
  \[
    |\{C_i\cap B: i\in I\}| \le \sum_{j=0}^{N-1} \binom{|B|}{j}.
  \]
  In the first case $\mathscr{F}$ is said to have “exponential complexity”, and in the second case “polynomial complexity”;
  every family of subsets of a \emph{finite} set $\Omega$ has polynomial complexity.
  The name “polynomial complexity” is justified by the following estimate, valid for fixed $N>0$ and for $|B| = m > 0$, which bounds the number of traces by a polynomial in~$m$:
  \begin{equation*}
  \sum_{j=0}^{N-1} \binom{m}{j}
  \le \sum_{j=0}^{N-1}\frac{m^{j}}{j!}
  \le \left( \sum_{j=0}^{N-1}\frac{1}{j!} \right)\cdot m^{N-1}
  < e\cdot m^{N-1}.
  \end{equation*}

\item[In machine learning]
  The Sauer--Shelah lemma is also the dichotomy that Vapnik and Chervonenkis had announced in 1968~\cite{Vapnik-Chervonenkis:1968}, in their study of uniform convergence in statistics.
  Their proof of the dichotomy appeared in 1971~\cite{Vapnik-Chervonenkis:1971}, a year before the papers of Sauer and Shelah (see also~\cite{Vapnik-Chervonenkis:1974}).
  The largest size of a set $B\subseteq \Omega$ with $2^{|B|}$ traces $C_i\cap B$ is called the \emph{VC dimension} of~$\mathscr{F}$.
  Equivalently, the VC dimension is the integer $N-1$ for the least integer $N$ as in the Sauer--Shelah lemma.
  If no such integer $N$ exists, the VC dimension of~$\mathscr{F}$ is infinite.
  Through the Vapnik--Chervonenkis theorem, the Sauer--Shelah lemma bounds the number of data points (the sample size) needed to “learn” a concept class~$\mathscr{F}$ of finite VC dimension in a statistical sense.
  Here “learning” is understood up to a preset error and with a preset confidence.
  The Fundamental Theorem of Statistical Learning concerns PAC-learnability, where PAC stands for “Probably Approximately Correct”~\cite{Valiant:1984}.
  This theorem states that a class of $\{0,1\}$-valued functions satisfying a mild measurability condition is PAC-learnable if and only if its VC dimension is finite~\cite{Blumer-Ehrenfeucht-Haussler-Warmuth:1989};
  see also~\cite[Theorem~6.7]{Shalev-Shwartz-Ben-David:2014}.

\end{description}

For first-order formulas the three dichotomies coincide.
To see this, fix a formula $\varphi(x;y)$ and consider the family of sets $\{b: M\models\varphi(a;b)\}$, with $a$ ranging over the $|x|$-tuples of elements of a model $M$ of~$T$.
By the description of $\varphi$-types given above, the formula $\varphi$ has IP exactly when this family has infinite VC dimension.
For an NIP formula, moreover, the Sauer--Shelah lemma gives the polynomial bound on the number of $\varphi$-types.
For an arbitrary family of sets, finiteness of the VC dimension implies that the characteristic functions of its members have the NIP of Definition~\ref{def:NIP-fns}.
The NIP does not, however, imply finite VC dimension: the finite subsets of $\mathbb{N}$ form a family of infinite VC dimension whose characteristic functions have the NIP.

\subsection*{Historical remarks}
The Baire hierarchy of functions was introduced by René-Louis Baire in his 1899 doctoral thesis, \emph{Sur les fonctions de variables réelles}.
The theory of spaces of continuous functions with the topology of pointwise convergence, $\Cp$-theory, was pioneered by A.~V.~Arhangel’skiĭ and his students in the 1970s and 1980s~\cite{Arkhangelskii:1992}.
This theory has since found numerous applications in model theory and in functional analysis;
Tkachuk’s problem book~\cite{tkachuk2011cp} develops the subject systematically, through problems with complete solutions.

\section{Compositional Computation Structures: \\
  A structural approach to arbitrary-precision computations}
\label{sec:CCS}

The theory of function spaces assembled in Section~\ref{sec:prelim} is brought to bear on computation through the framework of compositional computation structures, to which the present section is devoted.
Within this framework, deep computations become points of a function space.
A computation on real-valued data is observed through its features, each known to finite precision.
A deep computation is a limit of computations along an ultrafilter in the topology of this observation, namely the topology of pointwise convergence.
Following the first paper of this series~\cite{alva2024approximability}, we first describe states by their types, the tuples of values of all their features.
Computations act on these states and are composed;
deep computations are the pointwise limits of computations along ultrafilters, taken in the space of types.
Finally, we recall the stable sets of computations, studied in that paper, and close the section with two examples that motivate the classification of Section~\ref{sec:Classification}.

\subsection{Computation States Structures}
\label{sec:CSS}

The first ingredient of the framework, which we introduce here, is the set of states on which computations act, together with the real-valued features through which the states are observed.
A \emph{computation states structure} is a pair $(L,\mathcal{P})$ consisting of a set~$L$ of \emph{states} and a collection~$\mathcal{P}$ of real-valued functions on~$L$ that we call \emph{predicates} or \emph{features}.
Intuitively, the set $L$ consists of the states of a computation, and the predicates are the primitives that are given and accepted as computable.
For $P\in\mathcal{P}$ and $v\in L$, the number $P(v)$ is the value of the (observable) feature~$P$ at the state~$v$.
Maps $f\colon L\to L$ are called \emph{transitions} of $(L,\mathcal{P})$;
the computations of Definition~\ref{def:CCS} below are transitions.

The \emph{type} of a state $v\in L$ is the indexed family%
\footnote{As in Proposition~\ref{prop:Baire-identifications}, where $\mathcal{P}$ is merely a set of symbols, we regard each predicate as named by a formal label~$P$, a \emph{predicate symbol}.
  The meaning of this symbol in the computation states structure is the real function $P(\cdot)\colon v\mapsto P(v)$ on~$L$.
  The symbol $P$ thus names an observable feature, whose value at a state $v\in L$ is $P(v)$.
  Whenever $\mathcal{P}$ is used as an index set, as in the product $\mathbb{R}^{\mathcal{P}}$ below, it is regarded as a set of symbols.
  The distinction between a symbol and its meaning is usually suppressed in mathematics, but the model-theoretic framework requires awareness of it.
  When necessary, we write $P(\cdot)$ for the function on~$L$ and $P$ for the symbol;
  when no confusion can arise, we also write $P$ for the function.}
\[
  \operatorname{tp}(v)= (P(v))_{P\in \mathcal{P}}\in\mathbb{R}^\mathcal{P}
\]
of its features.

We assume that each state $v\in L$ is determined by its type $\operatorname{tp}(v)$.
Under this assumption we may identify $L$ with its image $\tp[L]\subseteq \mathbb{R}^{\mathcal{P}}$ under the type map and endow $L$ with the subspace topology.
This topology is the topology of pointwise convergence induced by the product topology of~$\mathbb{R}^\mathcal{P}$;
in particular, for each $P\in \mathcal{P}$, the map $v\mapsto P(v)$ is continuous on~$L$.

The principal examples are the state space $L=\mathbb{R}^\mathbb{N}$ with the predicates $P_i(v) = v_i$ for $i\in\mathbb{N}$ and, for a positive integer~$n$, the state space $L=\mathbb{R}^n$ with the predicates $P_i(v)=v_i$ for $0\le i<n$.
Their subspaces are the state spaces of primary interest in this paper.
Indeed, through the identification of each state with its type, every computation states structure with countably many predicates becomes a subspace of~$\mathbb{R}^\mathbb{N}$ or of some~$\mathbb{R}^n$ once its predicates are enumerated.

\begin{defn}
  Let $(L,\mathcal{P})$ be a computation states structure.
  The type $\tp(v)\in \mathbb{R}^{\mathcal{P}}$ of a state $v\in L$ is called a \emph{realized} (state) type.
  The closure of the set of realized types in~$\mathbb{R}^\mathcal{P}$, with the topology of pointwise convergence, is called the \emph{space of (state) types} of~$(L,\mathcal{P})$ and is denoted~$\mathcal{L}$.
  Its elements $\mathbf{v}\in \mathcal{L}$ are called \emph{state types}, and those in $\mathcal{L}\setminus L$ are called \emph{unrealized} state types.
\end{defn}

State types formalize the intuitive notion of a limit state, and the predicates extend to them continuously.
Indeed, for each $P\in \mathcal{P}$, the restriction to~$\mathcal{L}$ of the $P$-th coordinate map $\mathbb{R}^{\mathcal{P}}\to \mathbb{R}$ given by $\mathbf{v} = (\mathbf{v}_Q)_{Q\in \mathcal{P}}\mapsto \mathbf{v}_P$ is a continuous extension of the predicate $P(\cdot)$ on~$L$.
Abusing notation, we still denote this extension by~$P(\cdot)$, or just~$P$.

\begin{rem}\label{rem:poor-mans-ultraprod}
  One may regard unrealized state types as “hyper-states”, and those among them that are obtained as limits of sequences (or as ultralimits) of realized state types as “deep states”.
  For readers versed in model theory we note that, at least for realizing state types with parameters in~$L$, the type space $(\mathcal{L}, \mathcal{P})$ serves as a “poor man’s” ultraproduct of~$(L, \mathcal{P})$.
  Indeed, by the very definition of~$\mathcal{L}$, every type $\mathbf{v}$ is realized in this type space.
\end{rem}

\subsubsection{Bounds on precision and magnitude. Shards}
\label{sec:shards}

The analysis of this paper takes place on shards, the sets of states on which each feature obeys a prescribed bound;
before defining them, we explain why actual computation imposes such bounds.
Computation states structures are designed to model actual computation, in which operations on real-valued quantities can be carried out only approximately, to a preset degree of precision.
A bound on the magnitude of the quantities involved is no less important in practice, for no machine can store, let alone compute with, arbitrarily large numbers.
An actual computation therefore fixes at the outset bounds for its errors and for the absolute magnitudes of the quantities it handles.
These bounds are fixed even when the collection~$\mathcal{P}$ of predicates is infinite, possibly uncountable.
Moreover, different parts of a computation may require different error bounds and may deal with quantities of very different sizes.
Focusing on magnitudes, we consider sets of computations for which a bound for each predicate is fixed in advance.
We are thereby led to the notion of \emph{shard}, introduced in the first paper of this series~\cite{alva2024approximability}.

\begin{defn}
  Let $(L, \mathcal{P})$ be a computation states structure, with $L$ identified with the set $\tp[L]\subseteq \mathbb{R}^{\mathcal{P}}$ of its realized types.
  A \emph{sizer} is a family $r_{\bullet}=(r_P)_{P\in\mathcal{P}}$ of positive real numbers indexed by~$\mathcal{P}$.
  For a sizer $r_\bullet$, let $\mathbb{R}[\rb] \coloneqq \prod_{P\in\mathcal{P}}[-r_P,r_P]$, a compact space, and define the $r_\bullet$-\emph{shard} of the state space~$L$ by
  \[
    L[r_\bullet] = L\cap \mathbb{R}[\rb].
  \]

  The \emph{$r_{\bullet}$-type shard} is the closure $\mathcal{L}[r_\bullet]=\overline{L[r_\bullet]}$, a closed and hence compact subset of $\mathbb{R}[\rb]$.

  The union of all type shards, as the sizer $r_{\bullet}$ varies, is denoted $\mathcal{L}_{\sh}$ and called the space of \emph{shard-supported types}.
\end{defn}

We have $\mathcal{L}_{\sh}\subseteq \mathcal{L}$, and the inclusion may be proper;
equality holds in the important special case when $\mathcal{P}$ is countable (see~\cite{alva2024approximability}).

\subsection{Compositional Computation Structures}

In this subsection we recall the compositional computation structures of the first paper of this series~\cite{alva2024approximability}, within which deep computations are defined, beginning with the ingredients of their definition.
Computations transform states, and two computations performed one after the other form a computation.
The computations of a system therefore form a semigroup under composition, which we take to contain the identity, the computation that does nothing.
Deep computations are the limits of computations along ultrafilters in the topology of pointwise convergence, taken in the space of types.
Since an actual computation keeps the magnitudes that it handles within bounds, the shards of Section~\ref{sec:shards} enter the definition through the notion of confinement.
The following definition records these ingredients.

\begin{defn}
  \label{def:CCS}
  A \emph{compositional computation structure} (CCS) is a triple $(L,\mathcal{P},\Gamma)$, where
  \begin{itemize}
  \item
    the pair $(L,\mathcal{P})$ is a computation states structure, and
  \item
    the set $\Gamma\subseteq L^L$ is a semigroup under composition that contains the identity map $\text{id}$ of~$L$, that is, a monoid of transformations of~$L$.
  \end{itemize}
  The elements of $\Gamma$ are called the \emph{computations} of $(L,\mathcal{P},\Gamma)$.
  The state space $L$ is identified with a subset of the type space~$\mathcal{L} \subseteq \mathbb{R}^{\mathcal{P}}$ of~$(L, \mathcal{P})$.

  We topologize $\Gamma$ as a subset of the product $L^L$.
  In this product, as usual, the exponent~$L$ serves merely as an index set, whereas the base~$L\subseteq \mathbb{R}^{\mathcal{P}}$ carries the topology of pointwise convergence.
  The semigroup $\Gamma$ is thus a subset of the product~$\mathcal{L}^L \subseteq (\mathbb{R}^{\mathcal{P}})^L$, and its closure $\overline{\Gamma}$ is contained in~$\mathcal{L}^L$.
  The elements $\xi\in\overline{\Gamma}$ are called (real-valued) \emph{deep computations}, and the elements of the closure $\overline{\Delta}$ of a set $\Delta\subseteq\Gamma$ are called the deep computations of~$\Delta$.%
  \footnote{In the first paper of this series~\cite{alva2024approximability}, the deep computations of~$\Delta$ are the accumulation points of~$\Delta$ that satisfy a confinement condition.
    Counting the members of~$\Delta$ among them affects none of the properties studied in this paper.
    Indeed, under the Extendibility Axiom of Section~\ref{sec:computability}, the members of~$\Delta$ are continuous on the type shards, and those that are not accumulation points of~$\Delta$ are isolated points of~$\overline{\Delta}$.}

  A collection $\mathscr{R}$ of sizers is \emph{exhaustive} if $L = \bigcup_{\rb\in \mathscr{R}}L[\rb]$.
  A computation $\gamma\in\Gamma$ is \emph{$\mathscr{R}$-confined} if it maps $L[r_\bullet]$ into itself for every $r_\bullet\in \mathscr{R}$, and a subset $\Delta\subseteq\Gamma$ is \emph{$\mathscr{R}$-confined} if each $\gamma\in \Delta$ is.
\end{defn}

Unlike its subspace $L^L$, which is a semigroup under composition, the space $\mathcal{L}^L$ has no natural semigroup structure;
deep computations are therefore not composable in general.

\begin{prop}\label{prop:confined-compact}
  If $\Delta\subseteq\Gamma$ is confined by an exhaustive collection of sizers, then its closure $\overline{\Delta}$ is a compact subset of~$\Lshard^L$.
\end{prop}
\begin{proof}
  Let $\mathscr{R}$ be an exhaustive collection of sizers that confines~$\Delta$, and choose, for each $v\in L$, a sizer $\rb^{(v)}\in \mathscr{R}$ with $v\in L[\rb^{(v)}]$.
  Since every $\gamma\in\Delta$ maps $v$ into $L[\rb^{(v)}]$, the set $\Delta$ lies in the compact space $\prod_{v\in L}\mathcal{L}[\rb^{(v)}]\subseteq\Lshard^L$, which therefore contains~$\overline{\Delta}$.
\end{proof}

For a CCS $(L,\mathcal{P},\Gamma)$ we regard the elements of $\Gamma$ as “standard” finitary computations and the elements of $\overline{\Gamma}$, i.e., deep computations, as possibly infinitary limits of standard computations.

\begin{rems}(For readers versed in model theory.)%
  \footnote{The real-valued structures used here are those of Keisler’s model theory for real-valued structures~\cite{Keisler:2023}; see also~\cite[Section~6]{Duenez-Iovino:2017} for an introduction.}
  \begin{enumerate}
  \item
    Compositional computation structures, as defined above, are not literally real-valued structures, because nothing in the structure itself forces $\Gamma$ to be realized as a subset of the product space $L^L$.
    However, as in their original definition~\cite{alva2024approximability}, they may be regarded as real-valued structures $(L,\mathcal{P},\Gamma,\circ,\operatorname{ev})$, where
    \begin{itemize}
    \item
      the pair $(L,\mathcal{P})$ is a computation states structure,
    \item
      the pair $(\Gamma,\circ)$ is a semigroup,
    \item
      the map $\operatorname{ev}: \Gamma\times L\to L: (\gamma,v)\mapsto \operatorname{ev}(\gamma,v)$ is a semigroup action of $\Gamma$ on~$L$
      (i.e., $\operatorname{ev}(\gamma\circ\delta,v) = \operatorname{ev}(\gamma,\operatorname{ev}(\delta,v))$ for $\gamma,\delta\in\Gamma$ and $v\in L$).
    \end{itemize}
    The two definitions are reconciled by identifying each $\gamma\in \Gamma$ with the map $v\mapsto \operatorname{ev}(\gamma,v)$.

  \item
    Regarded as real-valued structures in this sense, the compositional computation structures with a fixed collection~$\mathcal{P}$ of predicate symbols form an elementary class.
    The notion of state type introduced above is essentially that of a quantifier-free type with parameters of sort $L$ only.
    This notion may be enriched by further formulas, for instance with parameters in~$\Gamma$ or with quantifiers, and extended to types of sort~$\Gamma$ (“transform-types”).
    As usual, all types over a structure $\underline{C}=(L,\mathcal{P}, \Gamma, \circ, \operatorname{ev})$
    are realized in some ultrapower~$\underline{C}^{(\mathcal{U})}$ of~$\underline{C}$, constructed as in the references above.%
    \footnote{By contrast, no ultrapower is needed to realize state types in a computation states structure (see Remark~\ref{rem:poor-mans-ultraprod}).}
    In Banach space theory, this approach to spaces of types goes back to Krivine's real-valued languages~\cite{Krivine:1972,Krivine:1974} and to the work of Krivine and Maurey on stable Banach spaces~\cite{Krivine-Maurey:1981}.
  \end{enumerate}
\end{rems}

\subsection{Computability of deep computations and the Extendibility Axiom}
\label{sec:computability}
The rest of the paper depends on two notions, the computability of transforms and the Extendibility Axiom, both of which we recall in this subsection.
The Extendibility Axiom requires every computation to have an extension that is continuous on each type shard.
As for computability, a machine evaluates only sums and products of finitely many observed features, and only to within a preset error on a preset range of magnitudes.
For this reason, computability is defined by uniform approximation on shards by polynomials in the features.
Let $(L,\mathcal{P})$ be a computation states structure.
A map $f\colon L\to\mathcal{L}$ from its state space to its space of types is called a \emph{transform}.
In the two definitions that follow, and in the characterization of computability after them, the collection~$\mathcal{P}$ of predicates is assumed to be \emph{countable}.

\begin{enumerate}
\item
  A \emph{computable predicate} is a function from~$L$ to~$\mathbb{R}$ whose restriction to each shard is a uniform limit of polynomials, with real coefficients, in the features~$P(\cdot)$.
\item
  A transform $f\colon L\to \mathcal{L}$ is \emph{computable} if, for each $Q\in\mathcal{P}$, the output feature $Q{\circ}f\colon L\to\mathbb{R}$ is a computable predicate.
\end{enumerate}

Alva, Dueñez, Iovino and Walton~\cite{alva2024approximability} showed that a transform $f\colon L\to\mathcal{L}$ is computable if and only if it extends to a continuous map $\tilde{f}\colon \mathcal{L} \to \mathcal{L}$.
This characterization motivates the following definition, in which continuity is required on each type shard.

A computation $\gamma\in \Gamma$ is \emph{extendable} if it has an extension $\tilde\gamma\colon\mathcal{L}_{\sh}\to\mathcal{L}_{\sh}$ such that,
for every sizer $r_{\bullet}$, there is a sizer $s_\bullet$ for which the restriction $\tilde\gamma|_{\mathcal{L}[r_\bullet]}\colon\mathcal{L}[r_\bullet]\to\mathcal{L}[s_\bullet]$ is continuous.
Such an extension is unique, since it is continuous on each type shard $\mathcal{L}[\rb]$, in which $L[\rb]$ is dense.
For a set $\Delta\subseteq\Gamma$ of extendable computations, we write $\tilde{\Delta}$ for the set $\{\tilde\gamma\colon\gamma\in\Delta\}$.

\begin{AxExt}
  \label{ax:extendibility-axiom}
  Every computation $\gamma\in\Gamma$ is extendable.
\end{AxExt}

The Extendibility Axiom is discussed at length in the first paper of this series~\cite{alva2024approximability}.
\emph{From now on, all compositional computation structures are assumed to satisfy the Extendibility Axiom.}

The following lemma extends deep computations to the type shards;
it is used in Theorems~\ref{thm:nip-baire1def} and~\ref{thm:NIP-Monte-Carlo}.

\begin{lem}\label{lem:extension}
  Let $\mathscr{R}$ be an exhaustive collection of sizers and let $\Delta\subseteq\Gamma$ be $\mathscr{R}$-confined.
  Every deep computation $f\in\overline{\Delta}$ has an extension $\tilde f:\Lshard\to\Lshard$ that, for every $\rb\in \mathscr{R}$, maps $\mathcal{L}[\rb]$ into itself and satisfies $\tilde f|_{\mathcal{L}[\rb]}\in\overline{\tilde\Delta|_{\mathcal{L}[\rb]}}$.
\end{lem}
\begin{proof}
  Write $f=\lim_{i,\mathcal{U}}\gamma_i$, where $(\gamma_i)_{i\in I}$ is a family in~$\Delta$ and $\mathcal{U}$ is an ultrafilter on the index set~$I$, and let $\rb\in \mathscr{R}$ be given.
  Each $\tilde\gamma_i$ maps $L[\rb]$ into itself and is continuous on $\mathcal{L}[\rb]$;
  hence it maps the closure $\mathcal{L}[\rb]$ of $L[\rb]$ into itself.
  Since $\mathcal{L}[\rb]$ is compact, the limit $\tilde f(v)\coloneqq\lim_{i,\mathcal{U}}\tilde\gamma_i(v)$ exists and lies in $\mathcal{L}[\rb]$ for every $v\in\mathcal{L}[\rb]$;
  in other words, the restriction $\tilde f|_{\mathcal{L}[\rb]}$ equals $\lim_{i,\mathcal{U}}\tilde\gamma_i|_{\mathcal{L}[\rb]}$.
  Since $\mathscr{R}$ is exhaustive, this defines $\tilde f$ on a set that contains~$L$, and there $\tilde f$ agrees with~$f$;
  at the remaining points of $\Lshard$ we let $\tilde f$ be the identity.
\end{proof}

Lemma~\ref{lem:extension} connects deep computations with the second notion of limit introduced in the first paper of this series~\cite{alva2024approximability}.
In that paper, this second notion is formulated in the space of maps from the space of types to itself, endowed with the topology of pointwise convergence.
An \emph{ultracomputation} of a set $\Delta\subseteq\Gamma$ of extendable computations is an accumulation point of~$\tilde\Delta$ in that space, subject to a confinement condition.
Thus a deep computation sends states to types, whereas an ultracomputation sends types to types.
Every deep computation is the restriction to~$L$ of an ultracomputation~\cite{alva2024approximability}.
On the type shards $\mathcal{L}[\rb]$ with $\rb\in\mathscr{R}$, the ultralimit $\lim_{i,\mathcal{U}}\tilde\gamma_i$ in the proof of Lemma~\ref{lem:extension} provides such an extension of the deep computation $f=\lim_{i,\mathcal{U}}\gamma_i$ considered there.
The classification of Section~\ref{sec:Classification} is based on the sets $\overline{\tilde\Delta|_{\mathcal{L}[\rb]}}$, whose members are the restrictions to $\mathcal{L}[\rb]$ of such ultralimits of extended computations.
Strictly speaking, the classification therefore concerns the restrictions of ultracomputations to type shards.
For the dynamical compositional computation structures of Section~\ref{sec:dynamical-CCS}, in which every type is realized and every computation is its own extension, deep computations and ultracomputations coincide.

\subsection{Stability}
\label{sec:stable}
The highest level of the hierarchy studied in this paper is occupied by the stable sets of computations, all of whose deep computations are computable when $\mathcal{P}$ is countable.
We recall these sets here, before descending to the levels below in Section~\ref{sec:Classification}.
By the theorem of Alva, Dueñez, Iovino and Walton recalled above, the computable transforms are exactly those that extend to continuous maps on the space of types.
Accordingly, the stable sets are defined by the continuity of their deep computations on every type shard.

\begin{defn}\label{def:stable}
  A set $\Delta\subseteq\Gamma$ of computations is \emph{stable} if, for every sizer $r_\bullet$, every deep computation in $\overline{\tilde\Delta|_{\mathcal{L}[r_\bullet]}}$ is continuous (that is, computable if $\mathcal{P}$ is countable).
\end{defn}

Stable sets of computations, the subject of the first paper of this series~\cite{alva2024approximability}, are characterized there by Grothendieck’s double limit criterion and identified with stability in the sense of model theory.

\subsection{Two motivating examples}
\label{sec:preview}
Before developing the classification, we describe two algorithms whose deep computations differ sharply in character, for the contrast between them indicates what the classification must distinguish.
Both algorithms are treated in detail in Section~\ref{sec:examples}.

\begin{figure}[htbp]
  \centering
  \includegraphics[width=0.78\textwidth]{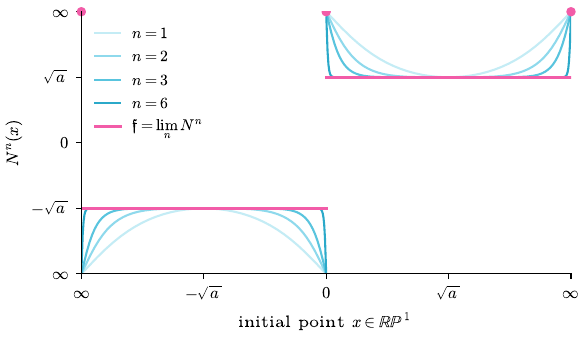}
  \caption{Newton iterates $N_p^n$ for $p(x)=x^2-a$, with $a>0$, on the real projective line, drawn as graphs whose axes run from $\infty$ through~$0$ back to~$\infty$.
    The continuous maps $N_p^n$ converge pointwise to the deep computation $\mathfrak f$, which sends the positive half line to $\sqrt a$, the negative half line to $-\sqrt a$, and $0$ and $\infty$ to~$\infty$.}
  \label{fig:sqrt}
\end{figure}

\emph{Square roots.}
Newton’s method for $p(x)=x^2-a$, with $a>0$, acts on the real projective line $\RP^1=\mathbb{R}\cup\{\infty\}$ by the map $N_p(x)=(x^2+a)/(2x)$, with $N_p(0)=N_p(\infty)=\infty$.
The iterates of~$N_p$ are continuous and converge pointwise to a map~$\mathfrak f$.
This map sends every positive number to~$\sqrt a$, every negative number to~$-\sqrt a$, and both $0$ and $\infty$ to~$\infty$ (see Figure~\ref{fig:sqrt}).
The deep computations are the iterates together with~$\mathfrak f$.
Since the map $\mathfrak f$ is discontinuous, it is not computable, and the set of iterates is not stable.
Nevertheless, the map $\mathfrak f$ is the limit of a \emph{sequence} of computations, and the space of deep computations is countable.

\emph{Newton’s method on the Riemann sphere.}
For $q(z)=z^3-2z+2$, Newton’s method on the Riemann sphere $\hatC=\mathbb{C}\cup\{\infty\}$ displays three distinct modes of behavior, each confined to its own part of the sphere (see Figure~\ref{fig:newton} in Section~\ref{sec:Newton-cubic}).
The iterates converge on the basins of the three roots.
On the basin of the cycle $0\mapsto1\mapsto0$, by contrast, the iterates do not converge, although the even iterates and the odd iterates converge there separately to different limits.
Finally, on the common boundary of these basins, a fractal set, the iterates are chaotic.
For the simpler polynomial $z^2-1$, chaos of this kind produces deep computations that are not even limits of sequences of computations (see Proposition~\ref{prop:Newton-wild} and Theorem~\ref{thm:nip-baire1def}).

The classification of Section~\ref{sec:Classification} is designed to separate such behaviors.
In Theorem~\ref{thm:nip-baire1def} below we characterize, for countably many features, the sets of computations whose deep computations are of the first Baire class on every type shard.
For these sets every deep computation is, on each shard, the limit of a sequence of computations, as in the case of square roots.
Todorčević’s trichotomy, which we recall in Section~\ref{sec:trichotomy}, then measures the complexity of the space of such deep computations.

\section{Classifying deep computations}
\label{sec:Classification}

The main theorems of the paper, proved in this section with the tools of the topology of function spaces, classify sets of computations by the complexity of their deep computations.
The resulting classification yields a new hierarchy of complexity classes.
The examples of Section~\ref{sec:preview} show why this classification is needed.
Indeed, some natural algorithms have discontinuous deep computations that are nevertheless limits of sequences of computations, whereas others have deep computations that are not even limits of sequences.
The main topological tool is Todorčević’s trichotomy for Rosenthal compacta~\cite{Todorcevic_1999_CompactSubsetsBaire}, and our model is the work of Glasner and Megrelishvili on tame semigroups in topological dynamics~\cite{glasner2022tame}.

\subsection{NIP, Rosenthal compacta, and deep computations}

The classification rests on Theorem~\ref{thm:nip-baire1def}, which we prove in this subsection.
This theorem characterizes, by a finitary condition and for countably many features, the sets of computations whose deep computations are of the first Baire class on every type shard.
Analogously, the sets of computations whose deep computations are continuous were characterized in the first paper of this series~\cite{alva2024approximability} by Grothendieck’s double limit criterion.

\begin{thm}\label{thm:nip-baire1def}
  Let $(L,\mathcal{P},\Gamma)$ be a compositional computation structure (see Definition~\ref{def:CCS}) satisfying the Extendibility Axiom, with $\mathcal{P}$ countable.
  Let $\mathscr{R}$ be an exhaustive collection of sizers,
  and suppose that $\Delta\subseteq\Gamma$ is $\mathscr{R}$-confined.
  Then the following conditions are equivalent:

  \begin{enumerate}[(i)]
  \item
    the inclusion $\overline{\tilde\Delta|_{\mathcal{L}[r_\bullet]}}\subseteq \mathrm{B}_1(\mathcal{L}[r_\bullet],\mathcal{L}[r_\bullet])$ holds for all $r_\bullet\in \mathscr{R}$;
  \item
    for all $P\in\mathcal{P}$ and $r_\bullet\in \mathscr{R}$, the set $\pi_P\circ \Delta|_{L[r_\bullet]}$ satisfies the NIP, where $\pi_P\colon\mathbb{R}^{\mathcal{P}}\to\mathbb{R}$ is the projection onto the coordinate~$P$.
    In other words, for all $P\in\mathcal{P}$ and $r_\bullet\in \mathscr{R}$, all reals $a<b$ and all sequences $(\gamma_n)_{n\in\mathbb{N}}$ in~$\Delta$, there are finite disjoint sets $I,J\subseteq\mathbb{N}$ such that
    \[
      L[r_\bullet]\cap\bigcap_{n\in I}(\pi_P\circ\gamma_n)^{-1}(-\infty,a]\cap\bigcap_{n\in J}(\pi_P\circ\gamma_n)^{-1}[b,\infty)=\emptyset.
    \]
  \end{enumerate}

  Furthermore, if these conditions hold, then every deep computation $f\in\overline{\Delta}$ has an extension $\tilde f\colon\mathcal{L}_{\sh}\to \mathcal{L}_{\sh}$ such that $\tilde f|_{\mathcal{L}[r_\bullet]}\in \mathrm{B}_1(\mathcal{L}[r_\bullet],\mathcal{L}[r_\bullet])$ for all $r_{\bullet}\in \mathscr{R}$.
  In particular, on each shard, every deep computation is the pointwise limit of a sequence of computations.
\end{thm}

\begin{proof}
  Since $\mathcal{P}$ is countable, each type shard $\mathcal{L}[r_\bullet]$ is a compact metrizable, hence Polish, subspace of~$\mathbb{R}^\mathcal{P}$.
  Let $r_\bullet\in \mathscr{R}$ be given;
  by the Extendibility Axiom and confinement, every $\tilde\gamma$ with $\gamma\in\Delta$ maps $\mathcal{L}[r_\bullet]$ continuously into itself (see the proof of Lemma~\ref{lem:extension}).
  Hence, for every $P\in\mathcal{P}$, the set $\pi_P\circ \tilde\Delta|_{\mathcal{L}[r_\bullet]}$ consists of continuous functions and is pointwise bounded.
  By Theorem~\ref{thm:Generalized-BFT} and Lemma~\ref{lem:NIP-and-closure}, condition~(ii) holds at $r_\bullet$ if and only if $\overline{\tilde\Delta|_{\mathcal{L}[r_\bullet]}}\subseteq\mathrm{B}_1(\mathcal{L}[r_\bullet],\mathbb{R}^{\mathcal{P}})$, which (i) implies.
  Conversely, suppose that this inclusion holds.
  For a Polish space~$X$, every compact subset of $\mathrm{B}_1(X)$ is Fréchet--Urysohn
  (see~\cite[Theorem~3F]{BFT_1978_PCompactBaire} or~\cite[Theorem~4.1]{debs2013rosenthal}).
  In a Fréchet--Urysohn space, every point in the closure of a subset is the limit of a sequence of points of that subset.
  By Proposition~\ref{prop:Baire-identifications}, the closure $\overline{\tilde\Delta|_{\mathcal{L}[r_\bullet]}}$ may be regarded as a compact subset of $\mathrm{B}_1(\mathcal{L}[r_\bullet]\times\mathcal{P})$, where the space $\mathcal{L}[r_\bullet]\times\mathcal{P}$ is Polish because $\mathcal{P}$ is countable and discrete.
  Every member of this closure is therefore the pointwise limit of a sequence of maps $\tilde\gamma_n|_{\mathcal{L}[r_\bullet]}$ with $\gamma_n\in\Delta$, which are continuous maps from $\mathcal{L}[r_\bullet]$ to itself;
  hence (i) holds at~$r_\bullet$.
  Finally, if (i) holds, then the extension $\tilde f$ of Lemma~\ref{lem:extension} satisfies
  $\tilde f|_{\mathcal{L}[r_\bullet]}\in\overline{\tilde\Delta|_{\mathcal{L}[r_\bullet]}}\subseteq \mathrm{B}_1(\mathcal{L}[r_\bullet],\mathcal{L}[r_\bullet])$ for all $r_\bullet\in \mathscr{R}$.
  By the argument just given, it follows that every deep computation is, on each shard, the pointwise limit of a sequence of computations.
\end{proof}

\subsection{Todorčević’s trichotomy and levels of NIP}
\label{sec:trichotomy}
Todorčević’s trichotomy now enables us to stratify the NIP sets of computations into three levels.
The trichotomy concerns separable Rosenthal compacta, and, by Theorem~\ref{thm:nip-baire1def} and under its hypotheses, the deep computations of a countable $\NIP$ set of computations form such a compactum on each type shard.
We measure the complexity of such a compactum by three topological properties, namely metrizability, hereditary separability and first countability.
The three properties form a chain:
every compact metrizable space is hereditarily separable, and R.~Pol proved that every hereditarily separable Rosenthal compactum is first countable (see~\cite[Section~10]{debs2013rosenthal}).
These properties define three levels of NIP.

\begin{defn}\label{def:NIP-levels}
  Let $(L,\mathcal{P},\Gamma)$ be a CCS satisfying the Extendibility Axiom, with $\mathcal{P}$ countable, and let $\mathscr{R}$ be an exhaustive collection of sizers.
  Let $\Delta\subseteq\Gamma$ be an $\mathscr{R}$-confined countable set of computations satisfying the NIP on shards and features (condition (ii) in Theorem \ref{thm:nip-baire1def}).
  We say that:
  \begin{enumerate}[(i)]
  \item
    the set $\Delta$ is $\NIP_1$ if $\overline{\tilde\Delta|_{\mathcal{L}[r_\bullet]}}$ is first countable for every $r_\bullet\in \mathscr{R}$;
  \item
    the set $\Delta$ is $\NIP_2$ if $\overline{\tilde\Delta|_{\mathcal{L}[r_\bullet]}}$ is hereditarily separable for every $r_\bullet\in \mathscr{R}$;
  \item
    the set $\Delta$ is $\NIP_3$ if $\overline{\tilde\Delta|_{\mathcal{L}[r_\bullet]}}$ is metrizable for every $r_\bullet\in \mathscr{R}$.
  \end{enumerate}
\end{defn}

By the chain of properties above, the levels satisfy $\NIP_3 \Rightarrow \NIP_2\Rightarrow \NIP_1\Rightarrow \NIP$.
Moreover, a countable stable set of computations (see Definition~\ref{def:stable}), which has the NIP by Theorem~\ref{thm:nip-baire1def}, is $\NIP_3$.
Indeed, on each type shard the deep computations of such a set form a separable compact subset of a countable product of spaces $\Cp(\mathcal{L}[r_\bullet])$, that is, a separable Eberlein compactum.
Separable Eberlein compacta, in turn, are metrizable~\cite{Arkhangelskii:1992}.
The five compacta of Section~\ref{sec:five-compacta} show that none of the implications between the three topological properties can be reversed for separable Rosenthal compacta.
Specifically, the extended Alexandroff compactification $\hat{A}(2^{\mathbb{N}})$ is separable but not first countable.
The extended duplicate $\hat{D}(S(2^{\mathbb{N}}))$ of the split Cantor set and the compactum $\hat{D}(2^{\mathbb{N}})$ that contains the Alexandroff duplicate are separable and first countable but not hereditarily separable.
Finally, the split Cantor set $S(2^{\mathbb{N}})$ is hereditarily separable but not metrizable.
Todorčević proved that these examples are canonical~\cite{Todorcevic_1999_CompactSubsetsBaire}, in the following sense.

\begin{thm}[Todorčević’s trichotomy~\cite{Todorcevic_1999_CompactSubsetsBaire}; see also {\cite[Theorem~3]{argyros2008rosenthal}}]\label{thm:Todorcevic-trichotomy}
  Let $K$ be a separable Rosenthal compactum.
  \begin{enumerate}[(i)]
  \item
    If $K$ is hereditarily separable but non-metrizable, then the split Cantor set $S(2^\mathbb{N})$ embeds into $K$.
  \item
    If $K$ is first countable but not hereditarily separable, then either the Alexandroff duplicate $D(2^\mathbb{N})$ or the extended duplicate $\hat{D}(S(2^\mathbb{N}))$ of the split Cantor set embeds into $K$.
  \item
    If $K$ is not first countable, then the Alexandroff compactification $A(2^\mathbb{N})$ embeds into $K$.
  \end{enumerate}
\end{thm}

Together with the chain of properties above, the trichotomy yields the following classification of separable Rosenthal compacta:

\begin{center}
  \begin{tikzpicture}[
      grow=down,
      sibling distance=10em,
      level distance=4em,
      edge from parent/.style={draw,-latex},
      every node/.style={font=\small,align=center}
    ]

    \node {$K$ is a separable Rosenthal compactum}
    child {node {$K$ is metrizable}
      }
    child {node {$K$ is not metrizable}
        child {node {$K$ is hereditarily separable \\ (copy of $S(2^\mathbb{N})$) \hspace{2cm}}}
        child {node {\hspace{2cm} $K$ is not hereditarily separable}
            child {node {$K$ is first countable \\ (copy of $D(2^\mathbb{N})$ or $\hat{D}(S(2^\mathbb{N}))$)}}
            child {node {$K$ is not first countable \\ (copy of $A(2^\mathbb{N})$)}}
          }
      };

  \end{tikzpicture}
\end{center}

The three levels admit a second description, in terms of the number of tests required to distinguish deep computations from one another.
Let $\Delta$ and $\mathscr{R}$ be as in Definition~\ref{def:NIP-levels}, fix $r_\bullet\in \mathscr{R}$, and write $K_{r_\bullet}=\overline{\tilde\Delta|_{\mathcal{L}[r_\bullet]}}$.
We identify each $f\in K_{r_\bullet}$ with the real-valued function $(x,P)\mapsto P(f(x))$ on $\mathcal{L}[r_\bullet]\times\mathcal{P}$, whose points we call \emph{tests}:
the test $(x,P)$ applies the feature $P$ to the output on the input~$x$.
The topology of pointwise convergence on $K_{r_\bullet}$ is then the topology of pointwise convergence on tests.
A set $\mathcal{T}$ of tests \emph{separates} a set of deep computations if any two distinct members of the set differ on some test in~$\mathcal{T}$.
A member $f$ of $K_{r_\bullet}$ is \emph{countably determined} if some countable set $\mathcal{T}_f$ of tests separates $f$ from every other member of $K_{r_\bullet}$.
A set $E\subseteq K_{r_\bullet}$ is \emph{finitely separated} if each $f\in E$ admits finitely many tests and some $\varepsilon_f>0$ such that every $g\in E\setminus\{f\}$ differs from $f$ by at least $\varepsilon_f$ on one of these tests.

\begin{prop}\label{prop:levels-tests}
  Let $\Delta$ and $\mathscr{R}$ be as in Definition~\ref{def:NIP-levels}, let $r_\bullet\in \mathscr{R}$ be given, and let $K_{r_\bullet}=\overline{\tilde\Delta|_{\mathcal{L}[r_\bullet]}}$.
  \begin{enumerate}[(i)]
  \item
    The compactum $K_{r_\bullet}$ is metrizable if and only if one countable set of tests separates $K_{r_\bullet}$.
  \item
    It is first countable if and only if every member of $K_{r_\bullet}$ is countably determined.
  \item
    It is hereditarily separable if and only if no finitely separated subset of $K_{r_\bullet}$, that is, no subset that is discrete in the topology of pointwise convergence on tests, has cardinality $2^{\aleph_0}$.
  \end{enumerate}
\end{prop}

\begin{proof}
  (i) If a countable set $\mathcal{T}$ of tests separates $K_{r_\bullet}$, then restriction to $\mathcal{T}$ is a continuous injection of the compact space $K_{r_\bullet}$ into the metrizable space $\mathbb{R}^{\mathcal{T}}$.
  This injection is therefore a homeomorphism onto its image.
  Conversely, if $K_{r_\bullet}$ is metrizable, then the square $K_{r_\bullet}\times K_{r_\bullet}$ is second countable, and hence its open subset $\{(f,g):f\ne g\}$ is Lindelöf.
  This open set is covered by the open sets $\{(f,g):f(t)\ne g(t)\}$ for tests~$t$, and a countable subcover gives the required set $\mathcal{T}$.
  (ii) If $f$ is countably determined, the singleton $\{f\}$ is a countable intersection of open sets, and in a compact Hausdorff space such a point has a countable base of neighborhoods.
  Conversely, let $(U_k)_{k\in\mathbb{N}}$ be a countable base of neighborhoods of $f$, and choose in each $U_k$ a basic neighborhood of $f$ defined by a finite set $F_k$ of tests.
  The union $\mathcal{T}_f=\bigcup_kF_k$ then separates $f$ from every other member of $K_{r_\bullet}$.
  (iii) A set is finitely separated exactly when it is a discrete subspace, and, by Theorem~\ref{thm:nip-baire1def}, the compactum $K_{r_\bullet}$ is a separable Rosenthal compactum.
  A hereditarily separable space has no uncountable discrete subspace, since a separable discrete space is countable.
  Conversely, suppose that $K_{r_\bullet}$ is not hereditarily separable;
  then Theorem~\ref{thm:Todorcevic-trichotomy} provides a copy of $D(2^{\mathbb{N}})$, of $\hat{D}(S(2^{\mathbb{N}}))$ or of $A(2^{\mathbb{N}})$ inside $K_{r_\bullet}$.
  Each of these three compacta contains a discrete subspace of cardinality $2^{\aleph_0}$ (see Section~\ref{sec:five-compacta}).
  For separable Rosenthal compacta, the existence of such a subspace in the absence of hereditary separability is due to R.~Pol (see~\cite{argyros2008rosenthal} and~\cite[Theorem~0.4]{glasner2022tame}).
\end{proof}

Thus the set $\Delta$ is $\NIP_3$ when, on each type shard, a single countable set of tests separates all the deep computations of~$\Delta$.
It is $\NIP_2$ when no type shard carries a finitely separated set of deep computations of the cardinality of the continuum.
Finally, the set $\Delta$ is $\NIP_1$ when, on each type shard, each deep computation is separated from all the others by countably many tests of its own.

In Section~\ref{sec:iterative} we exhibit computations at each level of the hierarchy and show that every implication between the classes stable, $\NIP_3$, $\NIP_2$, $\NIP_1$ and $\NIP$ is strict.
The definitions of $\NIP_1$, $\NIP_2$ and $\NIP_3$ are topological, and the descriptions of Proposition~\ref{prop:levels-tests} still involve infinitely many tests.
Whether the three levels nevertheless admit a finitary, combinatorial characterization, like that of $\NIP$ in Theorem~\ref{thm:nip-baire1def}, is Question~\ref{q:NIP-levels}.

\subsection{Canonical families and binary search}
\label{sec:heptachotomy}

Todorčević's trichotomy organizes the hierarchy of the preceding subsection.
For a finer interpretation of its canonical examples, we turn to the seven tree-indexed families studied by Argyros, Dodos and Kanellopoulos~\cite{argyros2008rosenthal}.
Their prefix and threshold functions admit a direct reading as primitives of binary search, and their closures illustrate distinct patterns of limiting behavior.
Recall from Section~\ref{sec:five-compacta} that a dyadic subtree of the binary tree $2^{<\mathbb{N}}$ is written $\{w_t:t\in2^{<\mathbb{N}}\}$, indexed through the canonical isomorphism preserving extension and the lexicographic order.

\begin{defn}
  Let $X$ and $Y$ be Polish spaces.

  \begin{enumerate}
  \item
    If $I$ is countable and $\{f_i:i\in I\}\subseteq\mathbb{R}^X$, $\{g_i:i\in I\}\subseteq\mathbb{R}^Y$ are two families indexed by $I$, we say that they are \emph{equivalent} if the map $f_i\mapsto g_i$ extends to a homeomorphism from $\overline{\{f_i:i\in I\}}$ to $\overline{\{g_i:i\in I\}}$.
  \item
    If $\{f_t:t\in 2^{<\mathbb{N}}\}\subseteq\mathbb{R}^X$ is a pointwise bounded family, we say that it is \emph{minimal} if, for every dyadic subtree $\{w_t:t\in 2^{<\mathbb{N}}\}$ of $2^{<\mathbb{N}}$, the family $\{f_{w_t}:t\in 2^{<\mathbb{N}}\}$ is equivalent to $\{f_t:t\in 2^{<\mathbb{N}}\}$.
  \end{enumerate}
\end{defn}

The seven canonical minimal families in the classification of Argyros, Dodos and Kanellopoulos~\cite{argyros2008rosenthal} are described below.
We retain their notation and the tree indexing, which distinguishes families even when their closures are homeomorphic.
The description is carried out along the regular dyadic subtree $\{s_t:t\in 2^{<\mathbb{N}}\}$ fixed in Section~\ref{sec:five-compacta} for the definition of the extended duplicate of the split Cantor set.
This subtree satisfies $s_t^\frown 0^{\infty}\neq s_{t'}^\frown 0^\infty$ and $s_t^\frown 1^{\infty}\neq s_{t'}^\frown 1^\infty$ whenever $t\neq t'$.
The families are built from the prefix tests $v_t$ and the threshold tests $\chi_{[a}$, $\chi_{(a}$ of Section~\ref{sec:five-compacta}.
Two further conventions are taken from the paper of Argyros, Dodos and Kanellopoulos~\cite{argyros2008rosenthal}.
First, a point $\sigma\in2^{\mathbb{N}}$, identified with a point of the Cantor set in~$[0,1]$, is regarded as the constant function on~$2^{\mathbb{N}}$ with value~$\sigma$.
Second, a pair $(f,g)$ of functions on $2^{\mathbb{N}}$ is regarded as a single function on the disjoint union $2^{\mathbb{N}}\sqcup2^{\mathbb{N}}$.

\begin{enumerate}
\item
  The family $D_1=\{\frac{1}{|t|+1}v_t:t\in 2^{<\mathbb{N}}\}$ has closure $D_1\cup\{\mathbf{0}\}$, the one-point compactification of a countable discrete set.
\item
  The family $D_2=\{s_t^\frown 0^{\infty}:t\in 2^{<\mathbb{N}}\}$ is a discrete subset of its closure $\overline{D_2}$, which Argyros, Dodos and Kanellopoulos denote by $2^{\leq\mathbb{N}}$.
  This closure is compact and metrizable, and its accumulation points form a copy of~$2^{\mathbb{N}}$.
\item
  The family $D_3=\{\chi_{[s_t^\frown 0^\infty}:t\in 2^{<\mathbb{N}}\}$ is a discrete subset of its closure $\overline{D_3}$, a closed subspace of~$S(2^{\mathbb{N}})$.
\item
  The family $D_4=\{\chi_{(s_t^\frown 1^\infty}:t\in 2^{<\mathbb{N}}\}$ is a discrete subset of its closure $\overline{D_4}$, a closed subspace of~$S(2^{\mathbb{N}})$.
\item
  The family $D_5=\{v_t:t\in 2^{<\mathbb{N}}\}$ is a discrete subset of its closure $\overline{D_5}=\hat{A}(2^{\mathbb{N}})$.
\item
  The family $D_6=\{(v_t,t^\frown 0^\infty):t\in 2^{<\mathbb{N}}\}$ is a discrete subset of its closure $\overline{D_6}=\hat{D}(2^{\mathbb{N}})$.
\item
  The family $D_7=\{(v_{s_t},\chi_{[s_t^\frown 0^\infty}):t\in 2^{<\mathbb{N}}\}$ is a discrete subset of its closure $\overline{D_7}=\hat{D}(S(2^{\mathbb{N}}))$.
\end{enumerate}

\begin{prop}[Canonical minimal families {\cite[Section~4.3]{argyros2008rosenthal}}]
  \label{prop:canonical-search-families}
  The seven families $D_1,\dots,D_7$ defined above are relatively compact in the first Baire class on their respective Polish domains.
  They are minimal and pairwise nonequivalent under the tree-indexed equivalence just defined.
\end{prop}

The building blocks of the families $D_1,\dots,D_7$ admit a computational interpretation as primitives of binary search on $2^{\mathbb{N}}$, the space of binary expansions.
A node $t\in2^{<\mathbb{N}}$ determines the dyadic interval $[t^\frown0^\infty,t^\frown1^\infty]$, the set of inputs that round to $t$ at precision $|t|$.
The prefix test $v_t$ asks whether the input rounds to $t$, and the point $t^\frown0^\infty$ is the rounded value itself.
The comparison $\chi_{[t^\frown0^\infty}$ asks whether the input lies at or to the right of the left endpoint of the interval.
Likewise, the comparison $\chi_{(t^\frown1^\infty}$ asks whether the input lies strictly to the right of the right endpoint.
The resulting reading of the seven families is recorded in Table~\ref{tab:hepta-search}.

\begin{table}[ht]
  \centering
  \small
  \renewcommand{\arraystretch}{1.25}
  \begin{tabular}{>{\raggedright\arraybackslash}p{1.1cm}>{\raggedright\arraybackslash}p{4.3cm}>{\raggedright\arraybackslash}p{5.9cm}}
    \hline
    Family & Function at node $t$                        & Limits                                                       \\
    \hline
    $D_1$  & prefix test weighted by $1/(|t|+1)$         & the test that never succeeds: confidence vanishes with depth \\
    $D_2$  & return the rounded value $t^\frown0^\infty$ & the exact input: successive approximation                    \\
    $D_3$  & compare with the left endpoint              & comparison with a limit of left endpoints, with ties broken upward and, at limits approached from above, also downward\\
    $D_4$  & compare with the right endpoint             & comparison with a limit of right endpoints, with ties broken downward and, at limits approached from below, also upward\\
    $D_5$  & prefix test                                 & exact equality, or no match                                  \\
    $D_6$  & prefix test and rounded value               & a dictionary lookup that returns its key                     \\
    $D_7$  & prefix test and comparison                  & a three-way comparison: below, inside or above               \\
    \hline
  \end{tabular}
  \vspace{6pt}
  \caption{The minimal families of Argyros, Dodos and Kanellopoulos as primitives of binary search.}
  \label{tab:hepta-search}
\end{table}

The translation to a concrete computation is made after fixing the magnitude bounds specified by its parameters.
Let $r_\bullet$ be the resulting sizer, and let $\{\gamma_t:t\in2^{<\mathbb{N}}\}$ be a family of computations in a CCS satisfying the Extendibility Axiom and having countably many features.
Suppose that the family preserves the shard $L[r_\bullet]$ and has NIP there for every feature.
For a chosen feature $P$, the argument of Theorem~\ref{thm:nip-baire1def} supplies a relatively compact family of scalar functions $\{P\circ\tilde\gamma_t|_{\mathcal{L}[r_\bullet]}:t\in2^{<\mathbb{N}}\}$ in the first Baire class on the type shard $\mathcal{L}[r_\bullet]$.
Restriction to a shard or observation through a single feature can identify distinct computations; a constant feature, for example, gives a constant family.
Accordingly, the seven families are used here as canonical patterns of observable limit behavior, whose occurrence in a particular computation requires a separate identification.

\section{Examples of  Compositional Computation Structures}
\label{sec:examples}
From the general theory we turn to the examples that motivate it, among them the two announced in Section~\ref{sec:preview}, and describe the deep computations of three compositional computation structures.
These three structures are Newton’s method on the Riemann sphere for two cubic polynomials (Sections~\ref{sec:Newton-cubic} and~\ref{sec:Newton-cyclotomic}) and Newton’s method for square roots on the real projective line (Section~\ref{exa:square-roots-compute}).
In Section~\ref{sec:iterative} we place the last of these, together with other classical algorithms, at the levels of the hierarchy of Section~\ref{sec:Classification}.

We begin by formalizing the Newton--Raphson root-finding method as a CCS~\cite{Galantai200025}, \cite[Section~2.3]{Burden2016numerical}.
Let $\hatC = \mathbb{C}\cup\{\infty\}$ be the extended complex plane, that is, the Alexandroff (one-point) compactification of~$\mathbb{C}$.
We identify $\hatC$ with the unit sphere $S^2\subseteq \mathbb{R}^3$, the Riemann sphere, through the inverse stereographic projection.
This projection is the homeomorphism $\hatC\to S^2$ given by $z \mapsto (P_1(z), P_2(z), P_3(z))$, where $(P_1(\infty), P_2(\infty), P_3(\infty)) = (0,0,1)$ and, for $z\ne\infty$,
\begin{align*}
P_1(z) & = \frac{2\,\Re{z}}{|z|^2+1}, &
P_2(z) & = \frac{2\,\Im{z}}{|z|^2+1}, &
P_3(z) & = \frac{|z|^2-1}{|z|^2+1}.
\end{align*}
The extended complex plane $\hatC$ thus carries the real predicates $P_1, P_2, P_3$ with values in~$[-1,1]$.
Through the same identification, the topology of $\hatC$ is metrized by the pullback of the Euclidean metric of~$S^2$.

\subsubsection*{Newton  CCSs:  steps and fractals}
\label{sec:Newton-map}
Here we describe the compositional computation structure generated by the Newton map, the single step that Newton's method iterates to approximate a root of a function.
The structure so obtained is the common setting of the three examples that follow.
We then recall the Fatou and Julia sets, on which the iterates of the Newton map behave regularly and chaotically, respectively.
Let $p = p(z)$ be a nonconstant rational function with complex coefficients, that is, a nonconstant meromorphic map $\hatC\to\hatC$.
The \emph{Newton step}, or \emph{Newton map}, of $p$ is the rational map
\begin{equation}\label{eq:Np}
N_p: z \mapsto z-\frac{p(z)}{p'(z)}.
\end{equation}
The iterates of the Newton map, that is, its compositional powers~$N_p^n$, may be regarded as the computations of a CCS with state space $L = \hatC$.
The natural predicates are the coordinates $P_1, P_2, P_3$ of the inverse stereographic projection, since the real and imaginary parts $z\mapsto\Re{z}$ and $z\mapsto\Im{z}$ are not defined at~$z = \infty$.

Let $\Gamma_p\coloneqq\{N_p^n:n\in\mathbb{N}\}$ be the semigroup of iterates of the Newton map~\eqref{eq:Np}.
The CCS $(\hatC, (P_1,P_2,P_3), \Gamma_p)$ is called \emph{Newton’s method for~$p(z)$}, or the \emph{Newton CCS for~$p(z)$}.
Since the state space $\hatC$ is compact and the Newton map $N_p$ is continuous on it, every Newton CCS satisfies the Extendibility Axiom.

Newton’s method succeeds at an initial estimate $z_0\in\hatC$ when the sequence $(z_n)_{n\in \mathbb{N}}$ of Newton steps, where $z_n \coloneqq N_p^n(z_0)$, converges to a root of~$p$.

\paragraph{Fatou and Julia sets}
Recall that the \emph{Fatou set $\mathcal{F}_g$} of a rational map~$g\colon\hatC\to\hatC$ is the set of points near which the iterates of~$g$ behave regularly.
More precisely, a point belongs to~$\mathcal{F}_g$ if it has a neighborhood on which every sequence of iterates of~$g$ has a subsequence that converges uniformly with respect to the spherical metric.
The complement $\mathcal{J}_g\coloneqq\hatC\setminus \mathcal{F}_g$ of the Fatou set is the \emph{Julia set} of $g$, a closed set that is typically a fractal.
The Julia set of the Newton map~$N_p$ is called the \emph{Newton fractal} of~$p$.

\subsection{Example: The Newton CCS for the cubic $q(z) = z^3-2z+2$}
\label{sec:Newton-cubic}
Our first example is the second motivating example of Section~\ref{sec:preview}, a cubic polynomial whose Newton map has, besides the three roots, an attracting cycle.
Newton's method therefore fails on an open set of initial estimates, and the iterates of the Newton map do not converge.
For the same reason, two distinct limits of the iterates along ultrafilters are among the deep computations of this method.
One is taken along an ultrafilter that contains the even numbers, and the other along an ultrafilter that contains the odd numbers.
The polynomial $q(z) = z^3-2z+2$ has Newton map
\begin{equation*}
N_q: z \mapsto z - \frac{z^3-2z+2}{3z^2-2} = \frac{2z^3-2}{3z^2-2}.
\end{equation*}
Denote the real root of $q$ by $r_0$, and the roots in the upper and the lower half plane by $r_{+}$ and $r_{-}$.
Starting from $z_0=0$, the Newton steps alternate $0, 1, 0, 1, \dots$, and neither $0$ nor $1$ is a root of~$q$.

\begin{figure}[htbp]
  \centering
  \includegraphics[width=\textwidth]{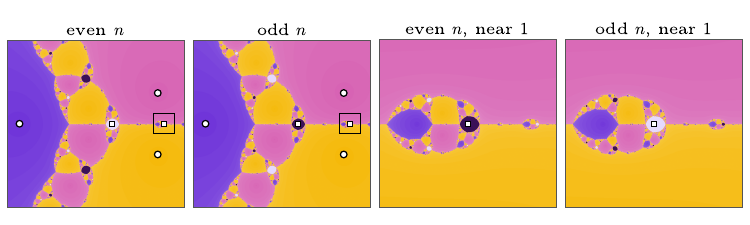}
  \caption{The deep computations of Newton’s method for $q(z)=z^3-2z+2$: the limits $f_0$ and $f_1$ of the iterates along $\mathcal{U}_0$ and along $\mathcal{U}_1$, on the whole plane and near the point~$1$.
    Points in the basins of the real root and of the roots in the upper and lower half planes (white circles) are purple, orchid and gold, paler where convergence is slower.
    The remaining open region is the basin of the superattracting cycle $0\mapsto1\mapsto0$ (white squares).
    A point there is lilac if the iterate in question is near~$0$ and plum if it is near~$1$, whence the two colors exchange between even and odd iterates.
    Accordingly, the two deep computations agree on the basins of the roots and differ on the basin of the cycle.
    The boundaries of the colored regions form the Newton fractal of~$q$.}
  \label{fig:newton}
\end{figure}

Let $B_0, B_+, B_-$ be the basins of the roots $r_0, r_+, r_-$, and let $B_{01}$ be the basin of the superattracting cycle $\{0,1\}$.
The Fatou set of $N_q$ is $B_{01} \cup B_0 \cup B_{+} \cup B_{-}$, and each of these four sets is nonempty and open.
Since $B_{01}$ is nonempty, the iterates $N_q^n$ do not converge pointwise.
Fix nonprincipal ultrafilters $\mathcal{U}_0$ and $\mathcal{U}_1$ on~$\mathbb{N}$ that contain the even and the odd numbers, respectively.
The limits of the iterates $N_q^n$ along~$\mathcal{U}_0$ and along~$\mathcal{U}_1$ are deep computations, which we denote by $f_0$ and~$f_1$, respectively.
On the Fatou set the even and the odd iterates converge pointwise, whence the two limits, on that set, do not depend on the choice of the ultrafilters.
Both limits send each basin of a root to that root; on~$B_{01}$, by contrast, they take only the values $0$ and~$1$ and exchange them.
The exchange of the two values follows from the identity $f_1=N_q\circ f_0$, which holds on the Fatou set by the continuity of~$N_q$.
In particular, the two deep computations are distinct.%
\footnote{The point $\infty$ is a repelling fixed point of~$N_q$: if $|z|$ is large, then $N_q(z) = \frac{2}{3}z + O(1)$.
  Hence the point $\infty$ lies on the Newton fractal of~$q$, and the two deep computations satisfy $f_0(\infty) = \infty = f_1(\infty)$.}
Figure~\ref{fig:newton} depicts $f_0$ and $f_1$, together with the Newton fractal of~$q$, on which the iterates behave chaotically.

\subsection{Example: The Newton CCS for the cubic $u(z) = z^{3}-1$}
\label{sec:Newton-cyclotomic}
\begin{figure}[htbp]
  \centering
  \includegraphics[width=\textwidth]{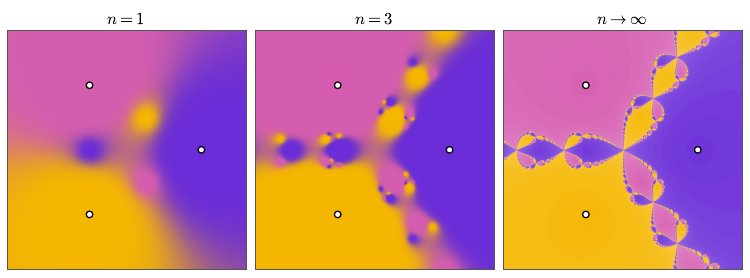}
  \caption{Newton’s method for $u(z)=z^3-1$ after one and three iterations, and the limit of its iterates on the Fatou set.
    The coloring blends the colors of the three cube roots of unity (white circles) according to the distance of $N_u^n(z)$ to each of them.
    This limit is constant on each basin and changes value across the Newton fractal;
    in the right panel, points are paler where convergence is slower.}
  \label{fig:newton2}
\end{figure}

For contrast, we turn to $u(z) = z^3-1$, whose roots are the cube roots of unity $r_0=1$ and $r_{\pm} = \exp(\pm2\pi i/3)$.
The Newton map of~$u$ is
\begin{equation*}
N_u: z\mapsto z-\frac{z^3-1}{3z^2}=\frac{2z^3+1}{3z^2}.
\end{equation*}
The critical points of this rational map of degree three are the three roots of~$u$, which are superattracting fixed points, and the point~$0$.
The map $N_u$ sends the point~$0$ to~$\infty$, a repelling fixed point of multiplier~$3/2$.
The immediate basin of every attracting or parabolic cycle contains a critical point whose orbit converges to that cycle~\cite[Sections~8 and~10]{Milnor2006}.
Unlike~$N_q$, the map $N_u$ therefore has no attracting cycle other than its three attracting fixed points, the roots of~$u$, and it has no parabolic cycle.
The boundary of a Siegel disk or of a Herman ring, moreover, is infinite and lies in the closure of the critical orbits~\cite[Sections~11 and~15]{Milnor2006}.
Since every critical orbit of~$N_u$ is finite, the map $N_u$ has no such domain.
Sullivan's theorem excludes wandering domains, and the classification of periodic Fatou components then shows that the Fatou set of~$N_u$ is the union of the basins of the three roots~\cite[Section~16 and Appendix~F]{Milnor2006}.
On each of these basins, the iterates of~$N_u$ converge to the corresponding root.
Figure~\ref{fig:newton2} depicts Newton’s method for~$u$.
To render the iterates themselves, and not merely their limits, we color each point $z_0$ by blending the colors of the three roots with relative weights $|z_n-r|^{-3}$ for $r\in\{r_0,r_+,r_-\}$, where $z_n\coloneqq N_u^n(z_0)$.
This coloring is a continuous function of~$z_0$, blurred for small~$n$.
As $n$ increases, the coloring sharpens into large regions of nearly constant color, separated by narrow regions in which the color varies so rapidly that the coloring is nearly discontinuous.

The chaos of Newton’s method on its Julia set is visible to the predicates $P_1,P_2,P_3$.
In Propositions~\ref{prop:Newton-wild} and~\ref{prop:rational-IP} below we show that, for the polynomials $z^2-1$, $q$ and~$u$, this chaos prevents the Newton CCS on the Riemann sphere from being $\NIP$.

\subsection{Example: Square roots}\label{exa:square-roots-compute}
We return to the first motivating example of Section~\ref{sec:preview}, Newton’s method on the real projective line for the polynomial $p(x) = x^2-a$ with real $a>0$, whose iterates approximate the square roots $\pm\sqrt a$.
The deep computations of this method, which we now work out, place it at the level $\NIP_3$ (see Section~\ref{sec:ex-NIP3}).
As for every rational function with real coefficients, the Newton map $N_p$ maps the real projective line $\RP^1 = \mathbb{R}\cup\{\infty\}$ into itself.
The Newton CCS for $p$ may therefore be regarded as having state space~$\RP^1$.

Let $\mathbb{R}\to S^1$, $x\mapsto (P_1(x),P_3(x))$, be the inverse stereographic projection of the real line onto $S^1\subseteq \mathbb{R}^2$, given by
\begin{align*}
P_1(x) & = \frac{2x}{x^2+1},    &
P_3(x) & = \frac{x^2-1}{x^2+1},
\end{align*}
and extended to a homeomorphism of the real projective line $\RP^1 = \mathbb{R}\cup\{\infty\}$ onto~$S^1$ by $\infty \mapsto (P_1(\infty), P_3(\infty)) = (0, 1)$.
These coordinates are the restrictions to $\RP^1$ of the predicates $P_1$ and $P_3$ of the Riemann sphere, whose names they keep.
As predicates, the two coordinates make $(\RP^1, (P_1, P_3))$ a computation states structure with a compact state space.
The Newton map $N_p$ is
\begin{equation*}
x \mapsto x - \frac{x^2-a}{2x}
=
\begin{cases}
\frac{x^2+a}{2x}, & x\ne0,\infty; \\
\infty,           & x=0, \infty.
\end{cases}
\end{equation*}
From every initial point $x_0\ne0,\infty$, the Newton steps converge to one of the square roots $\pm\sqrt{a}$.
More precisely, the iterates of $N_p$ converge pointwise to the map $\mathfrak{f}\colon \RP^1\to\RP^1$ given by
\[
  \mathfrak{f}(x)
  =
  \lim_{n\rightarrow\infty}N_p^n(x)
  =
  \begin{cases}
  \sqrt{a},  & x\in(0,\infty); \\
  -\sqrt{a}, & x\in(-\infty,0); \\
  \infty,    & x = 0, \infty.
  \end{cases}
\]
The map $\mathfrak{f}$ is therefore discontinuous at $0$ and at~$\infty$ (see Figure~\ref{fig:sqrt}).

\section{Computations across the hierarchy}
\label{sec:iterative}

The hierarchy of Section~\ref{sec:Classification} descends from the stable sets of computations through the classes $\NIP_3$, $\NIP_2$ and $\NIP_1$ to the class $\NIP$.
The implications of the hierarchy hold for countable sets of computations in compositional computation structures with countably many features.
Although its levels are defined through the topology of the space of deep computations, each of them is attained by classical algorithms;
some familiar algorithms lie beyond the hierarchy.
In this section we establish both facts, and with them the strictness of every implication.
In each example the computations are compositions of the steps of an algorithm.
At all levels but one, these computations are the powers of a single step map.
At the level $\NIP$ but not $\NIP_1$ they form instead a group generated by two steps, both in the Euclidean algorithm and in free reduction in the free group on two generators.
Moreover, a countable family of elements of the group generated by the steps of the Euclidean algorithm provides a second example at the level $\NIP_1$ but not $\NIP_2$.
The level of each algorithm is then a theorem about the closure of its computations, and for such computations the hierarchy coincides with the Glasner--Megrelishvili hierarchy of tame dynamical systems~\cite{glasner2022tame}.
We begin with this coincidence (see Proposition~\ref{prop:cascade-tame});
Table~\ref{tab:hierarchy} at the end of the section summarizes the examples.

\subsection{Iterative algorithms and tame dynamics}
\label{sec:dynamical-CCS}
The algorithms of this section act by continuous steps on a compact space of states.
Within this setting, which we now isolate, our hierarchy coincides with the Glasner--Megrelishvili hierarchy of tame dynamical systems.
The coincidence allows the results of Glasner and Megrelishvili to be brought to bear on the examples that follow.
A \emph{dynamical CCS} is a compositional computation structure $(L,\mathcal{P},\Gamma)$ in which $L$ is a compact metrizable space and the predicates in~$\mathcal{P}$ are continuous and separate the points of~$L$.
Moreover, the family~$\mathcal{P}$ and the semigroup~$\Gamma$ are countable, and the members of~$\Gamma$ are continuous maps $L\to L$, among them the identity.
A dynamical CCS is a \emph{cascade CCS} if $\Gamma=\{\gamma^n:n\in\mathbb{N}\}$ is the semigroup of iterates of a single continuous map $\gamma\colon L\to L$.
Since $\mathcal{P}$ separates the points of the compact space~$L$, the map $v\mapsto\tp(v)$ is a homeomorphism of $L$ onto the space of types.
Consequently a single shard covers~$L$, every type is realized, and the Extendibility Axiom holds with $\tilde\gamma=\gamma$ for every $\gamma\in\Gamma$.
The deep computations of a dynamical CCS form the pointwise closure
\[
E(L,\Gamma)=\overline{\Gamma}\subseteq L^L ,
\]
which is the \emph{enveloping semigroup}, introduced by Ellis, of the dynamical system $(L,\Gamma)$;
for a cascade CCS we write $E(L,\gamma)$.
The system $(L,\Gamma)$ is called \emph{tame} if, for every continuous function $f\colon L\to\mathbb{R}$, the family $\{f\circ\gamma:\gamma\in\Gamma\}$ has the $\NIP$ of Definition~\ref{def:NIP-fns}~\cite{Kohler1995,Glasner2006}.
We say that a dynamical CCS is \emph{stable} when its set~$\Gamma$ of computations has the property of Definition~\ref{def:stable}.
Similarly, the CCS is called $\NIP$ when $\Gamma$ satisfies condition~(ii) of Theorem~\ref{thm:nip-baire1def} with respect to a single sizer~$r_\bullet$ with $L[r_\bullet]=L$.
With respect to the same sizer, it is called $\NIP_1$, $\NIP_2$ or $\NIP_3$ when $\Gamma$ has the corresponding property of Definition~\ref{def:NIP-levels}.

A map $\xi\colon L\to L$ is \emph{of Baire class~$1$} if $f\circ \xi$ is of
Baire class~$1$ for every continuous $f\colon L\to\mathbb{R}$.
It suffices that $P\circ \xi$ be of Baire class~$1$ for every $P\in\mathcal{P}$.
Indeed, since the type map embeds $L$ into $\mathbb{R}^{\mathcal{P}}$ with $\mathcal{P}$ countable, every open subset of $L$ is a countable union of finite intersections of sets $P^{-1}(U)$ with $P\in\mathcal{P}$ and $U\subseteq\mathbb{R}$ open.
By Fact~\ref{baire}, the preimage of every open set under $\xi$ is then an $F_\sigma$ set, and so is the preimage of every open set under $f\circ \xi$.
A second appeal to Fact~\ref{baire} shows that the function $f\circ\xi$ is of Baire class~$1$.
In general, the definition of Baire class~$1$ given above for maps from~$L$ to~$L$ is less restrictive than the one adopted in Section~\ref{sec:prelim}.
The latter requires a map from~$L$ to~$L$ to be the pointwise limit of a sequence of continuous maps from~$L$ to~$L$.
For example, if $L=[0,1]\cup\{2\}\subseteq\mathbb{R}$, then the map that sends $[0,1/2]$ to~$0$ and the rest of~$L$ to~$2$ is of Baire class~$1$ in the sense defined above.
This map, however, is not the pointwise limit of a sequence of continuous maps from~$L$ to~$L$, because each such map sends the connected set $[0,1]$ into $[0,1]$ or to~$2$.
For the deep computations of a dynamical CCS, by contrast, the two definitions agree.
If every deep computation is of Baire class~$1$ in the sense defined above, then the CCS is $\NIP$, as we show in Proposition~\ref{prop:cascade-tame} below.
In that case the last assertion of Theorem~\ref{thm:nip-baire1def}, applied to the single sizer~$r_\bullet$ with $L=L[r_\bullet]$, shows that every deep computation is the pointwise limit of a sequence of computations.
Since computations are continuous maps from~$L$ to~$L$, every deep computation is then of Baire class~$1$ in the sense of Section~\ref{sec:prelim} as well.

\begin{prop}\label{prop:cascade-tame}
A dynamical CCS is $\NIP$ if and only if the dynamical system $(L,\Gamma)$ is
tame, and if and only if every deep computation is of Baire class~$1$.
In that case the CCS is $\NIP_1$, $\NIP_2$ or $\NIP_3$ if and only if $E(L,\Gamma)$ is
first countable, hereditarily separable or metrizable, respectively.
\end{prop}

\begin{proof}
The single sizer $r_\bullet$ with $L=L[r_\bullet]$ forms an exhaustive
collection that confines~$\Gamma$.
If $(L,\Gamma)$ is tame, the CCS is $\NIP$: take $f=P$ for $P\in\mathcal{P}$.
If the CCS is $\NIP$, then Theorem~\ref{thm:nip-baire1def} and
Proposition~\ref{prop:Baire-identifications} show that the function $P\circ \xi$ is of Baire
class~$1$ for every $\xi\in E(L,\Gamma)$ and $P\in\mathcal{P}$.
Every deep computation is therefore of Baire class~$1$.
Finally, suppose that every deep computation is of Baire class~$1$, and let
$f\colon L\to\mathbb{R}$ be continuous.
The map $\xi\mapsto f\circ \xi$ is continuous on $L^L$, so the closure of
$\{f\circ\gamma:\gamma\in\Gamma\}$ is the compact set
$\{f\circ \xi:\xi\in E(L,\Gamma)\}$, which lies in $\mathrm{B}_1(L)$.
By Theorem~\ref{thm:BFT2}, the family $\{f\circ\gamma:\gamma\in\Gamma\}$ has the
$\NIP$, so $(L,\Gamma)$ is tame.
The last statement is Definition~\ref{def:NIP-levels}, since
$\overline{\tilde\Gamma|_{L}}=E(L,\Gamma)$.
\end{proof}

For metrizable systems, Glasner and Megrelishvili call a tame system \emph{tame$_1$} when its enveloping semigroup is first countable, \emph{tame$_2$} when it is hereditarily separable, and \emph{hereditarily non-sensitive} when it is metrizable.
They prove the inclusions between these classes~\cite[Proposition~0.6]{glasner2022tame} and show by examples that the inclusions are proper.
It follows from Proposition~\ref{prop:cascade-tame} that a dynamical CCS is $\NIP_1$, $\NIP_2$ or $\NIP_3$ exactly when its dynamical system is tame$_1$, tame$_2$ or hereditarily non-sensitive, respectively.

The simplest iterative methods are those whose iterates converge pointwise.
In Proposition~\ref{prop:convergent} we show that all of them are $\NIP_3$, and that the continuity of the limit decides whether they are stable.

\begin{prop}\label{prop:convergent}
Let $(L,\mathcal{P},\Gamma)$ be a cascade CCS, with $\Gamma=\{\gamma^n:n\in\mathbb{N}\}$ for a continuous map $\gamma\colon L\to L$, and suppose that the iterates $\gamma^n$
converge pointwise to a map $\gamma_\infty\colon L\to L$.
Then $E(L,\gamma)=\{\gamma^n:n\in\mathbb{N}\}\cup\{\gamma_\infty\}$, and the CCS is
$\NIP_3$.
It is stable if and only if $\gamma_\infty$ is continuous.
\end{prop}

\begin{proof}
A convergent sequence together with its limit is a compact subset of the
Hausdorff space $L^L$, hence closed, and the enveloping semigroup $E(L,\gamma)$ therefore consists of the iterates
and~$\gamma_\infty$.
For each $P\in\mathcal{P}$, the map $P\circ\gamma_\infty$ is the pointwise limit of the
continuous functions $P\circ\gamma^n$ and is therefore of Baire class~$1$;
by Proposition~\ref{prop:cascade-tame}, the CCS is $\NIP$.
The space $E(L,\gamma)$ is a countable compact Hausdorff space, hence
metrizable, so the CCS is $\NIP_3$.
Its members other than $\gamma_\infty$ are continuous, which gives the statement about
stability.
\end{proof}

\subsection{Stable and metrizable}
\label{sec:ex-NIP3}
Our descent through the hierarchy starts at its top two levels, stable and $\NIP_3$, at which the deep computations form a metrizable space.
At these two levels we place a rotation of the circle, Newton’s method for square roots, the power method and the logistic map for $0<r\le3$.
By Proposition~\ref{prop:convergent}, every convergent method lies at one of these two levels.
Convergence, however, is not necessary for stability, as the first example shows.

\emph{A stable example: rotation.}
Let the space of states be the circle $L=\mathbb{T}=\mathbb{R}/\mathbb{Z}$, with the predicates $\cos2\pi x$ and $\sin2\pi x$, and let the step be the rotation $\gamma(x)=x+\alpha$ by an irrational number~$\alpha$.
Since a pointwise limit of rotations is a rotation and the points $n\alpha$ with $n\in\mathbb{N}$ are dense in~$\mathbb{T}$, the deep computations are all the rotations $x\mapsto x+c$ with $c\in\mathbb{T}$.
These rotations are continuous, so the rotation CCS is stable, and its space of deep computations is the circle group, a metrizable compactum.

\emph{Square roots: $\NIP_3$ but not stable.}
For $a>0$, the iterates of the Newton map $N_p$ of the polynomial $p(x)=x^2-a$ on the real projective line converge pointwise to the map $\mathfrak f$ defined in Section~\ref{exa:square-roots-compute}.
Since this limit is discontinuous at $0$ and~$\infty$, Proposition~\ref{prop:convergent} shows that the square-root CCS is $\NIP_3$ but not
stable.

\emph{The power method: $\NIP_3$ but not stable.}
The power method computes a dominant eigenvector of a matrix by repeated multiplication and normalization.
Let $d\ge2$ be an integer, and let $A=\operatorname{diag}(\lambda_1,\dots,\lambda_d)$ be a diagonal matrix with $\lambda_1>\lambda_2>\cdots>\lambda_d>0$.
The space of states is the real projective space $L=\RP^{d-1}$ of lines through the origin in~$\mathbb{R}^d$, in which $[v]$ denotes the line spanned by a nonzero vector~$v$.
The family $\mathcal{P}$ consists of the predicates $P_{ij}([v])=v_iv_j/\|v\|^2$ with $1\le i,j\le d$, which are the entries of the orthogonal projection onto the line~$[v]$.
The step of the algorithm is $\gamma([v])=[Av]$.
For $[v]\in\RP^{d-1}$, let $i(v)$ be the least index $i$ with $v_i\ne0$, and put $\ell([v])=[e_{i(v)}]$, where $e_1,\dots,e_d$ is the standard basis of~$\mathbb{R}^d$.

\begin{prop}\label{prop:power}
The deep computations of the power method are
$E(L,\gamma)=\{\gamma^n:n\in\mathbb{N}\}\cup\{\ell\}$.
Hence the power method is $\NIP_3$.
Its deep computation $\ell$, the only one that is not an iterate, is of Baire class~$1$ and discontinuous, so the power method is not stable.
\end{prop}

\begin{proof}
For $i=i(v)$, we have $\lambda_i^{-n}A^nv=v_ie_i+\sum_{j>i}(\lambda_j/\lambda_i)^nv_je_j\to v_ie_i$,
so $\gamma^n\to\ell$ pointwise, and Proposition~\ref{prop:convergent} applies.
The map $\ell$ is discontinuous at every $[v]$ with $i(v)>1$, since $[v+te_1]\to[v]$ as $t\to0$
while $\ell([v+te_1])=[e_1]\ne\ell([v])$.
\end{proof}

The deep computation~$\ell$ sends each input to the eigenline of the
largest eigenvalue whose eigenvector has a nonzero component in the input.
Being discontinuous on the space of types~$L$, the map~$\ell$ is
not computable in the sense of Section~\ref{sec:CCS}.
This deep computation is nevertheless the limit of a sequence of
computations, and the closure in which it lives is countable and metrizable.

\begin{figure}[htbp]
  \centering
  \includegraphics[width=\textwidth]{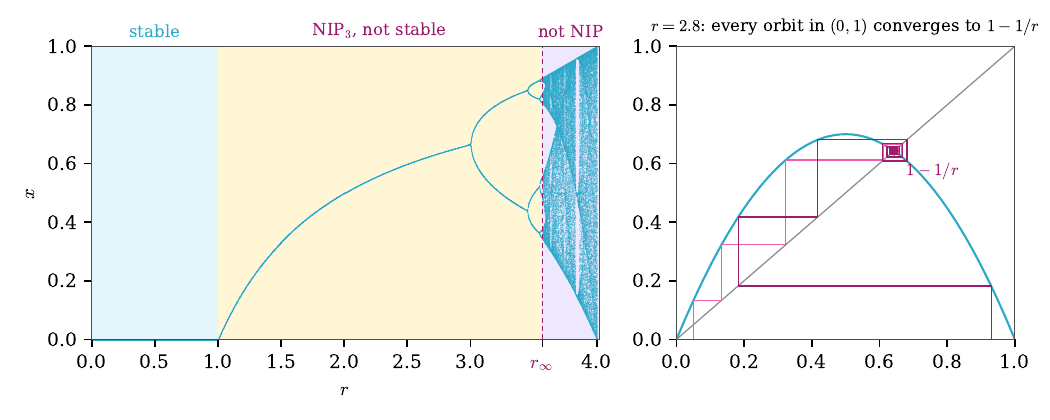}
  \caption{The logistic family $f_r(x)=rx(1-x)$.
    Left: the long-run values of the orbit of~$0.3$ under~$f_r$, for each~$r$.
    The logistic CCS is stable for $r\le1$.
    It is $\NIP_3$ but not stable for $1<r\le r_\infty$, and it is not $\NIP$ for $r_\infty<r\le4$ (see Propositions~\ref{prop:logistic} and~\ref{prop:logistic-all} and Remark~\ref{rem:logistic-feigenbaum}).
    Here $r_\infty\approx3.5699$ is the Feigenbaum parameter, at which the cascade of period doublings ends.
    Right: for $r=2.8$, two orbits converge to the fixed point $1-1/r$.
    The deep computation sends $0$ and~$1$ to~$0$ and every other state to $1-1/r$.}
  \label{fig:logistic}
\end{figure}

\emph{Both levels: the logistic map.}
The logistic map $f_r(x)=rx(1-x)$, with $0<r\le4$, maps $[0,1]$ into itself.
Its CCS has states $L=[0,1]$, the single predicate $x\mapsto x$, and the
semigroup of iterates of~$f_r$ (see Figure~\ref{fig:logistic}).
For $0<r\le1$, let $\ell_r$ be the constant map~$0$; for $1<r\le3$, let
$\ell_r$ be the map with $\ell_r(x)=1-1/r$ for $0<x<1$ and $\ell_r(0)=\ell_r(1)=0$.

\begin{prop}\label{prop:logistic}
For $0<r\le3$, the iterates of $f_r$ converge pointwise to~$\ell_r$.
Consequently the logistic CCS is stable for $0<r\le1$, and $\NIP_3$ but not
stable for $1<r\le3$.
\end{prop}

\begin{proof}
Write $f=f_r$ and $x^{*}=1-1/r$, and note that $f(x)/x=r(1-x)$ for $x>0$.
If $r\le1$, then $f(x)<x$ for $0<x\le1$, so every orbit decreases to a fixed
point of~$f$, and $0$ is the only fixed point in $[0,1]$.

Suppose that $1<r\le2$, in which case $0<x^{*}\le\frac12$.
The inequality $f(x)>x$ holds for $0<x<x^{*}$, and the reverse inequality for $x^{*}<x\le1$.
For $0<x<1$ we have $f(x)\in(0,\frac r4]\subseteq(0,\frac12]$, and the function $f$, which fixes~$x^{*}$, is
increasing on $[0,\frac12]$.
Hence the orbit of $f(x)$ is increasing and bounded by $x^{*}$ if $f(x)<x^{*}$, and
decreasing and bounded below by $x^{*}$ if $f(x)>x^{*}$.
In both cases, as well as trivially when $f(x)=x^{*}$, the orbit of $f(x)$ converges to a fixed point of $f$
in $(0,\frac12]$, namely to~$x^{*}$.

Suppose that $2<r\le3$, and put $c=\frac r4$, the maximum of $f$, and $b=f(c)$.
We have $x^{*}<c$, since $c-x^{*}=(r-2)^2/(4r)$, and therefore $b<f(x^{*})=x^{*}$, since $f$ is strictly decreasing on $[\frac12,1]$.
We also have $b\ge1/r$, because $b-1/r=(r^3(4-r)-16)/(16r)$ and the polynomial $r^3(4-r)$ increases from $16$ on $[2,3]$.
The concave function $f$ takes the value $x^{*}$ at $1/r$ and at $x^{*}$, so $f\ge x^{*}$ on $[1/r,x^{*}]$.
Hence the function $f$ maps $J_-=[b,x^{*}]$ into $J_+=[x^{*},c]$, and, being decreasing on $J_+$, it maps $J_+$ onto~$J_-$.
Fix $x$ with $0<x<1$; the orbit of $f(x)$ then lies in $(0,c]$, and it meets $J_-\cup J_+=[b,c]$.
Indeed, if the orbit did not meet $[b,c]$, all its points would lie in $(0,b)$, where $f(y)>y$ because $b<x^{*}$.
The orbit would then increase to a fixed point of $f$ in $(0,b]$.
No such fixed point exists, since the only positive fixed point of $f$ is $x^{*}$, which exceeds~$b$.
Since $f$ maps $J_-\cup J_+$ into itself, it remains to be shown that the orbit of every $y\in J_-\cup J_+$ converges to~$x^{*}$.
The points of period~$2$ of $f$ that are not fixed are the roots of $r^2y^2-r(r+1)y+(r+1)$, whose discriminant $r^2(r+1)(r-3)$ is negative for $r<3$; for $r=3$ the only root is~$x^{*}$.
Hence the function $f(f(y))-y$ vanishes on $[0,1]$ only at $0$ and $x^{*}$, and it is negative on $(x^{*},1]$, since it is negative at~$1$.
For $y\in J_+$ with $y>x^{*}$, the point $f(f(y))$ lies in $J_+$, so $x^{*}\le f(f(y))<y$.
The sequence $(f^{2n}(y))$ therefore decreases to a fixed point of $f\circ f$ in $[x^{*},y]$, that is, to $x^{*}$, whence $f^{2n+1}(y)\to f(x^{*})=x^{*}$.
Every $y\in J_-$ satisfies $f(y)\in J_+$, so its orbit converges to $x^{*}$ as well.

Finally, the equalities $f(0)=f(1)=0$ account for the endpoints, and in all cases the iterates therefore converge pointwise to~$\ell_r$.
The map $\ell_r$ is continuous for $r\le1$, and discontinuous at $0$ for
$r>1$, so Proposition~\ref{prop:convergent} gives the remaining statements.
\end{proof}

We return to the logistic family for $3<r\le4$ at the end of the section, where we complete its classification (see Proposition~\ref{prop:logistic-all} and Remark~\ref{rem:logistic-feigenbaum}).

\subsection{Hereditarily separable but not metrizable}
\label{sec:ex-NIP2}
Bresenham’s line algorithm in the plane, which we analyze in this subsection, lies at the level $\NIP_2$ but not $\NIP_3$, immediately below the metrizable level.
Below that level, no countable set of tests separates the deep
computations (see Proposition~\ref{prop:levels-tests}), a test being a pair of an input and a feature (see Section~\ref{sec:trichotomy}).
Rounding gives rise to spaces of deep computations of this kind, since an
algorithm that rounds at each step takes a new discrete decision at each step.

\emph{Drawing a line.}
Fix an irrational number $\alpha$ with $0<\alpha<1$.
To draw the line $y=\alpha x+\rho$ on a screen, Bresenham's algorithm~\cite{Bresenham1965} moves one pixel to the right at each step and decides whether to move one pixel up.
The $n$-th decision is $b_n(\rho)=\lfloor(n+1)\alpha+\rho\rfloor-\lfloor n\alpha+\rho\rfloor=\Brack{\{n\alpha+\rho\}\in[1-\alpha,1)}$, where $\{\cdot\}$ denotes the fractional part, and the state of the algorithm is the error term $\rho$ modulo~$1$.
The step $\rho\mapsto\rho+\alpha$ is a rotation of the circle $\mathbb{T}=\mathbb{R}/\mathbb{Z}$, but the decision bit is not a continuous function of~$\rho$.
The standard remedy~\cite{MorseHedlund1940,Lothaire2002} is to take the decision sequences themselves as states.

Let $\mathcal{O}=\{-k\alpha:k\in\mathbb{Z}\}\subseteq\mathbb{T}$, a countable dense set
containing both endpoints $0$ and $1-\alpha$ of the decision interval.
The \emph{split circle} $\mathbb{T}_{\mathcal{O}}$ is obtained from $\mathbb{T}$ by replacing
each $t\in \mathcal{O}$ by two points $t^-<t^+$, with the circular order topology; it is
a compact metrizable space, homeomorphic to the Cantor set.
The rotation by $\alpha$ preserves $\mathcal{O}$ and lifts to a homeomorphism $R$ of
$\mathbb{T}_{\mathcal{O}}$ that preserves the labels $\pm$.
The decision interval becomes clopen in $\mathbb{T}_{\mathcal{O}}$, and the coding map
$\mathbb{T}_{\mathcal{O}}\to2^{\mathbb{Z}}$, $y\mapsto(\Brack{R^ny\in[1-\alpha,1)})_{n\in\mathbb{Z}}$,
is a homeomorphism onto the \emph{Sturmian shift} $X_\alpha$ that conjugates $R$
to the shift~$\sigma$.
The \emph{Bresenham CCS} has states $L=X_\alpha$, the coordinate predicates $P_j(y)=y_j$ for $j\in\mathbb{Z}$, and
the semigroup $\Gamma=\{\sigma^n:n\in\mathbb{N}\}$.

A point of $\mathbb{T}_{\mathcal{O}}$ lies \emph{over} $t\in\mathbb{T}$ if the natural projection $\mathbb{T}_{\mathcal{O}}\to\mathbb{T}$ maps it to~$t$.
For $c\in\mathbb{T}$, let $R_c^+$ (respectively $R_c^-$) be the map of $\mathbb{T}_{\mathcal{O}}$ that sends a point over $t$ to $(t+c)^+$ (respectively $(t+c)^-$) when $t+c\in \mathcal{O}$, and to $t+c$ otherwise.
Through the coding map, we regard $R_c^+$ and $R_c^-$ as maps of~$X_\alpha$.
The \emph{double arrow circle} $\mathbb{D}$ is $\mathbb{T}\times\{-,+\}$ with the lexicographic circular order and its order topology, and we also write $c^-$ and $c^+$ for the points $(c,-)$ and $(c,+)$ of~$\mathbb{D}$.
The double arrow circle is compact, first countable, hereditarily separable and not metrizable~\cite[3.10.C]{Engelking:1989}.

\begin{thm}\label{thm:Bresenham}
The Bresenham CCS is $\NIP$.
Its deep computations are
\[
E(X_\alpha,\sigma)=\{\sigma^n:n\in\mathbb{N}\}\cup\{R_c^+,R_c^-:c\in\mathbb{T}\},
\]
where each $\sigma^n$ is isolated and $c^\pm\mapsto R_c^\pm$ is a homeomorphism
of $\mathbb{D}$ onto the remaining part.
In particular, the Bresenham CCS is $\NIP_2$ but not $\NIP_3$.
\end{thm}

\begin{proof}
To prove that the CCS is $\NIP$, note that the predicates $P_j\circ\sigma^n$ are $\{0,1\}$-valued.
Moreover, the $j$-th coordinate of $\sigma^n y$ is $1$ exactly when $R^{n+j}y$ lies in the clopen arc $[1-\alpha,1)$ of $\mathbb{T}_{\mathcal{O}}$, that is, when $y$ lies in an arc of $\mathbb{T}_{\mathcal{O}}$.
A family of arcs of a circularly ordered set has VC dimension at most~$3$.
By the remark at the end of Section~\ref{sec:NIP-poly-exp}, the characteristic functions of the members of such a family therefore have the $\NIP$.

For the deep computations, let $\gamma$ be the limit of $\sigma^{n_i}$ along an
ultrafilter, and let $c=\lim n_i\alpha$ in $\mathbb{T}$.
If $n_i$ is eventually constant, then $\gamma$ is a power of $\sigma$.
Otherwise, since the points $n\alpha$ are distinct, the ultrafilter decides
whether $n_i\alpha$ approaches $c$ from the right or from the left.
In the first case, for every point $y$ over $t$, the points $R^{n_i}y$ over
$t+n_i\alpha$ approach $t+c$ strictly from the right.
These points converge in $\mathbb{T}_{\mathcal{O}}$ to $(t+c)^+$ if $t+c\in \mathcal{O}$ and to $t+c$ otherwise, whence $\gamma=R_c^+$;
the second case gives $\gamma=R_c^-$.
Conversely, for every $N$, every $c$ is approached from both sides by the points $n\alpha$
with $n\ge N$;
hence every $R_c^\pm$ is a deep computation.
The maps $R_c^+$ and $R_c^-$ are distinct, since they differ at the points over $d-c$ with $d\in \mathcal{O}$.
The same argument shows that $c_i^{\varepsilon_i}\to c^+$ in $\mathbb{D}$ implies
$R_{c_i}^{\varepsilon_i}\to R_c^+$ pointwise, and symmetrically for $c^-$.
The map $c^\pm\mapsto R_c^\pm$ is therefore a continuous bijection from the compact space
$\mathbb{D}$ onto a Hausdorff space, hence a homeomorphism.
Finally, fix $d\in \mathcal{O}$ and consider the two points
$(d-n\alpha)^\pm$ of $\mathbb{T}_{\mathcal{O}}$, which $\sigma^n$ sends to $d^-$ and $d^+$.
Every map $R_c^\pm$ sends them to the same point, and so does every pointwise
limit of $\sigma^{m_i}$ with $m_i\to\infty$, by the argument above.
Hence no net of such maps, and no net of maps $\sigma^m$ with $m\ne n$ and $m$ bounded, converges to $\sigma^n$, so $\sigma^n$ is isolated.
Thus $E(X_\alpha,\sigma)$ is the union of a countable set and a copy of
$\mathbb{D}$; it is hereditarily separable and not metrizable.
\end{proof}

\begin{figure}[htbp]
  \centering
  \includegraphics[width=\textwidth]{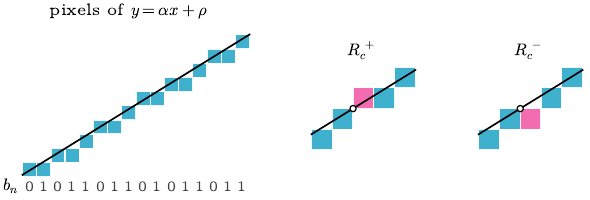}
  \caption{Left: the pixels drawn by Bresenham’s algorithm for $\alpha=(\sqrt5-1)/2$, with the decision bits $b_n$.
    Right: where the line passes exactly through a pixel corner, the deep computations $R_c^+$ and $R_c^-$ break the tie in opposite directions; they agree on every other pixel.}
  \label{fig:bresenham}
\end{figure}

The two deep computations $R_c^+$ and $R_c^-$ are the two ways of drawing the
line infinitely far away, where the accumulated error has reached~$c$.
The two maps agree except at the pixels where the line passes exactly through a corner,
and at those pixels they apply the two opposite rounding conventions.
Figure~\ref{fig:bresenham} shows the two conventions.

The space $E(X_\alpha,\sigma)$ is the closure of a countable set, and it is therefore a
separable Rosenthal compactum.
Since this compactum is hereditarily separable and not metrizable, it contains, by Todorčević's trichotomy
(see Theorem~\ref{thm:Todorcevic-trichotomy}), a copy of the split
Cantor set $S(2^{\mathbb{N}})$.
For the action of $\mathbb{Z}$ generated by the shift, Glasner and
Megrelishvili proved that Sturmian systems are tame$_2$ but not hereditarily
non-sensitive~\cite[Corollary~8.3]{glasner2022tame}.

\subsection{First countable but not hereditarily separable}
\label{sec:ex-NIP1}
When two rounding algorithms run in parallel, the computation descends a further level.
We establish this descent here in the case of Bresenham’s algorithm for lines in three-dimensional space, which lies at the level $\NIP_1$ but not $\NIP_2$.
At this level each deep computation is still separated from all the others by countably many tests of its own (see Proposition~\ref{prop:levels-tests}).
The same proposition shows, however, that some finitely separated set of deep computations has the cardinality of the continuum.

\emph{Drawing a line in space.}
Let $\alpha$ and $\beta$ be irrational numbers in $(0,1)$ such that $1,\alpha,\beta$
are linearly independent over $\mathbb{Q}$.
To draw a line in $\mathbb{Z}^3$ with direction $(1,\alpha,\beta)$, the
three-dimensional Bresenham algorithm runs the planar algorithm in two
coordinates at once.
Its CCS has states $L=X_\alpha\times X_\beta$, the coordinate predicates of
both factors, and the semigroup generated by $\sigma\times\sigma$.

The analysis of this algorithm uses a property of the double arrow circle~$\mathbb{D}$ of Section~\ref{sec:ex-NIP2} that keeps every compact metrizable subset of $\mathbb{D}\times\mathbb{D}$ countable.
Recall that a \emph{network} of a space is a family of subsets such that every neighborhood of a point contains a member of the family that contains the point.
In particular, every base is a network.

\begin{lem}\label{lem:double-arrow-metrizable}
Every subspace of the double arrow circle~$\mathbb{D}$ that has a countable network, in particular every metrizable subspace, is countable.
\end{lem}

\begin{proof}
A metrizable subspace of~$\mathbb{D}$ is separable, since $\mathbb{D}$ is hereditarily separable, and it therefore has a countable base, which is a countable network.
Let $M\subseteq\mathbb{D}$ have a countable network $\{N_k:k\in\mathbb{N}\}$.
For each $c\in\mathbb{T}$ with $c^+\in M$, the open arc $[c^+,(c+\frac12)^-)$ is a neighborhood of $c^+$ in~$\mathbb{D}$, and we may choose $k(c)\in\mathbb{N}$ with $c^+\in N_{k(c)}\subseteq[c^+,(c+\frac12)^-)$.
If $k(c)=k(c')$, then the point $c'^+$ lies in $[c^+,(c+\frac12)^-)$ and the point $c^+$ in $[c'^+,(c'+\frac12)^-)$, whence both $c'-c$ and $c-c'$ lie in $[0,\frac12)$ modulo~$1$, which forces $c=c'$.
Thus the map $c\mapsto k(c)$ is injective, and the set $M$ contains only countably many points~$c^+$.
The same argument, applied to the arcs $((c-\frac12)^+,c^-]$, shows that $M$ contains only countably many points~$c^-$.
\end{proof}

\begin{cor}\label{cor:DxD-compact-metrizable}
Every compact metrizable subspace of $\mathbb{D}\times\mathbb{D}$ is countable.
\end{cor}

\begin{proof}
Let $C\subseteq\mathbb{D}\times\mathbb{D}$ be compact and metrizable, and let $\mathcal{N}$ be a countable base of~$C$.
Since the two projections are continuous, the images of the members of~$\mathcal{N}$ under each projection form a countable network of the image of~$C$.
By Lemma~\ref{lem:double-arrow-metrizable}, both images are countable, and so is~$C$, which is contained in their product.
\end{proof}

For $c'\in\mathbb{T}$, we also write $R_{c'}^+$ and $R_{c'}^-$ for the maps of $X_\beta$ defined as in Section~\ref{sec:ex-NIP2}, with $\beta$ in place of~$\alpha$.

\begin{thm}\label{thm:Bresenham3}
The three-dimensional Bresenham CCS is $\NIP_1$ but not $\NIP_2$.
Its space of deep computations $E=E(X_\alpha\times X_\beta,\sigma\times\sigma)$ is the union of the set of iterates $(\sigma\times\sigma)^n$, $n\in\mathbb{N}$, and of the closed subspace
\[
\{(R_c^\varepsilon,R_{c'}^{\varepsilon'}):c,c'\in\mathbb{T},\ \varepsilon,\varepsilon'\in\{+,-\}\},
\]
which is homeomorphic to $\mathbb{D}\times\mathbb{D}$.
Each of these iterates is isolated in~$E$.
The space $E$ contains a copy of $\hat D(S(2^{\mathbb{N}}))$, and it contains no copy of $D(2^{\mathbb{N}})$, of $\hat D(2^{\mathbb{N}})$ or of~$A(2^{\mathbb{N}})$.
\end{thm}

\begin{proof}
Each predicate is a coordinate of one factor, so the CCS is $\NIP$ by Theorem~\ref{thm:Bresenham}.
Every deep computation is the limit $\xi$ of the iterates $(\sigma\times\sigma)^n$ along an ultrafilter $\mathcal{U}$ on~$\mathbb{N}$.
The coordinates of~$\xi$ are then the limits along~$\mathcal{U}$ of the iterates $\sigma^n$ on $X_\alpha$ and on~$X_\beta$.
If $\mathcal{U}$ is principal, the deep computation $\xi$ is an iterate;
otherwise the proof of Theorem~\ref{thm:Bresenham} shows that $\xi$ is a pair $(R_c^\varepsilon,R_{c'}^{\varepsilon'})$.
By Kronecker's theorem, the points $(n\alpha,n\beta)$ with $n\ge N$ are dense in $\mathbb{T}^2$ for every $N$.
Every point $(c,c')$ is therefore approached by these points from inside each of the four open quadrants at $(c,c')$.
As in the proof of Theorem~\ref{thm:Bresenham}, it follows that every pair $(R_c^\varepsilon,R_{c'}^{\varepsilon'})$ with $\varepsilon,\varepsilon'\in\{+,-\}$ is a deep computation.
By Theorem~\ref{thm:Bresenham}, applied to each factor, the set of these pairs is homeomorphic to $\mathbb{D}\times\mathbb{D}$, and it is therefore compact and closed in~$E$.
Each iterate $(\sigma\times\sigma)^n$ is isolated in~$E$, since $\sigma^n$ is isolated in $E(X_\alpha,\sigma)$ and is the first coordinate of no other deep computation.
The space $E$ is a closed subspace of $E(X_\alpha,\sigma)\times E(X_\beta,\sigma)$, which is first countable as a product of two first countable spaces.
Hence the CCS is $\NIP_1$, and no copy of the space $A(2^{\mathbb{N}})$, which is not first countable, lies in~$E$.
Consider the pairs $(R_c^+,R_{-c}^+)$, for $c\in\mathbb{T}$.
For $0<\delta<\frac12$, the set $[c^+,(c+\delta)^-)\times[(-c)^+,(-c+\delta)^-)$ is a neighborhood of $(c^+,(-c)^+)$ in $\mathbb{D}\times\mathbb{D}$ that contains $(s^+,(-s)^+)$ only if $c\le s<c+\delta$ and $-c\le -s<-c+\delta$, that is, only if $s=c$.
These pairs therefore form an uncountable discrete subspace, and the space of deep computations is not hereditarily separable;
thus the CCS is not $\NIP_2$.
A compact metrizable subspace of $E$ meets the copy of $\mathbb{D}\times\mathbb{D}$ in a compact metrizable set, which is countable by Corollary~\ref{cor:DxD-compact-metrizable}, and the subspace itself is therefore countable.
Since the non-isolated points of $D(2^{\mathbb{N}})$ form a copy of the Cantor set, neither $D(2^{\mathbb{N}})$ nor the compactum $\hat D(2^{\mathbb{N}})$ that contains it embeds into~$E$.
Finally, by Theorem~\ref{thm:nip-baire1def} the space $E$, being the closure of a countable set, is a separable Rosenthal compactum that is first countable but not hereditarily separable.
By Theorem~\ref{thm:Todorcevic-trichotomy}, this compactum contains a copy of $D(2^{\mathbb{N}})$ or of $\hat D(S(2^{\mathbb{N}}))$, necessarily of the latter.
\end{proof}

The uncountable discrete set consists of the deep computations that break ties upward in both coordinates at points where the two accumulated errors are opposite.
Each of the two planar algorithms is $\NIP_2$, and drawing in two directions at once leaves the class.
In just the same way, the Sorgenfrey plane shows that a product of two hereditarily separable spaces need not be hereditarily separable.
Glasner and Megrelishvili show similarly that a product of two tame$_2$ systems need not be tame$_2$~\cite[Proposition~4.1]{glasner2022tame}.%
\footnote{The double arrow space is the set of the points $p^-<p^+$, for $p\in[0,1]$, with its order topology.
  Glasner and Megrelishvili~\cite[Remark~4.2]{glasner2022tame} assert that the square of this space contains a copy of $D(2^{\mathbb{N}})$, namely the set of the points $(p^+,(1-p)^+)$ and $(p^-,(1-p)^+)$ with $p$ in the Cantor set.
  The set of these points is not closed, and hence not compact.
  For instance, if $p$ is the limit of a decreasing sequence $(x_k)_{k\in\mathbb{N}}$ of points of the Cantor set, then the points $(x_k^-,(1-x_k)^+)$ of this set converge to $(p^+,(1-p)^-)$, a point outside it.
  Moreover, the square of the double arrow space contains no copy of $D(2^{\mathbb{N}})$ at all.
  Indeed, with intervals in place of arcs, the arguments of Lemma~\ref{lem:double-arrow-metrizable} and Corollary~\ref{cor:DxD-compact-metrizable} apply to the double arrow space, and the non-isolated points of $D(2^{\mathbb{N}})$ form an uncountable compact metrizable set.
  Proposition~4.1 of Glasner and Megrelishvili does not depend on their Remark~4.2.}
A second example at this level arises from the Euclidean algorithm (see Proposition~\ref{prop:euclid-duplicate}).
This example is a countable family of elements of the group generated by the steps of the Euclidean algorithm, and its deep computations include a copy of $D(2^{\mathbb{N}})$.
Whether a single iterated map can realize the level with such a copy is Question~\ref{q:cascade-duplicate}.

\begin{rem}\label{rem:Bresenham3-explicit}
A copy of $\hat D(S(2^{\mathbb{N}}))$ in~$E$ whose dense set of isolated points consists of iterates can be exhibited explicitly.
For $t\in2^{<\mathbb{N}}$, let $x_t$ be the right endpoint of the closed interval of length $3^{-|t|}$ indexed by~$t$ in the construction of the middle-thirds Cantor set, and put $y_t=\frac32\,3^{-|t|}-x_t$.
By Kronecker's theorem, there are pairwise distinct $n_t\in\mathbb{N}$ such that, modulo~$1$, the point $n_t\alpha$ lies in $(x_t,x_t+3^{-|t|-1})$ and the point $n_t\beta$ in $(y_t,y_t+3^{-|t|-1})$.
For $b\in2^{\mathbb{N}}$, let $x_b$ be the point of the Cantor set coded by~$b$, the limit of the points $x_{b\restr n}$.
As $t$ runs through the initial segments of~$b$, the iterates $(\sigma\times\sigma)^{n_t}$ converge to $(R_{x_b}^+,R_{-x_b}^+)$.
Along words $t$ that converge to~$b$, are not initial segments of~$b$ and satisfy $t^\frown0^\infty<b$, the iterates converge to $(R_{x_b}^-,R_{-x_b}^+)$.
Finally, along words that converge to~$b$ and satisfy $t^\frown0^\infty>b$, the iterates converge to $(R_{x_b}^+,R_{-x_b}^-)$.
The three convergences just described are the convergence patterns of the family $D_7$ of Section~\ref{sec:heptachotomy}.
Since the iterates are isolated in~$E$, a lemma of Argyros, Dodos and Kanellopoulos~\cite[Lemma~19]{argyros2008rosenthal} shows that the family $\{(\sigma\times\sigma)^{n_t}:t\in2^{<\mathbb{N}}\}$ is equivalent to~$D_7$ and therefore has a closure homeomorphic to $\hat D(S(2^{\mathbb{N}}))$.
These three patterns are drawn in the left panel of Figure~\ref{fig:trichotomy-ii};
its right panel shows the copy of $\hat D(2^{\mathbb{N}})$ that maps built from the Euclidean algorithm produce in Section~\ref{sec:ex-NIP}.
\end{rem}

\subsection{NIP but not first countable: the Euclidean algorithm}
\label{sec:ex-NIP}
The last tame level, $\NIP$ but not $\NIP_1$, is reached by the Euclidean algorithm, acting on the real projective line through the group generated by its steps.
At this level some deep computation is not determined by countably
many tests of its own, that is, by countably many pairs of an input and a feature (see Proposition~\ref{prop:levels-tests}).
A constant deep computation has this property whenever uncountably many exact equality tests, maps each of which agrees with the constant map at every input but one, are also deep computations.
Indeed, countably many tests involve only countably many inputs and therefore cannot distinguish the constant map from an exact equality test whose exceptional input is not among those inputs.
The Euclidean algorithm gives rise to equality tests of precisely this kind among its deep computations.

Let $L=\RP^1=\mathbb{R}\cup\{\infty\}$, with the predicates $P_1,P_3$ of
Section~\ref{exa:square-roots-compute}, and let $\Gamma=\mathrm{PGL}(2,\mathbb{Z})$ be the
group of Möbius maps $x\mapsto(ax+b)/(cx+d)$ with integer coefficients and
$ad-bc=\pm1$.
It is generated by the steps $x\mapsto x+1$ and $x\mapsto1/x$ of the
Euclidean algorithm, and the maps that the continued fraction algorithm
composes are elements of~$\Gamma$.
We call $(\RP^1,(P_1,P_3),\Gamma)$ the \emph{Euclidean CCS}.
It is a dynamical CCS whose semigroup is a group, and its deep computations
form $E(\RP^1,\Gamma)$.

For an irrational number $\eta=[a_0;a_1,a_2,\dots]$, let $p_k/q_k$ be its
convergents, put $p_{-1}=1$ and $q_{-1}=0$, and let
$\eta_{k+1}=[a_{k+1};a_{k+2},\dots]$ be its complete quotients.
Let $M_k$ be the Möbius map $t\mapsto(p_kt+p_{k-1})/(q_kt+q_{k-1})$, which belongs
to $\Gamma$ because $p_kq_{k-1}-p_{k-1}q_k=(-1)^{k+1}$.
Writing $A_j(t)=a_j+1/t$, we have $M_k=A_0\circ A_1\circ\dots\circ A_k$, and hence the inverse map is
$M_k^{-1}=A_k^{-1}\circ\dots\circ A_0^{-1}$, with $A_j^{-1}(t)=1/(t-a_j)$.
Thus $M_k^{-1}$ is the map that the first $k+1$ steps of the continued fraction
algorithm for $\eta$ apply.
In these terms, the classical identity $\eta=M_k(\eta_{k+1})$ states that these steps send $\eta$ to its complete quotient~$\eta_{k+1}$.

\begin{lem}\label{lem:euclid-limit}
Let $\eta$ be an irrational number, and let $x\in\RP^1$ be distinct from~$\eta$.
Then $M_k^{-1}(x)+q_{k-1}/q_k\to0$ as $k\to\infty$, while $M_k^{-1}(\eta)=\eta_{k+1}>1$.
\end{lem}

\begin{proof}
Solving $M_k(t)=x$ gives $M_k^{-1}(\infty)=-q_{k-1}/q_k$ and, for real $x$,
\[
M_k^{-1}(x)=\frac{q_{k-1}x-p_{k-1}}{p_k-q_kx}=-\frac{q_{k-1}}{q_k}\cdot\frac{x-p_{k-1}/q_{k-1}}{x-p_k/q_k}
\]
for $k\ge1$ (for such $k$ we have $q_{k-1}\ge1$), where the value is $\infty$ if $x=p_k/q_k$.
The convergents of an irrational number are pairwise distinct, so $x$ equals
$p_k/q_k$ for at most one $k$; for all larger $k$ the formula gives a real number.
Both convergents tend to $\eta\ne x$, so the last fraction tends to $1$, and
$0<q_{k-1}/q_k\le1$ for $k\ge1$; hence $M_k^{-1}(x)+q_{k-1}/q_k\to0$.
The identity $M_k^{-1}(\eta)=\eta_{k+1}$ is the classical one, and
$\eta_{k+1}>a_{k+1}\ge1$ because $\eta_{k+1}$ is irrational.
\end{proof}

Let $\varphi=(1+\sqrt5)/2$ be the golden ratio.
Let $G$ be the set of irrational numbers $\eta=[0;a_1,a_2,\dots]$ whose sequence
of partial quotients is a concatenation $1^{1}d_11^{2}d_21^{3}d_3\cdots$, where
$1^n$ denotes $n$ consecutive ones and each $d_n\in\{1,2\}$.
The block $1^n$ occupies the indices $s_n,\dots,s_n+n-1$ with $s_n=n(n+1)/2$.
Distinct choices of $(d_n)$ give distinct numbers, by the uniqueness of
continued fraction expansions, so $G$ has the cardinality of the continuum.

\begin{prop}\label{prop:euclid}
\begin{enumerate}[(i)]
\item For every $\eta\in G$, the map $\delta_\eta$ that sends $\eta$ to $\varphi$
and every other point of $\RP^1$ to $-1/\varphi$ is a deep computation of the
Euclidean CCS.
\item The constant map $c$ with value $-1/\varphi$ is a deep computation, and it
has no countable base of neighborhoods in $E(\RP^1,\Gamma)$.
\end{enumerate}
Consequently $E(\RP^1,\Gamma)$ is not first countable.
\end{prop}

\begin{proof}
(i) Let $\eta\in G$ be given; for each $n$, let $k_n=s_n+\lfloor n/2\rfloor$ be the
middle index of the block $1^n$, and put $m_n=\lceil n/2\rceil-1$.
Then $a_{k_n}=a_{k_n-1}=\dots=a_{k_n-m_n}=1$ and $a_{k_n+1}=\dots=a_{k_n+m_n}=1$,
while $m_n\to\infty$.
The two quantities that occur in Lemma~\ref{lem:euclid-limit} are continued fractions: $q_{k-1}/q_k=[0;a_k,a_{k-1},\dots,a_1]$ and $\eta_{k+1}=[a_{k+1};a_{k+2},\dots]$.
Suppose that two continued fractions, finite or infinite, $[c_0;c_1,\dots]$ and
$[c_0';c_1',\dots]$ satisfy $c_j=c_j'$ for $j<m$, with $c_1=\dots=c_{m-1}=1$.
Then they share the convergents of index at most $m-1$, and both lie between
the last two of them (or equal one of them).
It follows that the two continued fractions differ by at most $1/(q_{m-2}q_{m-1})\le1/(F_{m-1}F_m)$, whatever their remaining partial quotients are.
Here $F_j$ is the $j$th Fibonacci number, with $F_1=F_2=1$.
The hypothesis holds with $m=m_n$ for the pair $[0;a_{k_n},\dots,a_1]$, $[0;1,1,\dots]$ and for the pair $[a_{k_n+1};a_{k_n+2},\dots]$, $[1;1,1,\dots]$.
Hence $q_{k_n-1}/q_{k_n}\to[0;1,1,\dots]=1/\varphi$ and
$\eta_{k_n+1}\to[1;1,1,\dots]=\varphi$, and Lemma~\ref{lem:euclid-limit} shows that the maps $M_{k_n}^{-1}\in\Gamma$ converge
pointwise to $\delta_\eta$, at $x=\infty$ as well.

(ii) Let $(\eta_j)$ be a sequence of distinct elements of $G$.
For each $x\in\RP^1$, at most one $j$ has $\eta_j=x$, so
$\delta_{\eta_j}(x)=-1/\varphi$ for all but at most one $j$.
Hence $\delta_{\eta_j}\to c$ pointwise, and $c\in E(\RP^1,\Gamma)$ because this set
is closed.
Now let $(U_i)_{i\in\mathbb{N}}$ be open neighborhoods of $c$ in
$E(\RP^1,\Gamma)$.
Each $U_i$ contains a basic neighborhood
$\{f:d(f(x),c(x))<\varepsilon_i\text{ for }x\in F_i\}$ with $F_i\subseteq\RP^1$
finite, $\varepsilon_i>0$ and $d$ a metric on~$\RP^1$.
Since the set $F=\bigcup_iF_i$ is countable, some $\eta\in G$ lies outside~$F$;
the map $\delta_\eta$ then agrees with $c$ on~$F$ and hence lies in $\bigcap_iU_i$, although $\delta_\eta\ne c$.
Thus no countable family of neighborhoods of $c$ has intersection $\{c\}$.
In a Hausdorff space, the intersection of a base of neighborhoods of $c$ is
$\{c\}$, so $c$ has no countable base of neighborhoods.
\end{proof}

\begin{figure}[htbp]
  \centering
  \includegraphics[width=\textwidth]{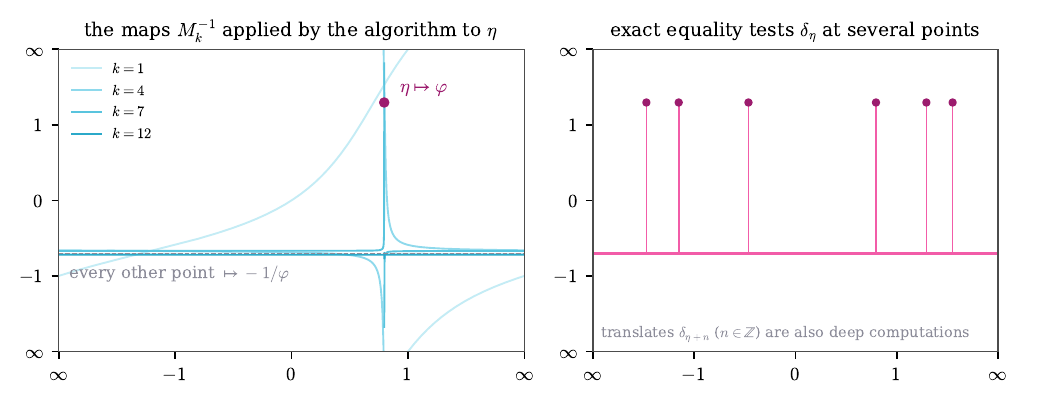}
  \caption{The Euclidean algorithm, drawn on the projective line (each axis runs from $\infty$ through $0$ back to $\infty$). Left: the maps $M_k^{-1}$ that the continued fraction algorithm applies to a number $\eta$ whose partial quotients contain longer and longer runs of ones. As $k$ grows along the middles of the runs, these maps send every point other than $\eta$ near $-1/\varphi$, and $\eta$ itself near $\varphi$. Right: the limits, exact equality tests $\delta_\eta$, at several points. Every neighborhood of the constant map $-1/\varphi$ contains all but finitely many of the uncountably many tests $\delta_\eta$, so the space of deep computations is not first countable (see Proposition~\ref{prop:euclid}).}
  \label{fig:euclid}
\end{figure}

\begin{cor}\label{cor:euclid}
The Euclidean CCS is $\NIP$ but not $\NIP_1$.
Its space of deep computations contains the copy
$\{\delta_\eta:\eta\in G\}\cup\{c\}$ of the Alexandroff compactification
$A(2^{\mathbb{N}})$.
\end{cor}

\begin{proof}
Ellis proved that the action of $\mathrm{GL}(2,\mathbb{R})$ on the projective line is
tame~\cite{Ellis1993} (see~\cite[Example~5.7]{glasner2022tame}).
For every continuous $f\colon\RP^1\to\mathbb{R}$, the family
$\{f\circ\gamma:\gamma\in\Gamma\}$ is contained in the corresponding family for
$\mathrm{GL}(2,\mathbb{R})$, so it has the $\NIP$.
Hence the action of $\Gamma$ is tame, and, by
Proposition~\ref{prop:cascade-tame}, the Euclidean CCS is $\NIP$.
By Proposition~\ref{prop:euclid}, it is not $\NIP_1$.
Each $\delta_\eta$ is isolated in $\{\delta_\eta:\eta\in G\}\cup\{c\}$, since $\delta_\eta$ is the only member of this set that sends $\eta$ near~$\varphi$.
Moreover, as in the proof of Proposition~\ref{prop:euclid}, every neighborhood of $c$ contains all but finitely many of the maps~$\delta_\eta$.
The set $\{\delta_\eta:\eta\in G\}\cup\{c\}$ is therefore the one-point compactification of a discrete set of cardinality $2^{\aleph_0}$, and hence it is homeomorphic to $A(2^{\mathbb{N}})$ (see Section~\ref{sec:five-compacta}).
\end{proof}

The deep computation $\delta_\eta$ is an exact equality test (see Figure~\ref{fig:euclid}): it takes one value if the input is exactly $\eta$ and another value otherwise.
It is the pointwise limit of maps that the continued fraction algorithm applies to $\eta$, when these maps are applied to all inputs.
The mechanism is that of convergence groups.
These are groups of homeomorphisms of a compact space in which every sequence of distinct elements has a subsequence converging to a constant map uniformly on the compact subsets of the complement of some point.
By Lemma~\ref{lem:euclid-limit}, the maps $M_k^{-1}$ draw the images of all
points other than $\eta$ together, near $-q_{k-1}/q_k\in[-1,0)$, while sending $\eta$ itself into $(1,\infty)$.
In the language of convergence groups, every irrational number~$\eta$ is thus a conical limit point of~$\Gamma$: some sequence in~$\Gamma$ sends~$\eta$ to one limit and every other point to another.
The choice of $G$ makes the two limits the same for uncountably many~$\eta$.
Glasner and Megrelishvili use the same device to prove that an infinite
countable uniform convergence group acting on a compact metric space with
infinite limit set is tame but not
tame$_1$~\cite[Proposition~6.1 and Theorem~6.5]{glasner2022tame}.
Their theorem does not cover $\Gamma$, which contains parabolic elements and is
not a uniform convergence group on~$\RP^1$.
Ellis's example, the projective action of the whole group
$\mathrm{GL}(2,\mathbb{R})$, is also tame but not
tame$_1$~\cite[Example~5.7]{glasner2022tame}.
The Euclidean CCS shows that the countable group generated by the two steps of the algorithm already reaches this level.

Iterating a single element of $\Gamma$ cannot produce $\delta_\eta$.
For $g\in\Gamma$, the closure of $\{g^n:n\in\mathbb{Z}\}$ in $(\RP^1)^{\RP^1}$ is
countable, hence metrizable, and it does not contain $\delta_\eta$.
Indeed, if $g$ has finite order, the set of its powers is finite.
If $g$ is parabolic, with fixed point $x_0$, then $g^n$ converges pointwise to the constant map $x_0$ as $n\to\pm\infty$.
If $g$ is hyperbolic, with attracting and repelling fixed points $x_a$ and $x_r$, then $g^n$ converges as $n\to+\infty$ to the map that fixes $x_r$ and sends every other point to~$x_a$.
As $n\to-\infty$, the two fixed points exchange their roles, and $g^n$ converges to the map that fixes $x_a$ and sends every other point to~$x_r$.
If $g$ reverses orientation and has infinite order, then its trace $t$ is
nonzero, whence the square $g^2$ has trace $t^2+2>2$ and is a hyperbolic map with the same fixed points as~$g$.
The odd powers $g^{2m+1}=g\circ g^{2m}$ have the same two limits as the even ones.
In every case the closure consists of the powers of~$g$, which are homeomorphisms, together with at most two limit maps.
Each of these limit maps is constant or fixes the unique point at which that map differs from a constant.
The map $\delta_\eta$ is neither a power of~$g$ nor such a limit map, since it is not injective and
$\delta_\eta(\eta)=\varphi\ne\eta$.

\emph{A copy of the Alexandroff duplicate.}
Among the deep computations of the Euclidean CCS, the device of Proposition~\ref{prop:euclid} also produces a copy of the Alexandroff duplicate $D(2^{\mathbb{N}})$.
With this copy comes a countable set of computations at the level $\NIP_1$ but not $\NIP_2$.
Instead of sending uncountably many irrational numbers to the same pair of values, we now let the value that a deep computation takes away from its exceptional point vary continuously with that point.
Let $\mathcal{C}$ be the set of irrational numbers $[0;\omega_1,\omega_2,\dots]$ with $\omega_n\in\{1,2\}$ for every $n\ge1$, a Cantor set homeomorphic to $\{1,2\}^{\mathbb{N}}$ through continued fraction expansions.
For $\omega=(\omega_1,\omega_2,\dots)\in\{1,2\}^{\mathbb{N}}$ we write $[0;\omega]$ for the point $[0;\omega_1,\omega_2,\dots]$ of~$\mathcal{C}$.
For a finite or infinite sequence $\omega$ over the alphabet $\{1,2\}$ and for each positive integer~$n$ not exceeding its length, let $\Lambda_n(\omega)$ be the block $(\omega_n,\omega_{n-1},\dots,\omega_1,1^n)$.
This block consists of the first $n$ terms of~$\omega$ in reverse order, followed by $n$ ones.
For $\omega\in\{1,2\}^{\mathbb{N}}$, let $\eta(\omega)$ be the irrational number in $(0,1)$ whose partial quotients after the initial~$0$ form the concatenation $\Lambda_1(\omega)\Lambda_2(\omega)\Lambda_3(\omega)\cdots$.
Since the terms $\omega_1,\dots,\omega_n$ can be recovered from $\Lambda_n(\omega)$, the map $\omega\mapsto\eta(\omega)$ is injective.
For a finite word $\mathsf{w}=(\mathsf{w}_1,\dots,\mathsf{w}_k)$ of positive integers, let $M_{\mathsf{w}}$ be the Möbius map $x\mapsto[0;\mathsf{w}_1,\dots,\mathsf{w}_{k-1},\mathsf{w}_k+1/x]$, and let $M_{\emptyset}(x)=1/x$ for the empty word.
For an irrational number $[0;a_1,a_2,\dots]$ in $(0,1)$, the map $M_k$ introduced before Lemma~\ref{lem:euclid-limit} is then $M_{(a_1,\dots,a_k)}$.
Write $\{1,2\}^{<\mathbb{N}}$ for the set of finite words over $\{1,2\}$.
For a word $t=(t_1,\dots,t_m)$ in this set, let $t^\sharp$ be the word $\Lambda_1(t)\cdots\Lambda_{m-1}(t)\,(t_m,t_{m-1},\dots,t_1)$ of length~$m^2$, and put $\gamma_t=M_{t^\sharp}^{-1}\in\Gamma$.
Finally, for $\omega\in\{1,2\}^{\mathbb{N}}$, let $\tau_\omega$ be the map that sends $\eta(\omega)$ to~$\varphi$ and every other point of~$\RP^1$ to~$-[0;\omega]$, and let $c_\omega$ be the constant map with value~$-[0;\omega]$.
For the constant sequence $\omega=(1,1,1,\dots)$, the map $c_\omega$ is the map~$c$ of Proposition~\ref{prop:euclid}.

\begin{prop}\label{prop:euclid-duplicate}
\begin{enumerate}[(i)]
\item
  The closure of the countable set $\Delta_{\mathcal{C}}=\{\gamma_t:t\in\{1,2\}^{<\mathbb{N}}\}$ is $\Delta_{\mathcal{C}}\cup\{\tau_\omega:\omega\in\{1,2\}^{\mathbb{N}}\}\cup\{c_\omega:\omega\in\{1,2\}^{\mathbb{N}}\}$, and it is homeomorphic to $\hat D(2^{\mathbb{N}})$ through a homeomorphism that maps $\{\tau_\omega:\omega\in\{1,2\}^{\mathbb{N}}\}\cup\{c_\omega:\omega\in\{1,2\}^{\mathbb{N}}\}$ onto $D(2^{\mathbb{N}})$.
  In particular, $E(\RP^1,\Gamma)$ contains a copy of $D(2^{\mathbb{N}})$, and the set $\Delta_{\mathcal{C}}$ of computations of the Euclidean CCS is $\NIP_1$ but not $\NIP_2$.
\item
  Every hereditarily separable compact subset of $E(\RP^1,\Gamma)$ is metrizable.
  In particular, $E(\RP^1,\Gamma)$ contains no copy of $S(2^{\mathbb{N}})$ and none of $\hat D(S(2^{\mathbb{N}}))$.
\end{enumerate}
\end{prop}

\begin{proof}
(i) We use two facts from the proofs of Lemma~\ref{lem:euclid-limit} and Proposition~\ref{prop:euclid}.
First, consider two continued fractions, finite or infinite, with positive integer partial quotients that agree in their terms of index less than~$m$.
Both lie between two consecutive convergents of index $m-2$ and $m-1$, and hence differ by at most $1/(F_{m-1}F_m)$.
Second, the formula for $M_k^{-1}$ in the proof of Lemma~\ref{lem:euclid-limit} uses only the terms $a_0,\dots,a_k$:
if $k\ge1$ and $p_i/q_i$ are the convergents of $[a_0;a_1,\dots,a_k]$, the map $M(x)=[a_0;a_1,\dots,a_{k-1},a_k+1/x]$ satisfies
\begin{equation}\label{eq:Mobius-inverse}
M^{-1}(x)=-\frac{q_{k-1}}{q_k}\cdot\frac{x-p_{k-1}/q_{k-1}}{x-p_k/q_k}\quad(x\in\mathbb{R}),
\qquad
M^{-1}(\infty)=-\frac{q_{k-1}}{q_k},
\end{equation}
where $q_{k-1}/q_k=[0;a_k,a_{k-1},\dots,a_1]$.

\emph{Limits along initial segments.}
Let $\omega\in\{1,2\}^{\mathbb{N}}$ be given, and put $\eta=\eta(\omega)=[0;a_1,a_2,\dots]$.
For $m\ge1$, the word $(\omega\restr m)^\sharp$ consists of the first $m^2$ partial quotients of~$\eta$, whence $\gamma_{\omega\restr m}=M_{m^2}^{-1}$.
The terms $a_{m^2},a_{m^2-1},\dots,a_{m^2-m+1}$ are $\omega_1,\omega_2,\dots,\omega_m$, and $a_{m^2+1}=\dots=a_{m^2+m}=1$.
By the first fact, as $m\to\infty$, the quotient $q_{m^2-1}/q_{m^2}=[0;a_{m^2},\dots,a_1]$ tends to~$[0;\omega]$ and the complete quotient $\eta_{m^2+1}$ tends to $[1;1,1,\dots]=\varphi$.
Lemma~\ref{lem:euclid-limit} therefore shows that the maps $\gamma_{\omega\restr m}$ converge pointwise to~$\tau_\omega$.

\emph{Limits along other words.}
Let $(t_j)_{j\in\mathbb{N}}$ be a sequence of words that converges to~$\omega$, none of which is an initial segment of~$\omega$.
For each $j$, let $n_j$ be the length of the longest common initial segment of $t_j$ and~$\omega$, and observe that $n_j\to\infty$ because $t_j\to\omega$.
To show that $\gamma_{t_j}\to c_\omega$ pointwise, we write $t_j^\sharp=(\mathsf{w}_1,\dots,\mathsf{w}_k)$, where $k=k_j=|t_j|^2$, and we let $\ell=\ell_j=n_j(n_j+1)+1$.
The blocks $\Lambda_1(t_j),\dots,\Lambda_{n_j}(t_j)$ coincide with $\Lambda_1(\omega),\dots,\Lambda_{n_j}(\omega)$, and the next term of $t_j^\sharp$ is the term of index $n_j+1$ of~$t_j$.
Hence $\mathsf{w}_i=a_i$ for $i<\ell$, while $\mathsf{w}_\ell$ and $a_\ell$ are distinct elements of $\{1,2\}$.
Moreover, since $t_j^\sharp$ ends with the terms of~$t_j$ in reverse order, the terms $\mathsf{w}_k,\mathsf{w}_{k-1},\dots$ begin with $\omega_1,\dots,\omega_{n_j}$, and $k-\ell\ge(n_j+1)^2-\ell=n_j$.
Let $x\in\RP^1$ be distinct from~$\eta$.
The convergents of $[0;\mathsf{w}_1,\dots,\mathsf{w}_k]$ of index at least $\ell-1$ agree with~$\eta$ in their terms of index less than~$\ell$.
By the first fact, the convergents of index $k-1$ and~$k$ therefore tend to~$\eta$ as $j\to\infty$, while $q_{k-1}/q_k=[0;\mathsf{w}_k,\dots,\mathsf{w}_1]$ tends to~$[0;\omega]$.
Formula~\eqref{eq:Mobius-inverse} then shows that $\gamma_{t_j}(x)\to-[0;\omega]$ for every such~$x$.
At the point~$\eta$ we write $\gamma_{t_j}$ as the composite of the map $x\mapsto1/x$ followed by the maps $x\mapsto1/(x-\mathsf{w}_i)$, for $i=1,\dots,k$ in this order.
We then apply these $k+1$ factors to~$\eta$ one at a time.
The first $\ell$ factors send $\eta$ to its complete quotient $\eta_\ell=[a_\ell;a_{\ell+1},\dots]$, which lies in $(a_\ell,a_\ell+1)$.
The next factor sends $\eta_\ell$ to $z_j=1/(\eta_\ell-\mathsf{w}_\ell)$, which lies in $(\frac12,1)$ if $(a_\ell,\mathsf{w}_\ell)=(2,1)$ and in $(-\infty,-1)$ if $(a_\ell,\mathsf{w}_\ell)=(1,2)$;
in either case $z_j<1$.
The remaining factors form the inverse of the map $M'(x)=[\mathsf{w}_{\ell+1};\mathsf{w}_{\ell+2},\dots,\mathsf{w}_{k-1},\mathsf{w}_k+1/x]$, and therefore $\gamma_{t_j}(\eta)=M'^{-1}(z_j)$.
Let $p'_i/q'_i$ be the convergents of $[\mathsf{w}_{\ell+1};\mathsf{w}_{\ell+2},\dots,\mathsf{w}_k]$, the last of which has index $i=k-\ell-1\ge n_j-1$.
Every convergent of index at least~$1$ has the form $\mathsf{w}_{\ell+1}+1/y$ with $1\le y\le3$, since all its partial quotients lie in $\{1,2\}$.
Each of these convergents is therefore at least~$\frac43$, and its distance to~$z_j$ exceeds~$\frac13$.
The convergents of index $i-1$ and~$i$ differ by $1/(q'_{i-1}q'_i)\le1/(F_iF_{i+1})$, and $q'_{i-1}/q'_i=[0;\mathsf{w}_k,\dots,\mathsf{w}_{\ell+2}]$ tends to~$[0;\omega]$, since its first $n_j-1$ terms are $\omega_1,\dots,\omega_{n_j-1}$.
By~\eqref{eq:Mobius-inverse}, applied to~$M'$, we conclude that $\gamma_{t_j}(\eta)\to-[0;\omega]$, and hence that $\gamma_{t_j}\to c_\omega$.

\emph{The closure.}
By Corollary~\ref{cor:euclid} and Theorem~\ref{thm:nip-baire1def}, the space $E(\RP^1,\Gamma)$ is a Rosenthal compactum.
Consequently, every point of the closure of $\Delta_{\mathcal{C}}$ that does not belong to $\Delta_{\mathcal{C}}$ is the limit of a sequence $(\gamma_{t_j})$ with pairwise distinct words~$t_j$~\cite[Theorem~3F]{BFT_1978_PCompactBaire}.
A compactness argument yields a subsequence of $(t_j)$ that converges to some $\omega\in\{1,2\}^{\mathbb{N}}$.
A further subsequence consists either of initial segments of~$\omega$ or of words that are not initial segments of~$\omega$.
By the two preceding paragraphs, the corresponding maps converge to~$\tau_\omega$ in the first case and to~$c_\omega$ in the second.
Both kinds of limit occur for every~$\omega$, and the closure of $\Delta_{\mathcal{C}}$ is therefore the set described in~(i).

\emph{The homeomorphism.}
Identify $\{1,2\}^{\mathbb{N}}$ with $2^{\mathbb{N}}$ and $\{1,2\}^{<\mathbb{N}}$ with $2^{<\mathbb{N}}$ by replacing the letters $1$ and~$2$ with $0$ and~$1$.
Under this identification, define a map $\Theta$ on $\hat D(2^{\mathbb{N}})$ by $\Theta(v_t,t^\frown0^\infty)=\gamma_t$, $\Theta(\delta_\omega,\omega)=\tau_\omega$ and $\Theta(\mathbf{0},\omega)=c_\omega$.
The map $\Theta$ is surjective by the preceding paragraph.
To see that $\Theta$ is injective, note first that the constant maps $c_\omega$ are pairwise distinct, because distinct sequences~$\omega$ give distinct numbers $[0;\omega]$.
For the same reason, the maps $\tau_\omega$ are pairwise distinct, and none of them is constant.
Neither kind of map is injective, whereas the maps $\gamma_t$ are homeomorphisms.
Finally, the maps $\gamma_t$ are pairwise distinct, because $t^\sharp$ determines~$t$ and the Möbius map $M_{\mathsf{w}}$ determines the word~$\mathsf{w}$.
Indeed, a rational number in $[0,1]$ has at most two expansions $[0;\mathsf{w}_1,\dots,\mathsf{w}_k]$, whose lengths differ by one.
The number $M_{\mathsf{w}}(\infty)=[0;\mathsf{w}_1,\dots,\mathsf{w}_k]$ and the determinant $(-1)^{k+1}$ of a matrix that represents~$M_{\mathsf{w}}$ therefore determine~$\mathsf{w}$.
Since $\hat D(2^{\mathbb{N}})$ is first countable, the continuity of $\Theta$ may be tested on sequences.
A sequence that converges to a point $(\delta_\omega,\omega)$ consists, apart from that point, of points $(v_t,t^\frown0^\infty)$ with $t$ an initial segment of~$\omega$, whose lengths tend to infinity.
By the paragraph on limits along initial segments, the image of this sequence under~$\Theta$ converges to~$\tau_\omega$.
A sequence that converges to a point $(\mathbf{0},\omega)$ consists, apart from that point, of points of at most three kinds.
Points of the first kind have the form $(\mathbf{0},\omega')$, and points of the second kind the form $(\delta_{\omega'},\omega')$, with $\omega'\to\omega$ and $\omega'\neq\omega$.
Points of the third kind have the form $(v_t,t^\frown0^\infty)$, with $t$ converging to~$\omega$ and eventually not an initial segment of~$\omega$.
The images of the points of the first kind converge to~$c_\omega$ because $[0;\omega']$ depends continuously on~$\omega'$.
The same holds for the images of the points of the third kind, by the paragraph on limits along other words.
For the second kind, observe that the equality $\tau_{\omega'}(x)=-[0;\omega']$ fails at a given point $x\in\RP^1$ only if $\eta(\omega')=x$, which happens for at most one~$\omega'$.
The images of the points of the second kind therefore converge to~$c_\omega$ as well.
Hence $\Theta$ is a continuous bijection of a compact space onto a Hausdorff space, and therefore a homeomorphism;
it maps $D(2^{\mathbb{N}})$ onto $\{\tau_\omega:\omega\in\{1,2\}^{\mathbb{N}}\}\cup\{c_\omega:\omega\in\{1,2\}^{\mathbb{N}}\}$.
The set $\Delta_{\mathcal{C}}$ is countable, and it has the NIP because $\Gamma$ has it (see Corollary~\ref{cor:euclid}).
The closure of $\Delta_{\mathcal{C}}$, a copy of $\hat D(2^{\mathbb{N}})$, is first countable but not hereditarily separable (see Section~\ref{sec:trichotomy}), and $\Delta_{\mathcal{C}}$ is therefore $\NIP_1$ but not $\NIP_2$.

(ii) We first describe the deep computations that are not computations.
For $y,z,w\in\RP^1$, let $\zeta_{y,z,w}$ be the map that sends $z$ to~$w$ and every other point to~$y$.
Let $\xi\in E(\RP^1,\Gamma)\setminus\Gamma$ be given;
since $E(\RP^1,\Gamma)$ is a Rosenthal compactum, the map~$\xi$ is the limit of a sequence $(g_n)_{n\in\mathbb{N}}$ of pairwise distinct elements of~$\Gamma$.
Represent each $g_n$ by an integer matrix $\mathbf{G}_n$ with determinant $\pm1$.
Since only finitely many integer matrices have norm below a given bound, the norms $\|\mathbf{G}_n\|$ tend to infinity.
Passing to a subsequence, we may assume that the normalized matrices $\mathbf{G}_n/\|\mathbf{G}_n\|$ converge to a matrix $\mathbf{G}$ of norm~$1$.
Since $\det(\mathbf{G}_n/\|\mathbf{G}_n\|)=\pm\|\mathbf{G}_n\|^{-2}\to0$, the matrix $\mathbf{G}$ has rank one, and $\mathbf{G}\mathbf{e}=\langle\mathbf{e},\mathbf{e}'\rangle\,\mathbf{e}''$ for all $\mathbf{e}\in\mathbb{R}^2$, for some nonzero vectors $\mathbf{e}',\mathbf{e}''\in\mathbb{R}^2$.
If $x\in\RP^1$ is the line spanned by a vector~$\mathbf{e}$ with $\langle\mathbf{e},\mathbf{e}'\rangle\neq0$, then $\mathbf{G}_n\mathbf{e}/\|\mathbf{G}_n\|\to\langle\mathbf{e},\mathbf{e}'\rangle\,\mathbf{e}''\neq0$, and $g_n(x)$ converges to the line $y$ spanned by~$\mathbf{e}''$.
Thus, if $z$ is the line orthogonal to~$\mathbf{e}'$ and $w=\xi(z)$, then $\xi=\zeta_{y,z,w}$.

Now let $K\subseteq E(\RP^1,\Gamma)$ be compact and hereditarily separable, and let $Z_K$ be the set of points $z\in\RP^1$ for which some map $\zeta_{y,z,w}$ with $y\neq w$ belongs to~$K$;
such a map determines~$z$.
Choose, for each $z\in Z_K$, one such map $\zeta_z=\zeta_{y_z,z,w_z}\in K$.
We claim that the set $\{\zeta_z:z\in Z_K\}$ is discrete.
Indeed, fix a metric~$d$ on~$\RP^1$ and a point $z\in Z_K$, let $\varepsilon=d(y_z,w_z)/2$, and choose distinct points $x,x'\in\RP^1\setminus\{z\}$.
Consider the neighborhood of~$\zeta_z$ formed by the maps~$g$ with $d(g(z),w_z)<\varepsilon$, $d(g(x),y_z)<\varepsilon$ and $d(g(x'),y_z)<\varepsilon$.
Given $z'\in Z_K\setminus\{z\}$, one of the points $x$ and~$x'$ differs from~$z'$, and we may suppose that $x$ does.
Then $\zeta_{z'}(z)=\zeta_{z'}(x)=y_{z'}$, since neither $z$ nor $x$ equals~$z'$.
Were $\zeta_{z'}$ in the neighborhood of~$\zeta_z$ just described, we would have $d(y_z,w_z)\le d(y_z,y_{z'})+d(y_{z'},w_z)<2\varepsilon$, which is impossible.
Since $K$ is hereditarily separable, the discrete set $\{\zeta_z:z\in Z_K\}$ is countable, and so is~$Z_K$.
Evaluation at the points of the countable set $Z_K\cup\mathbb{Q}$ separates the members of~$K$.
To see this, note first that two members of~$\Gamma$ that agree on the dense set~$\mathbb{Q}$ agree everywhere.
A member of~$\Gamma$, moreover, is injective, whereas a map $\zeta_{y,z,w}$ takes the value~$y$ at all points of~$\mathbb{Q}$ but at most one.
Two maps $\zeta_{y,z,w}$ and $\zeta_{y',z',w'}$ with $y\neq y'$ differ at a point of $\mathbb{Q}\setminus\{z,z'\}$.
Lastly, let $\zeta_{y,z,w}$ and $\zeta_{y,z',w'}$ be two distinct maps, labelled so that the first satisfies $w\neq y$.
These two maps differ at the point $z\in Z_K$, at which the second takes the value~$y$ if $z'\neq z$ and the value $w'\neq w$ if $z'=z$.
Hence restriction to $Z_K\cup\mathbb{Q}$ is a continuous injection of the compact space~$K$ into the metrizable space $(\RP^1)^{Z_K\cup\mathbb{Q}}$, and $K$ is metrizable.

Finally, the split Cantor set $S(2^{\mathbb{N}})$ is hereditarily separable and not metrizable, and so is the closed subspace of $\hat D(S(2^{\mathbb{N}}))$ formed by its points with first coordinate~$\mathbf{0}$.
In view of the preceding paragraph, neither space embeds into $E(\RP^1,\Gamma)$.
The subspace just described is hereditarily separable, since the projection onto the second coordinate maps this subspace homeomorphically onto a subspace of $S(2^{\mathbb{N}})$.
To see that this subspace is not metrizable, let $\varpi\in2^{\mathbb{N}}$ have infinitely many zeros, and let $a$ be the branch of the subtree $\{s_t:t\in2^{<\mathbb{N}}\}$ of Section~\ref{sec:five-compacta} that extends the words $s_{\varpi\restr n}$.
Let $u_1,u_2,\dots$ be the words $s_{(\varpi\restr n)^\frown1}$ with $\varpi(n)=0$, listed by increasing~$n$.
These words are pairwise incomparable and converge to~$a$, and the sequences $u_k^\frown0^\infty$ decrease, since the canonical isomorphism preserves the lexicographic order.
By the description of $\hat D(S(2^{\mathbb{N}}))$ in Section~\ref{sec:five-compacta}, the point $(\mathbf{0},\chi_{(a})$ belongs to the subspace just described, and distinct sequences~$\varpi$ give distinct branches.
In $S(2^{\mathbb{N}})$, each point $\chi_{(a}$ has the neighborhoods $\{\chi_{(a}\}\cup\{\chi_{[x},\chi_{(x}:a<x<b\}$, for $b>a$, and the argument of Lemma~\ref{lem:double-arrow-metrizable} shows that a metrizable subspace of $S(2^{\mathbb{N}})$ contains only countably many of the points~$\chi_{(a}$.
\end{proof}

\begin{figure}[htbp]
  \centering
  \includegraphics[width=\textwidth]{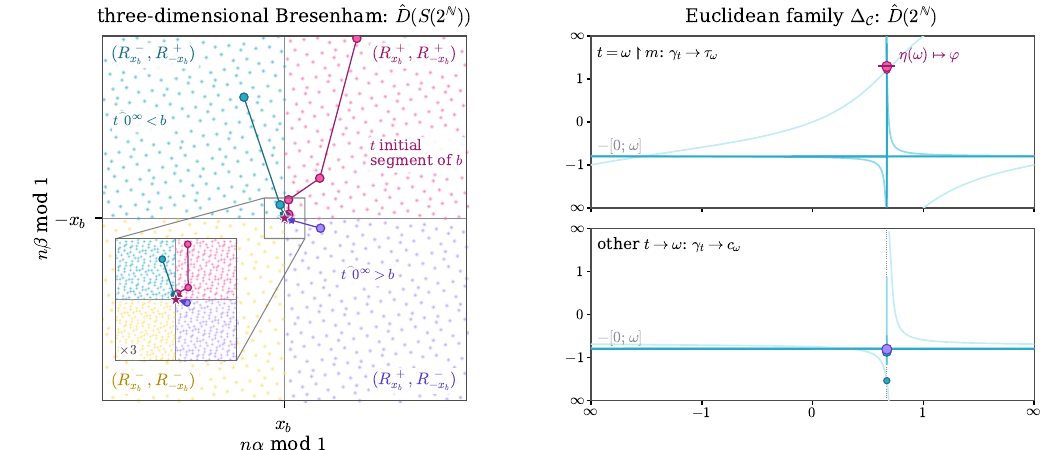}
  \caption{Case~(ii) of Todorčević’s trichotomy, realized in two ways (see Theorem~\ref{thm:Todorcevic-trichotomy}).
    Left: the three-dimensional Bresenham algorithm, with $\alpha$ the reciprocal of the plastic number and $\beta=\alpha^2$.
    The small dots are the points $(n\alpha,n\beta)$ modulo~$1$ near $(x_b,-x_b)$, for $b=(1,0,1,0,\dots)$ and $x_b=\frac34$;
    the inset enlarges the center three times.
    Iterates $(\sigma\times\sigma)^n$ whose points approach $(x_b,-x_b)$ within one quadrant converge to the pair written there (see Theorem~\ref{thm:Bresenham3}).
    Each entry $R^\pm$ of a pair is a deep computation of the planar algorithm;
    it rotates by its subscript and breaks ties upward~($+$) or downward~($-$) (see Theorem~\ref{thm:Bresenham}).
    The large dots mark the iterates of Remark~\ref{rem:Bresenham3-explicit}, in the three patterns of the family~$D_7$.
    Right: the Möbius maps $\gamma_t$ of Proposition~\ref{prop:euclid-duplicate} on the projective line, for $\omega=(1,2,1,2,\dots)$.
    Above, the word~$t$ runs through initial segments of~$\omega$, and the maps converge to~$\tau_\omega$.
    This limit sends a point $\eta(\omega)$ defined from~$\omega$ to the golden ratio~$\varphi$, and every other point to $-[0;\omega]=1-\sqrt3$.
    Below, the word~$t$ branches off from~$\omega$ below or above it (colors as on the left), and the maps converge to the constant map~$c_\omega$ with value $-[0;\omega]$.
    The large dots mark the values at~$\eta(\omega)$.
    Thus the two sides of~$b$ have distinct limits, while the two sides of~$\omega$ share one.}
  \label{fig:trichotomy-ii}
\end{figure}

Together with Corollary~\ref{cor:euclid}, Proposition~\ref{prop:euclid-duplicate} shows that the space of deep computations of the Euclidean CCS contains copies of $A(2^{\mathbb{N}})$ and of $D(2^{\mathbb{N}})$, but no copy of $S(2^{\mathbb{N}})$ or of $\hat D(S(2^{\mathbb{N}}))$.
At the level $\NIP_1$ but not $\NIP_2$, the set $\Delta_{\mathcal{C}}$ thus realizes the alternative of Todorčević's trichotomy opposite to the one realized by the three-dimensional Bresenham algorithm (see Theorem~\ref{thm:Bresenham3}).
The two alternatives are drawn side by side in Figure~\ref{fig:trichotomy-ii}.

\emph{Free reduction.}
A second example at the level $\NIP$ but not $\NIP_1$ is the free group $F_2$ on two generators $a,b$,
acting on its boundary~$\partial F_2$.
This boundary, a Cantor set, is the space of infinite reduced words in $a,a^{-1},b,b^{-1}$.
A generator acts on an infinite word by prepending a letter and cancelling it
against the first letter of the word if the two letters are inverse, as a stack does.
With the predicates $w\mapsto\Brack{w_n=x}$, for $n\in\mathbb{N}$ and each letter $x$, the boundary $\partial F_2$ and the action of~$F_2$ form a dynamical CCS.
Glasner and Megrelishvili proved that this system is tame but not
tame$_1$~\cite[Corollary~6.6 and Example~6.7]{glasner2022tame}.
By Proposition~\ref{prop:cascade-tame}, this CCS is therefore $\NIP$ but not $\NIP_1$.

\emph{An open question.}
Both examples at this level are actions of groups with two generators.
Whether a single map can reach this level, that is, whether some cascade CCS
is $\NIP$ but not $\NIP_1$, is Question~\ref{q:cascade-NIP1}.
For actions of $\mathbb{Z}$, and more generally of abelian groups, the
existence of a tame metrizable system whose enveloping semigroup is not first
countable is Question~6.8 of Glasner and Megrelishvili~\cite{glasner2022tame}.
To our knowledge, this question remains open.

\subsection{Beyond NIP}
\label{sec:beyond-NIP}
Chaotic iteration leads beyond the hierarchy, to sets of computations with the independence property.
In a set of computations with the independence property, the computations realize, for some feature and on suitable inputs, every finite pattern of high and
low values.
We close the section with examples of chaotic iteration, among them Newton’s method on the Riemann sphere, the doubling map of the circle and the logistic map beyond the Feigenbaum parameter.

\emph{Newton's method on the Riemann sphere.}
On its Julia set the Newton map is chaotic, and the chaos is visible to the
predicates of the CCS.
We exhibit the independence in the simplest case, the square roots of~$1$.

\begin{prop}\label{prop:Newton-wild}
Let $p(z)=z^2-1$, let $N(z)=\frac12\bigl(z+\frac1z\bigr)$ be its Newton map on~$\hatC$, and let $P_2$ be the second coordinate of the inverse stereographic projection $\hatC\to S^2$ of Section~\ref{sec:examples}.
The functions $P_2\circ N^{3k}$, $k\in\mathbb{N}$, form an independent sequence:
for all finite disjoint sets $I,J\subseteq\mathbb{N}$ there is $z\in\hatC$ with
$P_2(N^{3k}z)\ge\frac12$ for $k\in I$ and $P_2(N^{3k}z)\le-\frac12$ for $k\in J$.
In particular, the Newton CCS of $p$ is not $\NIP$.
\end{prop}

\begin{proof}
The Möbius transformation $\psi(z)=(z-1)/(z+1)$ conjugates $N$ to the squaring
map: $\psi(N(z))=\psi(z)^2$.
The map $\psi$ sends the imaginary axis, together with $\infty$, onto the unit circle, with
$\psi^{-1}(e^{i\theta})=i\cot(\theta/2)$.
Hence, for $z=i\cot(\theta/2)$, we have $N^n(z)=i\cot(2^{n}\theta/2)$ and
\[
P_2\bigl(N^n(z)\bigr)=\frac{2\cot(2^{n-1}\theta)}{\cot^2(2^{n-1}\theta)+1}=\sin(2^n\theta).
\]
Writing $\theta=2\pi\sum_{j\ge1}x_j2^{-j}$ with bits $x_j$, we find that the fractional part of $2^{3k}\theta/2\pi$ has binary expansion
$0.x_{3k+1}x_{3k+2}x_{3k+3}\cdots$.
If $(x_{3k+1},x_{3k+2},x_{3k+3})=(0,1,0)$, this fractional part lies in
$[\frac14,\frac38]$, and hence $\sin(2^{3k}\theta)\ge\sin\frac{3\pi}{4}>\frac12$.
If the block is $(1,1,0)$, the fractional part lies in $[\frac34,\frac78]$, whence $\sin(2^{3k}\theta)<-\frac12$.
Choosing the block $(0,1,0)$ for $k\in I$ and $(1,1,0)$ for $k\in J$ gives the
required~$z$.
\end{proof}

The proof also shows that the doubling map $x\mapsto2x$ of the circle
$\mathbb{T}=\mathbb{R}/\mathbb{Z}$, with the predicates $\cos2\pi x$ and
$\sin2\pi x$, is not $\NIP$: with $x=\theta/2\pi$, the functions
$\sin(2\pi\,2^{3k}x)$ form an independent sequence.
The doubling map shifts the binary digits of~$x$, and the independence comes
from the freedom to choose those digits.

On the real projective line, by contrast, the same algorithm is tame.
The Julia set of $N$, the imaginary axis together with $\infty$, meets
$\RP^1$ only at~$\infty$, and on $\RP^1$ the iterates converge to the deep
computation of Section~\ref{exa:square-roots-compute} for $a=1$.
By Proposition~\ref{prop:convergent}, the corresponding CCS is therefore $\NIP_3$ (see Section~\ref{sec:ex-NIP3}).
Whether Newton's method is tame is thus decided by the space of states on which it acts, and not by the
formula that defines its step.

\emph{Rational maps of the sphere.}
The chaos exhibited in Proposition~\ref{prop:Newton-wild} is not peculiar to the square roots of~$1$.
Chaos of this kind, and with it the independence property, is present on the Julia set of every rational map of degree at least~$2$.
The proof of Proposition~\ref{prop:rational-IP} below rests on two classical facts~\cite[Sections~4 and~14]{Milnor2006} about the Julia set $\mathcal{J}_g$ of such a map~$g$ (see Section~\ref{sec:Newton-map}).
The first is that $\mathcal{J}_g$ is a nonempty compact set without isolated points that satisfies $g^{-1}(\mathcal{J}_g)=\mathcal{J}_g$.
The second is that, for every neighborhood $U$ of a point of $\mathcal{J}_g$, the image $g^m(U)$ contains $\mathcal{J}_g$ for all sufficiently large~$m$.
This second fact is a consequence of Montel's theorem and of the density of repelling cycles in~$\mathcal{J}_g$.

\begin{prop}\label{prop:rational-IP}
Let $g\colon\hatC\to\hatC$ be a rational map of degree at least~$2$, and consider the CCS with states~$\hatC$, the predicates $P_1,P_2,P_3$ of Section~\ref{sec:examples}, and the semigroup of iterates of~$g$.
There are $j\in\{1,2,3\}$, a positive integer~$n$ and reals $a<b$ such that the functions $P_j\circ g^{nk}$, $k\in\mathbb{N}$, form an independent sequence:
for all finite disjoint sets $I,J\subseteq\mathbb{N}$ there is $z\in\hatC$ with $P_j(g^{nk}(z))\le a$ for $k\in I$ and $P_j(g^{nk}(z))\ge b$ for $k\in J$.
In particular, the CCS is not~$\NIP$, and neither are the Newton CCSs for $q(z)=z^3-2z+2$ and for $u(z)=z^3-1$.
\end{prop}

\begin{proof}
Choose distinct points $z_0,z_1\in\mathcal{J}_g$ and an index $j$ with $P_j(z_0)\neq P_j(z_1)$, say $P_j(z_0)<P_j(z_1)$, and reals $a<b$ strictly between these two values.
By continuity, the points $z_0$ and $z_1$ have closed neighborhoods $U_0$ and~$U_1$ such that $P_j\le a$ on~$U_0$ and $P_j\ge b$ on~$U_1$, and we may choose $n$ with $g^n(U_0)\supseteq\mathcal{J}_g$ and $g^n(U_1)\supseteq\mathcal{J}_g$.
Let $I,J\subseteq\mathbb{N}$ be finite and disjoint, and let $m$ exceed all their elements;
put $U^{(k)}=U_0$ for $k\in I$ and $U^{(k)}=U_1$ for $k\notin I$.
Let $C_m=U^{(m)}\cap\mathcal{J}_g$ and, for $k=m-1,\dots,0$, let $C_k=U^{(k)}\cap\mathcal{J}_g\cap g^{-n}(C_{k+1})$.
Each $C_k$ is nonempty: since $g^n(U^{(k)})\supseteq\mathcal{J}_g\supseteq C_{k+1}$, some point of $U^{(k)}$ is mapped by $g^n$ into~$C_{k+1}$, and that point lies in $\mathcal{J}_g$ because $g^{-1}(\mathcal{J}_g)=\mathcal{J}_g$.
Every $z\in C_0$ satisfies $g^{nk}(z)\in U^{(k)}$ for $k\le m$, and hence $P_j(g^{nk}(z))\le a$ for $k\in I$ and $P_j(g^{nk}(z))\ge b$ for $k\in J$.
For the last assertion, note that the Newton maps $N_q(z)=(2z^3-2)/(3z^2-2)$ and $N_u(z)=(2z^3+1)/(3z^2)$ of Section~\ref{sec:examples} have degree~$3$, since in each quotient the numerator and the denominator have no common zero.
\end{proof}

\begin{figure}[htbp]
  \centering
  \includegraphics[width=\textwidth]{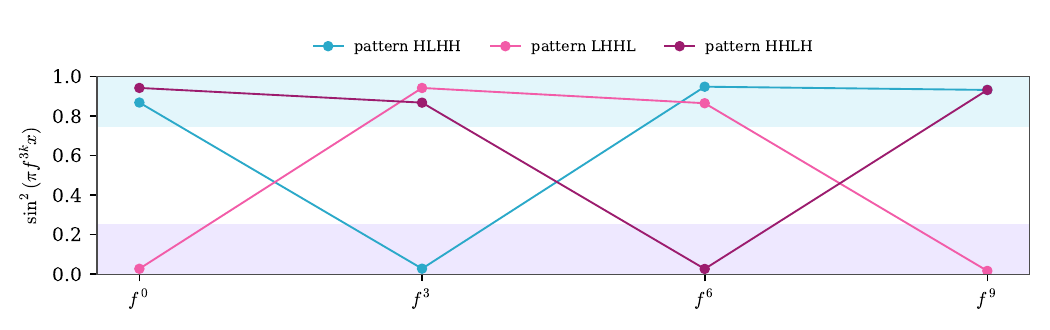}
  \caption{Independence for the doubling map $x\mapsto2x$ of the circle and for the logistic map $f_4$. Along an orbit of the doubling map, the values of $\sin^2(\pi x)=(1-\cos2\pi x)/2$ form an orbit of~$f_4$. At the times $0,3,6,9$, these values follow any prescribed pattern of high (H, at least $3/4$) and low (L, at most $1/4$) values. For this purpose, the binary digits of the initial state are chosen in blocks $011$ for H and $000$ for L.}
  \label{fig:doubling}
\end{figure}

\emph{The logistic map at $r=4$.}
The doubling map also governs the logistic map $f_4(x)=4x(1-x)$ of
Section~\ref{sec:ex-NIP3}, as Figure~\ref{fig:doubling} illustrates.

\begin{prop}\label{prop:logistic-four}
The functions $f_4^{3k}$, $k\in\mathbb{N}$, form an independent sequence on
$[0,1]$: for all finite disjoint sets $I,J\subseteq\mathbb{N}$ there is
$x\in[0,1]$ with $f_4^{3k}(x)\ge\frac34$ for $k\in I$ and
$f_4^{3k}(x)\le\frac14$ for $k\in J$.
In particular, the logistic CCS for $r=4$ is not $\NIP$.
\end{prop}

\begin{proof}
Let $h(\theta)=\sin^2(\theta/2)$ for $\theta\in\mathbb{R}$, and note the identity $f_4(h(\theta))=4\sin^2(\theta/2)\cos^2(\theta/2)=\sin^2\theta=h(2\theta)$.
Iterated, this identity yields $f_4^{n}(h(\theta))=h(2^n\theta)=\sin^2(2^{n-1}\theta)$.
Writing $\theta=2\pi\sum_{j\ge1}x_j2^{-j}$ with bits $x_j$, we find that the fractional part of $2^{3k}\theta/2\pi$ has binary expansion $0.x_{3k+1}x_{3k+2}x_{3k+3}\cdots$.
Moreover, the angle $2^{3k-1}\theta$ differs from $\pi$ times this fractional part by an integer multiple of~$\pi$.
If $(x_{3k+1},x_{3k+2},x_{3k+3})=(0,1,1)$, this fractional part lies in $[\frac38,\frac12]$, and hence $f_4^{3k}(h(\theta))=\sin^2(2^{3k-1}\theta)\ge\sin^2\frac{3\pi}8>\frac34$.
If the block is $(0,0,0)$, the fractional part lies in $[0,\frac18]$, whence $f_4^{3k}(h(\theta))\le\sin^2\frac{\pi}8<\frac14$.
Choosing the block $(0,1,1)$ for $k\in I$ and $(0,0,0)$ for $k\in J$, and putting $x=h(\theta)$, gives the required point.
\end{proof}

Thus the logistic family, which is stable for $0<r\le1$ and $\NIP_3$ but not stable for $1<r\le3$ (see Proposition~\ref{prop:logistic}), leaves the hierarchy at $r=4$.
The passage takes place at the Feigenbaum parameter $r_\infty\approx3.5699$, the limit of the cascade of period doublings, and the family skips the intermediate levels.
Indeed, in Proposition~\ref{prop:logistic-all} below we show that the family remains $\NIP_3$ below~$r_\infty$ and has the independence property beyond~$r_\infty$, periodic windows included.
Moreover, in Remark~\ref{rem:logistic-feigenbaum} we sketch the argument that the family is still $\NIP_3$ at~$r_\infty$.
The arguments rest on the kneading theory of Milnor and Thurston~\cite{MilnorThurston1988} and on the theory of interval maps all of whose periodic points have periods that are powers of~$2$~\cite{BlockCoppel1992,deMeloVanStrien1993}.

\begin{prop}\label{prop:logistic-all}
Let $r_\infty$ be the supremum of the parameters $r\in(0,4]$ for which every periodic point of~$f_r$ has a power of~$2$ as its period.
\begin{enumerate}[(i)]
\item
  For $3<r<r_\infty$, the logistic CCS is $\NIP_3$ but not stable.
\item
  For $r_\infty<r\le4$, there are a positive integer~$n$ and reals $a<b$ such that the functions $f_r^{nk}$, $k\in\mathbb{N}$, form an independent sequence:
  for all finite disjoint sets $I,J\subseteq\mathbb{N}$ there is $x\in[0,1]$ with $f_r^{nk}(x)\le a$ for $k\in I$ and $f_r^{nk}(x)\ge b$ for $k\in J$.
  In particular, the logistic CCS is not~$\NIP$.
\end{enumerate}
\end{prop}

\begin{proof}
Writing $f$ for~$f_r$, we use three classical facts~\cite{MilnorThurston1988,BlockCoppel1992,deMeloVanStrien1993}.
First, by the monotonicity of the kneading invariant in the logistic family, every periodic point of~$f$ has a power of~$2$ as its period if and only if $r\le r_\infty$.
Moreover, for $r<r_\infty$ only finitely many periods occur, whence $f$ has only finitely many periodic points.
Second, for $2<r<4$, every limit point of an orbit that starts in $(0,1)$ lies in the interval $[f^2(\frac12),f(\frac12)]$, which does not contain~$0$.
Third, a continuous map of a compact interval with only finitely many periodic points has every orbit asymptotic to a periodic orbit~\cite{BlockCoppel1992}.

(i) Suppose that $3<r<r_\infty$, and let $N$ be a common multiple of the periods of~$f$, which are powers of~$2$.
If the orbit of~$x$ is asymptotic to the periodic orbit~$\mathcal{Q}$, then $f^{Nk}(x)$ converges to a point of~$\mathcal{Q}$.
Indeed, the map $f^N$ fixes~$\mathcal{Q}$ pointwise, and by uniform continuity the point of~$\mathcal{Q}$ nearest to $f^i(x)$ is eventually carried by~$f$ to the point of~$\mathcal{Q}$ nearest to $f^{i+1}(x)$.
Let $\ell_r(x)$ be the limit of the sequence $(f^{Nk}(x))_{k\in\mathbb{N}}$.
Then $f^{Nk+i}\to f^i\circ\ell_r$ pointwise for $0\le i<N$, and every limit of the iterates along a nonprincipal ultrafilter on~$\mathbb{N}$ is one of these $N$ maps.
Each of the maps $f^i\circ\ell_r$ is of Baire class~$1$, being the pointwise limit of a sequence of iterates.
Hence $E([0,1],f)$ is countable and therefore metrizable, and by Proposition~\ref{prop:cascade-tame} the CCS is $\NIP_3$.
The map $\ell_r$, which fixes~$0$ and maps $(0,1)$ into the interval $[f^2(\frac12),f(\frac12)]$, is thus discontinuous at~$0$, and the CCS is not stable.

(ii) Suppose that $r_\infty<r\le4$; by the first fact, the map $f$ then has a periodic point whose period is not a power of~$2$.
Hence, by a theorem of Block and Coppel~\cite[Chapter~II]{BlockCoppel1992}, some iterate $f^n$ is strictly turbulent.
Strict turbulence of~$f^n$ means that there are disjoint closed intervals $H_0$ and~$H_1$, with $H_0$ to the left of~$H_1$, such that $f^n(H_0)$ and $f^n(H_1)$ both contain $H_0\cup H_1$.
Let $a=\max H_0$ and $b=\min H_1$, let $I,J\subseteq\mathbb{N}$ be finite and disjoint, and let $m$ exceed all their elements;
put $H^{(k)}=H_0$ for $k\in I$ and $H^{(k)}=H_1$ for $k\notin I$.
Let $C_m=H^{(m)}$ and, for $k=m-1,\dots,0$, let $C_k=H^{(k)}\cap f^{-n}(C_{k+1})$.
Each set $C_k$ is nonempty, since $f^n(H^{(k)})\supseteq C_{k+1}$.
Every $x\in C_0$ satisfies $f^{nk}(x)\in H^{(k)}$ for $k\le m$, and hence $f^{nk}(x)\le a$ for $k\in I$ and $f^{nk}(x)\ge b$ for $k\in J$.
\end{proof}

\begin{rem}\label{rem:logistic-feigenbaum}
At the Feigenbaum parameter itself the logistic CCS is again $\NIP_3$ but not stable; we sketch the argument, which rests on the structure of infinitely renormalizable maps~\cite{BlockCoppel1992,deMeloVanStrien1993}.
For $r=r_\infty$ and every $m\in\mathbb{N}$, the map $f=f_r$ has a cycle of $2^m$ pairwise disjoint closed intervals, permuted cyclically by~$f$ and contained in the intervals of the cycle of level $m-1$.
Moreover, the lengths of the intervals of level~$m$ tend to~$0$ as $m\to\infty$.
Finally, every orbit that starts in $(0,1)$ either lands on a periodic orbit or enters the cycle of every level.
Label the intervals of the cycle of level~$m$ by the residues modulo~$2^m$, in such a way that $f$ adds~$1$ to the label.
Let $\mathcal{U}$ be a nonprincipal ultrafilter on~$\mathbb{N}$, and let $\varsigma$ be the limit of $n$ along~$\mathcal{U}$ in the ring $\mathbb{Z}_2$ of $2$-adic integers.
If the orbit of~$x$ lands on a periodic orbit, of period~$2^m$, then $f^n(x)$ is eventually determined by the residue of $n$ modulo~$2^m$.
Otherwise, for each~$m$, the point $f^n(x)$ eventually lies in the interval of level~$m$ whose label differs from~$n$ by a certain residue modulo~$2^m$.
This residue depends on $x$ and~$m$ but not on~$n$, and the intervals obtained for successive~$m$ shrink to a point.
In either case, the points $f^n(x)$ converge along~$\mathcal{U}$ to a limit $\ell^{\varsigma}(x)$ that depends only on~$\varsigma$.
The map $\ell^{\varsigma}$ is the pointwise limit of the iterates $f^{n_j}$ along any sequence of natural numbers $n_j$ that tends to infinity and converges to~$\varsigma$ in~$\mathbb{Z}_2$.
Each map $\ell^{\varsigma}$ is therefore of Baire class~$1$.
Hence $E([0,1],f)$ consists of the iterates and the maps $\ell^{\varsigma}$, $\varsigma\in\mathbb{Z}_2$.
The rational points of $[0,1]$, together with the point~$\frac12$, separate the members of $E([0,1],f)$:
distinct iterates are continuous and differ at a rational point.
Moreover, an iterate $f^n$ and a map $\ell^{\varsigma}$ differ at every rational point $x>0$ so small that $f^n(x)<f^2(\frac12)$, since $\ell^{\varsigma}$ maps $(0,1)$ into $[f^2(\frac12),f(\frac12)]$.
Finally, two maps $\ell^{\varsigma}$ and $\ell^{\varsigma'}$ with $\varsigma\neq\varsigma'$ differ at~$\frac12$, whose orbit lands on no periodic orbit.
Indeed, whenever $\varsigma\not\equiv\varsigma'$ modulo~$2^m$, the points $\ell^{\varsigma}(\frac12)$ and $\ell^{\varsigma'}(\frac12)$ lie in distinct intervals of level~$m$.
As in the proof of Proposition~\ref{prop:levels-tests}(i), the space $E([0,1],f)$ is therefore metrizable, and by Proposition~\ref{prop:cascade-tame} the CCS is $\NIP_3$.
The CCS is not stable, since each map $\ell^{\varsigma}$ fixes~$0$ and is discontinuous there.
\end{rem}

\emph{The Thue--Morse sequence.}
Even a sequence produced by a simple substitution rule can be wild.
The Thue--Morse sequence $t_n$ is the parity of the number of ones in the binary
expansion of~$n$.
The same sequence is the limit of the words obtained from $0$ by the
substitution $0\mapsto01$, $1\mapsto10$.
The two-sided sequences all of whose finite blocks occur in $(t_n)$ form the
Thue--Morse subshift $X_{TM}\subseteq2^{\mathbb{Z}}$, on which the shift $\sigma$ is a
homeomorphism.
Glasner and Megrelishvili observed that the Morse minimal system, which consists of the subshift $X_{TM}$ and the action of~$\mathbb{Z}$ generated by~$\sigma$, is not tame~\cite{GlasnerMegrelishvili2018}.
In Appendix~\ref{app:thue-morse} we give a self-contained proof, which relies only on the structure theorem for tame minimal systems.
By Proposition~\ref{prop:cascade-tame}, the CCS with states $X_{TM}$, the coordinate predicates and the group
$\{\sigma^n:n\in\mathbb{Z}\}$ is therefore not $\NIP$.
The cascade CCS of~$\sigma$, whose computations are only the iterates $\sigma^n$ with $n\in\mathbb{N}$, is not $\NIP$ either.
Indeed, the reversal of two-sided sequences maps $X_{TM}$ onto itself and conjugates~$\sigma$ to its inverse.
The reversal therefore turns an independent sequence of functions $f\circ\sigma^{n_k}$ with negative exponents~$n_k$ into one with positive exponents.
Both facts are stated, with complete proofs, in Proposition~\ref{prop:thue-morse}.

\subsection{Summary}
Table~\ref{tab:hierarchy} collects the examples of this section.
The Euclidean algorithm and free reduction are actions of groups with two
generators, and the family $\Delta_{\mathcal{C}}$ consists of elements of the first of these groups.
Each of the remaining entries of the table is a single map iterated.

\begin{table}[ht]
\centering
\small
\renewcommand{\arraystretch}{1.25}
\begin{tabular}{>{\raggedright\arraybackslash}p{2.1cm}>{\raggedright\arraybackslash}p{4.4cm}>{\raggedright\arraybackslash}p{4.6cm}}
\hline
Level & Computation & Deep computations; source\\
\hline
stable & rotation of the circle & the circle group\\
& logistic map, $0<r\le1$ & iterates and one continuous limit; Prop.~\ref{prop:logistic}\\
$\NIP_3$, not stable & power method & countable; Prop.~\ref{prop:power}\\
& Newton's method for $x^2-a$ on $\RP^1$ & countable; Prop.~\ref{prop:convergent}\\
& logistic map, $1<r\le r_\infty$ & countable for $r<r_\infty$, metrizable for $r=r_\infty$; Props.~\ref{prop:logistic}, \ref{prop:logistic-all}, Rem.~\ref{rem:logistic-feigenbaum}\\
$\NIP_2$, not $\NIP_3$ & Bresenham, plane (Sturmian shift) & isolated iterates and a double arrow circle; Thm.~\ref{thm:Bresenham}\\
$\NIP_1$, not $\NIP_2$ & Bresenham, space (two Sturmian shifts) & isolated iterates and a copy of $\mathbb{D}\times\mathbb{D}$; contains $\hat D(S(2^{\mathbb{N}}))$, not $D(2^{\mathbb{N}})$; Thm.~\ref{thm:Bresenham3}\\
& the family $\Delta_{\mathcal{C}}$ in $\mathrm{PGL}(2,\mathbb{Z})$ & a copy of $\hat D(2^{\mathbb{N}})$; Prop.~\ref{prop:euclid-duplicate}\\
& a single map, with $D(2^{\mathbb{N}})$ & open; Question~\ref{q:cascade-duplicate}\\
$\NIP$, not $\NIP_1$ & Euclidean algorithm: $\mathrm{PGL}(2,\mathbb{Z})$ on $\RP^1$ (group) & contains $A(2^{\mathbb{N}})$, of exact equality tests, and $D(2^{\mathbb{N}})$; Cor.~\ref{cor:euclid}, Prop.~\ref{prop:euclid-duplicate}\\
& free reduction: $F_2$ on $\partial F_2$ (group) & not first countable; \cite{glasner2022tame}\\
& a single map & open; Question~\ref{q:cascade-NIP1}\\
not $\NIP$ & doubling map of the circle & Prop.~\ref{prop:Newton-wild}\\
& Newton's method for $z^2-1$ on $\hatC$ & Prop.~\ref{prop:Newton-wild}\\
& Newton's method for $q$ and $u$ on $\hatC$ & Prop.~\ref{prop:rational-IP}\\
& logistic map, $r_\infty<r\le4$ & Props.~\ref{prop:logistic-four}, \ref{prop:logistic-all}\\
& Thue--Morse shift & Prop.~\ref{prop:thue-morse}\\
\hline
\end{tabular}
\vspace{6pt}
\caption{Computations at each level of the hierarchy and beyond it, with the cases left open.}
\label{tab:hierarchy}
\end{table}

\section{Randomized versions of NIP and Monte Carlo computability of deep computations}
\label{sec:measure-theoretic-NIP}

Replacing topology by measure, we dispense in this section with the countability of the set $\mathcal{P}$ of predicates, on which the theory of Sections~\ref{sec:Classification}--\ref{sec:iterative} depends.
In those sections we studied deep computations as \emph{topological} limits of computations.
There we showed that, for countable $\mathcal{P}$ and under the NIP, each deep computation is of the first Baire class on the type shards.
We also showed that, on each shard, each such deep computation is the limit of a sequence of computations.
That theory rests on the type shards being Polish, and hence on the countability of~$\mathcal{P}$, which is used in the proof of Theorem~\ref{thm:nip-baire1def}.
In place of membership in the first Baire class, we ask that a deep computation be \emph{universally measurable}, that is, measurable with respect to every Radon probability measure on the type shards.
For every such measure, each output feature of a universally measurable deep computation is then, by Luzin's theorem, continuous on a closed subset of the shard whose complement has arbitrarily small measure.
Hence, by the Stone--Weierstrass theorem, each such output feature is a uniform limit of polynomials in the features on that subset.
We call a deep computation that is universally measurable in this sense \emph{universally Monte Carlo computable} (see Definition~\ref{def:universally-essentially-computable}).
The role of the Bourgain--Fremlin--Talagrand dichotomy is then taken by its measure-theoretic form (see Theorem~\ref{thm:BFT-2F}).
The role of $\NIP$, for its part, is shared with \emph{Talagrand stability}, a measure-theoretic condition introduced by Talagrand in his study of Pettis integration (see Definition~\ref{def:Talagrand-stable}).

The first main result of the section, Theorem~\ref{thm:NIP-Monte-Carlo}, characterizes the NIP by the universal measurability of the restrictions of ultracomputations to the type shards.
The theorem derives Monte Carlo computability from the NIP thus characterized; the second main result, Corollary~\ref{cor:Talagrand-stab-Monte-Carlo}, derives the same conclusion from universal Talagrand stability.
By Theorem~\ref{thm:NIP-Talagrand-countable}, these two hypotheses, the NIP and universal Talagrand stability, coincide for countable sets of computations (see Corollary~\ref{cor:NIP-Talagrand-computations}).

\subsection{NIP and Monte Carlo computability of deep computations}

Of the two routes to universal Monte Carlo computability followed in this section, the first starts from the NIP and passes through universal measurability, which we introduce here.
Universal measurability is the measure-theoretic counterpart of membership in the first Baire class.
This first route leads to Theorem~\ref{thm:NIP-Monte-Carlo}, by which the deep computations of a confined set of computations with the NIP on shards and features are universally Monte Carlo computable.
By Fact~\ref{baire}, a function $f\colon X\to\mathbb{R}$ on a perfectly normal space $X$ is of Baire class~1 if and only if $f^{-1}[U]$ is an $F_\sigma$ subset of $X$ for every open $U\subseteq \mathbb{R}$.
The following definition replaces the $F_\sigma$ condition by measurability with respect to every Radon measure.

\begin{defn}
  \label{def:universal-measurability}
  Let $X$ be a Hausdorff space, and let $(Y,\Sigma)$ be a measurable space.
  A map $f:X\rightarrow Y$ is \emph{universally measurable} (with respect to $\Sigma$) if $f^{-1}(E)$ is $\mu$-measurable for every Radon measure $\mu$ on $X$ and every $E\in\Sigma$.
  When $Y=\mathbb{R}$, we will always take $\Sigma=\mathcal{B}(\mathbb{R})$, the Borel $\sigma$-algebra of $\mathbb{R}$.
\end{defn}

If $X$ is a compact Hausdorff space, then Radon measures on $X$ are finite, and every nonzero Radon measure on $X$ becomes a probability measure after normalization.
Moreover, the sets $E\subseteq\mathbb{R}$ whose preimages under a given function are $\mu$-measurable form a $\sigma$-algebra, which contains the Borel sets as soon as it contains the open sets.
Together, these two observations yield the following simplification of Definition~\ref{def:universal-measurability} for compact spaces.

\begin{rem}
  \label{rem:probability}
  If $X$ is compact, then a function $f:X\rightarrow\mathbb{R}$ is universally measurable if and only if $f^{-1}(U)$ is $\mu$-measurable for every Radon probability measure $\mu$ on $X$ and every open set $U\subseteq\mathbb{R}$.
\end{rem}

Universal measurability thus asks for measurability with respect to every Radon measure on the domain, that is, with respect to every admissible way of measuring its subsets.
An account of universally measurable sets and functions may be found in the \emph{Handbook of measure theory}~\cite[Chapters 1 and 2]{Pap:2002}.

\begin{notn}
  Following Bourgain, Fremlin and Talagrand~\cite{BFT_1978_PCompactBaire}, we denote by $M_r(X)$ the collection of all universally measurable real-valued functions on~$X$.
  For a Radon measure $\mu$ on $X$, we further denote by $\mathscr{M}^0(X,\mu)$ the collection of all $\mu$-measurable real-valued functions on~$X$.
\end{notn}

The maps that arise in the study of deep computations, transitions and deep computations themselves among them, send a state space $L\subseteq \mathbb{R}^\mathcal{P}$ into $\mathbb{R}^\mathcal{P}$.
The product space $\mathbb{R}^\mathcal{P}$ carries two natural $\sigma$-algebras: the Borel $\sigma$-algebra, generated by the open subsets of $\mathbb{R}^\mathcal{P}$, and the cylinder $\sigma$-algebra.
The latter, generated by the subbasic open sets $\pi_P^{-1}(U)$ with $P\in\mathcal{P}$ and $U\subseteq \mathbb{R}$ open, is in general strictly smaller than the Borel $\sigma$-algebra, but the two $\sigma$-algebras coincide when $\mathcal{P}$ is countable.
We use the cylinder $\sigma$-algebra to define universally measurable maps into $\mathbb{R}^\mathcal{P}$, because of the following characterization.

\begin{prop}\label{prop:coordinated_univ_meas}
  Let $X$ be a Hausdorff space and $Y=\prod_{i\in I}Y_i$ be any product of measurable spaces $(Y_i,\Sigma_i)$ for $i\in I$.
  Let $\Sigma_Y$ be the cylinder $\sigma$-algebra generated by the measurable spaces $(Y_i,\Sigma_i)$, and let $\pi_i\colon Y\to Y_i$ be the projections.
  The following are equivalent for  $f:X\rightarrow Y$:
  \begin{enumerate}[(i)]
  \item
    The map $f:X\rightarrow Y$ is universally measurable (with respect to $\Sigma_Y$).
  \item
    The map $\pi_i\circ f:X\rightarrow Y_i$ is universally measurable (with respect to $\Sigma_i$) for all $i\in I$.
  \end{enumerate}
\end{prop}

\begin{proof}
  The implication (i)$\Rightarrow$(ii) holds because each projection $\pi_i$ is measurable and $(\pi_i\circ f)^{-1}(E)=f^{-1}(\pi_i^{-1}(E))$ for every $E\in\Sigma_i$.
  To prove (ii)$\Rightarrow$(i), fix a Radon measure $\mu$ on $X$, and note that the sets $C\subseteq Y$ for which $f^{-1}(C)$ is $\mu$-measurable form a $\sigma$-algebra.
  This $\sigma$-algebra contains every measurable cylinder $C=\prod_{i\in I}C_i$:
  if $J$ is the finite set of indices $i\in I$ with $C_i\neq Y_i$, then $f^{-1}(C)=\bigcap_{i\in J}(\pi_i\circ f)^{-1}(C_i)$ is $\mu$-measurable by~(ii).
  Since the measurable cylinders generate $\Sigma_Y$, the set $f^{-1}(C)$ is $\mu$-measurable for every $C\in\Sigma_Y$.
\end{proof}

By Proposition~\ref{prop:coordinated_univ_meas}, the universal measurability of a map into $\mathbb{R}^\mathcal{P}$ can be checked feature by feature, just as membership in the first Baire class can (see Proposition~\ref{prop:Baire-identifications}).

We saw in Section~\ref{sec:prelim} that the NIP corresponds to relative compactness in the space of functions of Baire class~1 (see Theorems~\ref{thm:BFT}, \ref{thm:BFT2} and~\ref{thm:Generalized-BFT}).
In Section~\ref{sec:Classification} we saw, moreover, that this correspondence extends to deep computations when the set $\mathcal{P}$ of predicates is countable (see Theorem~\ref{thm:nip-baire1def}).
The following theorem, which holds on every Hausdorff space, characterizes the NIP in terms of universal measurability.
Applied to the type shards, the theorem therefore requires no countability of~$\mathcal{P}$.

\begin{thm}
  [Bourgain--Fremlin--Talagrand {\cite[Theorem~2F]{BFT_1978_PCompactBaire}}]
  \label{thm:BFT-2F}
  Let $X$ be a Hausdorff space, and let $A\subseteq \Cp(X)$ be pointwise bounded.
  The following are equivalent:
  \begin{enumerate}[(i)]
  \item
    The closure $\overline{A}$ is contained in~$M_r(X)$.
  \item
    For every compact $K\subseteq X$, the restriction $A|_K$ satisfies the NIP.
  \item
    For every Radon measure $\mu$ on $X$, the set $A$ is relatively countably compact in~$\mathscr{M}^0(X,\mu)$, i.e., every countable subset of $A$ has a limit point in $\mathscr{M}^0(X,\mu)$.
  \end{enumerate}
\end{thm}

Theorem~\ref{thm:BFT-2F} leads to a measure-theoretic counterpart of computability, the universal Monte Carlo computability of Definition~\ref{def:universally-essentially-computable} below, whose motivation lies in two classical facts.

\begin{enumerate}
\item
  Littlewood’s second principle, that every Lebesgue measurable function is “nearly continuous”, is made precise by Luzin’s theorem.
  Let $(X,\Sigma ,\mu )$ be a finite Radon measure space, and let $Y$ be a second-countable topological space (e.g., $Y=\mathbb{R}^\mathcal{P}$ with $\mathcal{P}$ countable), equipped with its Borel $\sigma$-algebra.
  Then a function $f\colon X\to Y$ is measurable if and only if, for every $E\in\Sigma$ and every $\varepsilon>0$, there is a closed set $F\subseteq E$ such that the restriction $f|_F$ is continuous and $\mu(E\setminus F) < \varepsilon$.
\item
  Computability of a transform is characterized by its continuous extendibility to the space of types (see Section~\ref{sec:computability});
  this characterization is at the core of the first paper of this series~\cite{alva2024approximability}.
\end{enumerate}
Taken together, the two facts justify Definition~\ref{def:universally-essentially-computable} below.
Indeed, suppose that an extension of a transform is universally measurable on a type shard, and fix a Radon probability measure on that shard.
By Luzin’s theorem, each output feature of that extension is then continuous on a closed subset of the shard whose measure is as close to~$1$ as desired.
On that subset, by the Stone--Weierstrass theorem, the output feature in question is a uniform limit of polynomials in the features.

\begin{defn}
  \label{def:universally-essentially-computable}
  Let $(L,\mathcal{P},\Gamma)$ be a CCS, and let $\mathscr{R}$ be a collection of sizers.
  A transform $f\colon L\to\mathcal{L}_{\sh}$ is \emph{universally Monte Carlo computable on~$\mathscr{R}$} if it has an extension $\tilde f\colon\mathcal{L}_{\sh}\to \mathcal{L}_{\sh}$ such that, for every $r_\bullet\in \mathscr{R}$, there is a sizer $s_\bullet$ for which the restriction $\tilde f|_{\mathcal{L}[r_{\bullet}]}\colon\mathcal{L}[r_\bullet]\to\mathcal{L}[s_\bullet]$ is universally measurable.
  Here universal measurability means that, for every Radon probability measure $\mu$ on $\mathcal{L}[r_\bullet]$ and every $P\in\mathcal{P}$, the map $\pi_P\circ \tilde f|_{\mathcal{L}[r_{\bullet}]}\colon\mathcal{L}[r_\bullet]\to [-s_P,s_P]$ is $\mu$-measurable.
\end{defn}

The deep computations of a set of computations confined by an exhaustive collection of sizers are transforms $L\to\mathcal{L}_{\sh}$ (see Proposition~\ref{prop:confined-compact}).
Since each type shard $\mathcal{L}[r_\bullet]$ is compact, Remark~\ref{rem:probability} shows that Radon probability measures suffice in Definition~\ref{def:universally-essentially-computable}, in place of general Radon measures.

\begin{thm}
  \label{thm:NIP-Monte-Carlo}
  Let $(L,\mathcal{P},\Gamma)$ be a CCS satisfying the Extendibility Axiom, let $\mathscr{R}$ be an exhaustive collection of sizers, and let $\Delta\subseteq\Gamma$ be $\mathscr{R}$-confined.
  The following are equivalent:
  \begin{enumerate}[(i)]
  \item
    for all $P\in\mathcal{P}$ and all $r_{\bullet}\in \mathscr{R}$, the set $\pi_P\circ\Delta|_{L[r_\bullet]}$ satisfies the NIP, where $\pi_P\colon\mathbb{R}^{\mathcal{P}}\to\mathbb{R}$ is the projection onto the coordinate~$P$;
  \item
    for every $r_\bullet\in\mathscr{R}$, every member of $\overline{\tilde\Delta|_{\mathcal{L}[r_\bullet]}}$ is universally measurable, as a map into $\mathbb{R}^{\mathcal{P}}$ with the cylinder $\sigma$-algebra.
  \end{enumerate}
  If these conditions hold, then every deep computation in $\overline{\Delta}$ is universally Monte Carlo computable on~$\mathscr{R}$.
\end{thm}

\begin{proof}
  Fix $r_\bullet\in \mathscr{R}$ and $P\in\mathcal{P}$, and write $K=\overline{\tilde\Delta|_{\mathcal{L}[r_\bullet]}}$, a compact subset of $\mathcal{L}[r_\bullet]^{\mathcal{L}[r_\bullet]}$, and $A_P=\pi_P\circ\tilde\Delta|_{\mathcal{L}[r_\bullet]}$.
  By the Extendibility Axiom and confinement, the set $A_P$ consists of continuous functions on the compact space $\mathcal{L}[r_\bullet]$ and is pointwise bounded.
  Moreover, since the map $g\mapsto\pi_P\circ g$ is continuous, the closure of~$A_P$ is the compact set $\{\pi_P\circ g:g\in K\}$.
  Assume first that (i) holds; then the restriction of $A_P$ to $L[r_\bullet]$, namely $\pi_P\circ\Delta|_{L[r_\bullet]}$, has the NIP.
  Since $L[r_\bullet]$ is dense in $\mathcal{L}[r_\bullet]$, Lemma~\ref{lem:NIP-and-closure} shows that the set $A_P$ itself has the NIP, and so does its restriction to every compact subset of~$\mathcal{L}[r_\bullet]$.
  Theorem~\ref{thm:BFT-2F} then gives the inclusion $\overline{A_P}\subseteq M_r(\mathcal{L}[r_\bullet])$, that is, every member of $\overline{A_P}$ is universally measurable.
  Since $P$ is arbitrary, Proposition~\ref{prop:coordinated_univ_meas} shows that every member of~$K$ is universally measurable, and (ii) follows.
  Conversely, if (ii) holds, then every member of $\overline{A_P}$ is universally measurable by Proposition~\ref{prop:coordinated_univ_meas}.
  By the implication (i)$\Rightarrow$(ii) of Theorem~\ref{thm:BFT-2F}, applied to the compact set $\mathcal{L}[r_\bullet]$ itself, the set $A_P$ has the NIP.
  Consequently, so does its restriction $\pi_P\circ\Delta|_{L[r_\bullet]}$, and (i) holds.
  Finally, let $f\in\overline{\Delta}$, and let $\tilde f$ be the extension of Lemma~\ref{lem:extension}.
  For every $r_\bullet\in\mathscr{R}$, the restriction $\tilde f|_{\mathcal{L}[r_\bullet]}$ belongs to $\overline{\tilde\Delta|_{\mathcal{L}[r_\bullet]}}$, and it is therefore universally measurable when (ii) holds.
  Since $\tilde f$ maps $\mathcal{L}[r_\bullet]$ into itself, the deep computation $f$ satisfies Definition~\ref{def:universally-essentially-computable} with $s_\bullet=r_\bullet$ for every $r_\bullet\in \mathscr{R}$.
\end{proof}

\begin{rem}\label{rem:doubling-rationals}
Condition~(ii) of Theorem~\ref{thm:NIP-Monte-Carlo} concerns all the members of the sets $\overline{\tilde\Delta|_{\mathcal{L}[r_\bullet]}}$, that is, the restrictions of ultracomputations to the type shards.
By contrast, the equivalence of the theorem fails if (ii) is weakened to the universal Monte Carlo computability of every deep computation, which requires only that some extension be universally measurable.
To see this, let $L=\mathbb{Q}/\mathbb{Z}\subseteq\mathbb{T}$ carry the predicates $\cos2\pi x$ and $\sin2\pi x$, and let $\Gamma$ consist of the iterates of the doubling map $x\mapsto2x$.
The space of types is the circle~$\mathbb{T}$, and each computation extends to the continuous map $x\mapsto2^nx$ of~$\mathbb{T}$.
Thus the Extendibility Axiom holds, and the single sizer with both bounds equal to~$1$ forms an exhaustive collection that confines~$\Gamma$.
By the paragraph that follows Proposition~\ref{prop:Newton-wild}, the iterates of the doubling map of~$\mathbb{T}$ do not have the NIP;
by Lemma~\ref{lem:NIP-and-closure}, neither do their restrictions to the dense set~$L$.
Every deep computation $f\colon L\to\mathbb{T}$, however, extends to the map that agrees with $f$ on~$L$ and takes the value~$0$ on $\mathbb{T}\setminus L$.
Since $L$ is countable, this extension is Borel, whence every deep computation is universally Monte Carlo computable.
\end{rem}

\subsection{Talagrand stability and Monte Carlo computability of deep computations}

The second route to universal Monte Carlo computability starts from Talagrand stability, a notion that Talagrand introduced in his study of Pettis integration~\cite{talagrand1984pettis} and that is closely related to the NIP.
We follow this route here and compare the two routes at the end of the subsection.
Let $X$ be a compact Hausdorff space, let $A\subseteq \mathbb{R}^X$ be a set of functions, and let $\mu$ be a Radon probability measure on $X$.
For a $\mu$-measurable set $E\subseteq X$, a positive integer $k$ and real numbers $a<b$, we write the points of $E^{2k}$ as $x=(x_0,\dots,x_{2k-1})$ and define
\[
  D_k(A,E,a,b)=\bigcup_{f\in A}\{x\in E^{2k}:f(x_{2i})\leq a, \hspace{1mm} f(x_{2i+1})\geq b \hspace{1mm}{\text{ for all }i<k} \}.
\]

\begin{defn}
  \label{def:Talagrand-stable}
  The set $A$ is \emph{Talagrand $\mu$-stable} if, for every $\mu$-measurable set $E\subseteq X$ of positive measure and all reals $a<b$, there is an integer $k\geq 1$ such that
  \[
    (\mu^{2k})^*(D_k(A,E,a,b))<(\mu(E))^{2k}.
  \]
  Here $(\mu^{2k})^*$ denotes the outer measure associated with the product measure $\mu^{2k}$ on~$X^{2k}$.
\end{defn}

The outer measure is needed in Definition~\ref{def:Talagrand-stable} because $D_k(A,E,a,b)$ need not be measurable.
The set $D_k(A,E,a,b)$ is measurable when $A$ is a countable set of $\mu$-measurable functions, since it is then a countable union of measurable sets.

\begin{rems}
  \label{rem:Talagrand-stability}
  Let $X$, $\mu$, $A\subseteq \mathbb{R}^X$ and $E\subseteq X$ be as above, and let $k\ge1$ be an integer and $a<b$ be reals.
  The following hold.
  \begin{enumerate}
  \item
    If $A\subseteq B\subseteq \mathbb{R}^X$, then
    \[
      D_k(A,E,a,b)\subseteq D_k(B,E,a,b).
    \]

  \item
    If $a<a'<b'<b$, then
    \[
      D_k(\overline{A},E,a,b)\subseteq D_k(A,E,a',b').
    \]
    [\emph{Proof}:
      If $x\in D_k(\overline{A},E,a,b)$, then there is $f\in\overline{A}$ such that $f(x_{2i})\leq a<a'$ and $f(x_{2i+1})\geq b>b'$ for all $i<k$.
      By the definition of the topology of pointwise convergence, there is $g\in A$ such that $g(x_{2i})<a'$ and $g(x_{2i+1})>b'$ for all $i<k$;
      hence $x\in D_k(A,E,a',b')$.]

  \item
    Every subset of a Talagrand $\mu$-stable set is Talagrand $\mu$-stable, by~(1).

  \item
    If $A$ is Talagrand $\mu$-stable, then so is $\overline{A}$, by~(2).
  \end{enumerate}
\end{rems}

A set of functions on a compact Hausdorff space is \emph{universally Talagrand stable} if it is Talagrand $\mu$-stable for every Radon probability measure~$\mu$ on that space (see Definition~\ref{def:universal-Talagrand} below).
In Corollary~\ref{cor:Talagrand-stab-Monte-Carlo} below we consider a confined set $\Delta$ of computations such that, for every feature $P$ and every type shard $\mathcal{L}[r_\bullet]$, the set $\pi_P\circ\tilde\Delta|_{\mathcal{L}[r_\bullet]}$ is universally Talagrand stable.
We show there that the deep computations of every such set are universally Monte Carlo computable.
The proof of the corollary rests on Lemma~\ref{lem:talagrand-stab-relcpct}, by which the pointwise closure of a Talagrand $\mu$-stable set consists of $\mu$-measurable functions.
The lemma, in turn, rests on the following characterization of measurable functions (compare Fact~\ref{baire}).

\begin{fact}[{\cite[Lemma~413G]{fremlin2003vol4}}]\label{fact:meas-fns-charn}
  Let $(X,\Sigma,\mu)$ be a measure space, and let $\mathcal{K}\subseteq\Sigma$ be a collection of measurable sets such that the following conditions hold:
  \begin{itemize}
  \item
    The measure space $(X,\Sigma,\mu)$ is complete, i.e., for all $E\in\Sigma$ with $\mu(E)=0$ and $F\subseteq E$ we have $F\in\Sigma$.
  \item
    It is semi-finite, i.e., for all $E\in\Sigma$ with $\mu(E)=\infty$ there exists $F\subseteq E$ such that $F\in\Sigma$ and $0<\mu(F)<\infty$.
  \item
    It is saturated, i.e., a set $E\subseteq X$ belongs to~$\Sigma$ if and only if $E\cap F\in\Sigma$ for all $F\in\Sigma$ with $\mu(F)<\infty$.
  \item
    It is inner regular with respect to $\mathcal{K}$, i.e., for all $E\in\Sigma$
    \[
      \mu(E)=\sup\{\mu(K):K\in\mathcal{K} \text{ and } K\subseteq E\}.
    \]
  \end{itemize}
  Then a function $f\colon X\to\mathbb{R}$ is measurable if and only if,
  for every $K\in\mathcal{K}$ with $0<\mu(K)<\infty$ and all reals $a<b$, either $\mu^*(K_{\le a})<\mu(K)$ or $\mu^*(K_{\ge b})<\mu(K)$,
  where $K_{\le a}=\{x\in K:f(x)\leq a\}$ and $K_{\ge b}=\{x\in K:f(x)\geq b\}$.
  For instance, let $X$ be a compact Hausdorff space, and let $\mu$ be a Radon probability measure on~$X$.
  If $\Sigma$ is the completion by~$\mu$ of the Borel $\sigma$-algebra of~$X$ and $\mathcal{K}$ is the collection of compact subsets of~$X$, then all four conditions hold.
\end{fact}

We first show that Talagrand stability passes to the pointwise closure and forces measurability.

\begin{lem}\label{lem:talagrand-stab-relcpct}
  Let $X$ be a compact Hausdorff space, let $\mu$ be a Radon probability measure on~$X$, and let $A\subseteq\mathbb{R}^X$ be a set of functions.
  If $A$ is Talagrand $\mu$-stable, then $\overline{A}$ is also Talagrand $\mu$-stable and $\overline{A}\subseteq\mathscr{M}^0(X,\mu)$.
\end{lem}

\begin{proof}
  By Remarks~\ref{rem:Talagrand-stability}, the closure $\overline{A}$ is Talagrand $\mu$-stable, and hence so is $\{f\}$ for every $f\in\overline{A}$.
  It therefore suffices to show that if $\{f\}$ is Talagrand $\mu$-stable, then $f$ is $\mu$-measurable.
  Suppose that $f$ is not $\mu$-measurable.
  Then Fact~\ref{fact:meas-fns-charn} gives a compact set $E$ of positive measure and reals $a<b$ with $\mu^*(E_{\le a})=\mu^*(E_{\ge b})=\mu(E)$, where $E_{\le a}=\{x\in E: f(x)\leq a\}$ and $E_{\ge b}=\{x\in E: f(x)\geq b\}$.
  For every $k\ge1$ we have $(E_{\le a}\times E_{\ge b})^k\subseteq D_k(\{f\},E,a,b)$, and this product has outer measure $(\mu(E))^{2k}$ (see~\cite{talagrand1984pettis});
  hence $\{f\}$ is not Talagrand $\mu$-stable.
\end{proof}

\begin{defn}\label{def:universal-Talagrand}
  Let $X$ be a compact Hausdorff space.
  A set $A\subseteq\mathbb{R}^X$ is \emph{universally Talagrand stable} if it is Talagrand $\mu$-stable for every Radon probability measure $\mu$ on~$X$.
\end{defn}

For pointwise bounded sets of continuous functions on a compact space, universal Talagrand stability implies the NIP.

\begin{prop}\label{prop:Talagrand-implies-NIP}
  Let $X$ be a compact Hausdorff space, and let $A\subseteq \Cp(X)$ be pointwise bounded.
  If $A$ is universally Talagrand stable, then $A$ satisfies the NIP.
\end{prop}

\begin{proof}
  By Lemma~\ref{lem:talagrand-stab-relcpct}, the closure $\overline{A}$ is contained in $\mathscr{M}^0(X,\mu)$ for every Radon probability measure $\mu$ on~$X$, that is, $\overline{A}\subseteq M_r(X)$.
  By the implication (i)$\Rightarrow$(ii) of Theorem~\ref{thm:BFT-2F}, applied with $K=X$, the set $A$ satisfies the NIP.
\end{proof}

\begin{cor}
  \label{cor:Talagrand-stab-Monte-Carlo}
  Let $(L,\mathcal{P},\Gamma)$ be a CCS satisfying the Extendibility Axiom, let $\mathscr{R}$ be an exhaustive collection of sizers, and let $\Delta\subseteq\Gamma$ be $\mathscr{R}$-confined.
  If $\pi_P\circ\tilde\Delta|_{\mathcal{L}[r_\bullet]}$ is universally Talagrand stable for all $P\in\mathcal{P}$ and all $r_{\bullet}\in \mathscr{R}$, then every deep computation in $\overline{\Delta}$ is universally Monte Carlo computable on~$\mathscr{R}$.
\end{cor}

\begin{proof}
  Let $P\in\mathcal{P}$ and $r_\bullet\in \mathscr{R}$ be given.
  By the Extendibility Axiom and confinement, the set $\pi_P\circ\tilde\Delta|_{\mathcal{L}[r_\bullet]}$ is a pointwise bounded family of continuous functions on the compact space~$\mathcal{L}[r_\bullet]$.
  Proposition~\ref{prop:Talagrand-implies-NIP} therefore shows that this set has the NIP, and hence so does its restriction $\pi_P\circ\Delta|_{L[r_\bullet]}$.
  Theorem~\ref{thm:NIP-Monte-Carlo} now completes the proof.
\end{proof}

In the context of deep computations we have thus found two routes to universal Monte Carlo computability, the NIP and universal Talagrand stability, and the second implies the first (see Proposition~\ref{prop:Talagrand-implies-NIP}).
For countable sets of measurable functions, the following two results show that whether the NIP implies universal Talagrand stability depends on the axioms of set theory.
On the one hand, it is consistent with the Zermelo--Fraenkel axioms of set theory with the axiom of choice (ZFC) that, for such sets, the NIP implies universal Talagrand stability.

\begin{thm}[Talagrand {\cite[Theorem~9-3-1(a)]{talagrand1984pettis}}]
  \label{thm:NIP-implies-Talagrand-stab}
  Let $X$ be a compact Hausdorff space, and let $A\subseteq M_r(X)$ be countable and pointwise bounded.
  Assume that $[0,1]$ is not the union of fewer than $2^{\aleph_0}$ closed sets of Lebesgue measure zero.
  If $A$ satisfies the NIP, then $A$ is universally Talagrand stable.
\end{thm}

(The hypothesis on $[0,1]$ follows, for example, from the Continuum Hypothesis.)

On the other hand, for Lebesgue measure, the failure of this implication is also consistent.

\begin{thm}[Fremlin and Shelah~\cite{fremlin1993pointwise}]
  It is consistent with ZFC that there exists a countable pointwise bounded set of Lebesgue measurable functions with the NIP that is not Talagrand stable with respect to Lebesgue measure.
\end{thm}

Both results concern countable sets of measurable functions.
For countable sets of \emph{continuous} functions on a compact space, by contrast, the NIP implies universal Talagrand stability in ZFC, without any further set-theoretic hypothesis (see Theorem~\ref{thm:NIP-Talagrand-countable} below).
Consequently, since the Extendibility Axiom makes computations continuous on the type shards, the two routes to universal Monte Carlo computability coincide for countable sets of computations.
After a reduction to a compact metrizable space, this implication for countable sets of continuous functions can be derived from two results on the Glivenko--Cantelli property of a family of functions.
The Glivenko--Cantelli property is the almost sure convergence of empirical means to expectations, uniformly over the family.
The first result is a theorem of van Handel~\cite[Theorem~1.3]{vanHandel2013}, which concerns separable, uniformly bounded families of measurable functions on a standard Borel space.
According to this theorem, such a family has the Glivenko--Cantelli property for every probability measure if and only if it contains no Boolean $\sigma$-independent sequence.
The second result is a theorem of Talagrand~\cite{talagrand1987glivenko}, which identifies, among uniformly bounded families, the Glivenko--Cantelli families for a measure~$\mu$ with the Talagrand $\mu$-stable ones (see also~\cite[\S465]{fremlin2003vol4}).
We give instead a direct proof, which reduces the question to a compact metrizable space and then applies Fact~\ref{baire} to a cluster point of the reduced family.

\begin{thm}\label{thm:NIP-Talagrand-countable}
Let $X$ be a compact Hausdorff space, and let $A\subseteq\Cp(X)$ be countable and pointwise bounded.
If $A$ satisfies the NIP, then $A$ is universally Talagrand stable.
\end{thm}

\begin{proof}
Suppose that $A$ satisfies the NIP, and let $\mu$ be a Radon probability measure on~$X$.
Towards a contradiction, suppose further that $A$ is not Talagrand $\mu$-stable.
Then there are a $\mu$-measurable set $E\subseteq X$ of positive measure and reals $a<b$ such that $\mu^{2k}(D_k(A,E,a,b))=(\mu(E))^{2k}$ for every $k\ge1$.
The sets $D_k(A,E,a,b)$ are measurable, since $A$ is countable and its members are continuous.
We first transfer the situation to a compact metrizable space.
The map $\Psi\colon X\to\mathbb{R}^A$ given by $\Psi(x)=(f(x))_{f\in A}$ is continuous, and its image $Y=\Psi(X)$ is compact and metrizable because $A$ is countable.
Moreover, each $f\in A$ factors as $f=\hat f\circ\Psi$, where $\hat f\colon Y\to\mathbb{R}$ is the restriction of the coordinate projection.
The set $\hat A=\{\hat f:f\in A\}\subseteq\Cp(Y)$ is pointwise bounded, and it satisfies the NIP, because $Y^{\hat f_I}_{\le a'}\cap Y^{\hat f_J}_{\ge b'}=\Psi(X^{f_I}_{\le a'}\cap X^{f_J}_{\ge b'})$ for every sequence $(f_n)_{n\in\mathbb{N}}$ in~$A$, all finite sets $I,J\subseteq\mathbb{N}$ and all reals $a'<b'$.
Let $\nu$ be the Borel probability measure on~$Y$ given by $\nu(B)=\mu(E\cap\Psi^{-1}(B))/\mu(E)$.
Since $D_k(A,E,a,b)=E^{2k}\cap(\Psi^{2k})^{-1}(D_k(\hat A,Y,a,b))$, where $\Psi^{2k}$ applies $\Psi$ to each coordinate, we obtain $\nu^{2k}(D_k(\hat A,Y,a,b))=1$ for every $k\ge1$.

Now consider the product measure $\nu^{\mathbb{N}}$ on $Y^{\mathbb{N}}$, under which the coordinates $y_n$ of a sequence $(y_n)_{n\in\mathbb{N}}$ are independent, each with distribution~$\nu$.
Let $Z$ denote the support of~$\nu$, which is a nonempty closed, and hence Polish, subspace of~$Y$.
For each $k\ge1$, almost every sequence admits some $\hat f\in\hat A$ with $\hat f(y_{2i})\le a$ and $\hat f(y_{2i+1})\ge b$ for all $i<k$, since this event has probability $\nu^{2k}(D_k(\hat A,Y,a,b))=1$.
Moreover, since $Y$ has a countable base and every basic open set that meets~$Z$ has positive measure, both $\{y_{2i}:i\in\mathbb{N}\}$ and $\{y_{2i+1}:i\in\mathbb{N}\}$ are dense subsets of~$Z$ for almost every sequence.
Fix a sequence $(y_n)$ with all these properties, and, for each $k\ge1$, choose $\hat f_k\in\hat A$ with $\hat f_k(y_{2i})\le a$ and $\hat f_k(y_{2i+1})\ge b$ for all $i<k$.
Since $\hat A$ is pointwise bounded, the sequence $(\hat f_k)$ has a cluster point $g$ in~$\mathbb{R}^Y$, and $g\in\mathrm{B}_1(Y)$ by Theorem~\ref{thm:BFT2}.
For every~$i$, the inequalities $\hat f_k(y_{2i})\le a$ and $\hat f_k(y_{2i+1})\ge b$ hold for all $k>i$, and therefore $g(y_{2i})\le a$ and $g(y_{2i+1})\ge b$.
By condition~(3) of Fact~\ref{baire}, applied to the Polish space~$Y$ and its closed subset~$Z$, the restriction $g|_Z$ has a point of continuity~$z$.
Every neighborhood of $z$ in~$Z$, however, contains points $y_{2i}$ with $g(y_{2i})\le a$ and points $y_{2i+1}$ with $g(y_{2i+1})\ge b$; this contradicts the continuity of $g|_Z$ at~$z$.
\end{proof}

\begin{cor}\label{cor:NIP-Talagrand-computations}
Let $(L,\mathcal{P},\Gamma)$ be a CCS satisfying the Extendibility Axiom, let $\mathscr{R}$ be an exhaustive collection of sizers, and let $\Delta\subseteq\Gamma$ be countable and $\mathscr{R}$-confined.
For all $P\in\mathcal{P}$ and $r_\bullet\in\mathscr{R}$, the set $\pi_P\circ\tilde\Delta|_{\mathcal{L}[r_\bullet]}$ is universally Talagrand stable if and only if $\pi_P\circ\Delta|_{L[r_\bullet]}$ satisfies the NIP.
\end{cor}

\begin{proof}
By the Extendibility Axiom and confinement, the set $\pi_P\circ\tilde\Delta|_{\mathcal{L}[r_\bullet]}$ is a countable, pointwise bounded family of continuous functions on the compact space~$\mathcal{L}[r_\bullet]$, in which $L[r_\bullet]$ is dense.
If this set is universally Talagrand stable, it has the NIP by Proposition~\ref{prop:Talagrand-implies-NIP}, and so does its restriction to~$L[r_\bullet]$.
Conversely, suppose that $\pi_P\circ\Delta|_{L[r_\bullet]}$ has the NIP.
By Lemma~\ref{lem:NIP-and-closure}, so does the set $\pi_P\circ\tilde\Delta|_{\mathcal{L}[r_\bullet]}$, to which Theorem~\ref{thm:NIP-Talagrand-countable} therefore applies.
\end{proof}

For countable sets of computations, Theorem~\ref{thm:NIP-Talagrand-countable} thus supersedes Theorem~\ref{thm:NIP-implies-Talagrand-stab} and its set-theoretic hypothesis.
Accordingly, the example of Fremlin and Shelah, a set of measurable functions, has no counterpart among countable sets of computations.

\section{Open questions}
\label{sec:open-questions}

We close with the questions that remain open, in the order in which they arise from the sections of the paper.
The first concerns the hierarchy of Section~\ref{sec:trichotomy}, whose levels were defined through the topology of the space of deep computations.
Proposition~\ref{prop:levels-tests} describes these levels by the number of tests needed to tell deep computations apart, but the descriptions so obtained still involve infinitely many tests.

\begin{question}\label{q:NIP-levels}
  Is there a finitary, combinatorial characterization of $\NIP_1$, $\NIP_2$ and $\NIP_3$, in the spirit of the characterization of $\NIP$ by finite patterns in Theorem~\ref{thm:nip-baire1def}?
\end{question}

In Section~\ref{sec:iterative} we place a single iterated map at every level of the hierarchy except one.
The exception is the level $\NIP$ but not $\NIP_1$, at which both of our examples are actions of groups with two generators (see Section~\ref{sec:ex-NIP}).

\begin{question}\label{q:cascade-NIP1}
  Is there a cascade CCS that is $\NIP$ but not $\NIP_1$?
  Equivalently, is there a continuous map $\gamma$ of a compact metrizable space $L$ such that $(L,\gamma)$ is tame and $E(L,\gamma)$ is not first countable?
\end{question}

For actions of $\mathbb{Z}$, and more generally of abelian groups, the analogous question is Question~6.8 of Glasner and Megrelishvili~\cite{glasner2022tame}.

At the level $\NIP_1$ but not $\NIP_2$, the deep computations of the three-dimensional Bresenham algorithm contain a copy of $\hat D(S(2^{\mathbb{N}}))$ and none of $D(2^{\mathbb{N}})$ (see Theorem~\ref{thm:Bresenham3}).
By contrast, the countable family $\Delta_{\mathcal{C}}$ of elements of the group generated by the steps of the Euclidean algorithm realizes the other alternative of Todorčević’s trichotomy (see Proposition~\ref{prop:euclid-duplicate}).

\begin{question}\label{q:cascade-duplicate}
  Is there a cascade CCS that is $\NIP_1$ but not $\NIP_2$ and whose space of deep computations contains a copy of $D(2^{\mathbb{N}})$?
\end{question}

\appendix

\section{The Thue--Morse system is not tame}
\label{app:thue-morse}

In this appendix we prove that the shift on the Thue--Morse subshift generates an action of~$\mathbb{Z}$ that is not tame, a fact observed by Glasner and Megrelishvili~\cite{GlasnerMegrelishvili2018}.
From this fact we derive that the two dynamical CCSs defined by this shift, which we use in Section~\ref{sec:beyond-NIP} as examples beyond the hierarchy, are not $\NIP$.
Our proof is self-contained except for the structure theorem for tame minimal systems under actions of abelian groups, due to Huang, to Kerr and Li, and to Glasner.
We recall that theorem below and combine it with an elementary computation of the correlations of the Thue--Morse sequence.

We first recall the objects involved.
The Thue--Morse sequence $(t_n)_{n\in\mathbb{N}}$ is defined by letting $t_n$ be the parity of the number of ones in the binary expansion of~$n$, whence $t_{2k}=t_k$ and $t_{2k+1}=1-t_k$ for every $k\in\mathbb{N}$.
The same sequence is the limit of the words $W_m$, for $m\in\mathbb{N}$, obtained from $0$ by $m$ applications of the substitution $0\mapsto01$, $1\mapsto10$.
The two-sided sequences all of whose finite blocks occur in $(t_n)$ form the Thue--Morse subshift $X_{TM}\subseteq2^{\mathbb{Z}}$, a closed subset of the compact metrizable space $2^{\mathbb{Z}}$.
On this subshift the shift~$\sigma$, defined by $\sigma(x)_n=x_{n+1}$, is a homeomorphism.
As in Section~\ref{sec:dynamical-CCS}, a system $(L,\Gamma)$ consists of a compact metrizable space~$L$ and a countable semigroup~$\Gamma$ of continuous maps $L\to L$.
Such a system is \emph{tame}, in the sense of that section, if for every continuous function $f\colon L\to\mathbb{R}$ the family $\{f\circ\gamma:\gamma\in\Gamma\}$ has the $\NIP$ of Definition~\ref{def:NIP-fns}.
We call a sequence $(f\circ\gamma_k)_k$ of members of that family \emph{independent} if, for some reals $a<b$, it witnesses the finitary independence property of Definition~\ref{def:NIP-fns}.
Hence the system $(L,\Gamma)$ is tame exactly when no independent sequence exists.
The CCS with states $X_{TM}$, the coordinate predicates $x\mapsto x_n$ for $n\in\mathbb{Z}$, and the group $\{\sigma^n:n\in\mathbb{Z}\}$ is a dynamical CCS in the sense of Section~\ref{sec:dynamical-CCS}.
The cascade CCS of~$\sigma$, in turn, is the dynamical CCS with the same states and predicates but with the semigroup $\{\sigma^n:n\in\mathbb{N}\}$.
By Proposition~\ref{prop:cascade-tame}, each of these two CCSs is $\NIP$ if and only if its dynamical system is tame.

That the dynamics of~$\sigma$ on~$X_{TM}$ is wild follows from the structure theory of tame minimal systems under actions of abelian groups.
Huang~\cite{Huang2006} and, independently, Kerr and Li~\cite{KerrLi2007} proved that every metrizable minimal tame system under an action of an abelian group is uniquely ergodic.
Unique ergodicity means that the system carries exactly one invariant Borel probability measure.
The same authors proved, in addition, that such a system is an almost one-to-one extension of its maximal equicontinuous factor.
Glasner~\cite{Glasner2007} showed, moreover, that the factor map onto that factor is an isomorphism of measure-preserving systems.
For a minimal homeomorphism, the maximal equicontinuous factor is a rotation $y\mapsto y+\alpha$ of a compact abelian group~$Y$ in which the multiples of~$\alpha$ are dense, equipped with its Haar measure.
The characters of~$Y$ are eigenfunctions of this rotation and span a dense subspace of the space of square-integrable functions on~$Y$.
Consider now a minimal system given by a homeomorphism~$\sigma$ that generates a tame action.
The Koopman operator of this system sends each function that is square-integrable with respect to the invariant measure to its composite with~$\sigma$.
By the results just recalled, every such function is a limit in mean square of linear combinations of eigenfunctions of the Koopman operator.
The Thue--Morse subshift fails this test, because the correlations of the Thue--Morse sequence are, on average, too small.

\begin{prop}\label{prop:thue-morse}
Let $X_{TM}\subseteq2^{\mathbb{Z}}$ be the Thue--Morse subshift, and let $\sigma$ be the shift on~$X_{TM}$.
The action of $\mathbb{Z}$ on $X_{TM}$ generated by~$\sigma$ is not tame.
Consequently, neither the CCS with states $X_{TM}$, the coordinate predicates and the group $\{\sigma^n:n\in\mathbb{Z}\}$ nor the cascade CCS of~$\sigma$ is $\NIP$.
\end{prop}

\begin{proof}
We first show that $X_{TM}$ is minimal.
Every finite block~$w$ of $(t_n)$ occurs in some word~$W_m$.
If $\overline{W}_m$ denotes the word obtained from $1$ by $m$ applications of the substitution, then the words $W_{m+1}=W_m\overline{W}_m$ and $\overline{W}_{m+1}=\overline{W}_mW_m$ both contain~$W_m$.
Moreover, since $t_{2k}=t_k$ and $t_{2k+1}=1-t_k$, the block of $(t_n)$ of length $2^{m+1}$ that begins at the index $2^{m+1}j$ is $W_{m+1}$ or $\overline{W}_{m+1}$ according as $t_j$ is $0$ or~$1$.
Hence the block~$w$ occurs in every block of length $2^{m+2}$ of $(t_n)$, and therefore in every point of~$X_{TM}$;
every orbit of~$\sigma$ is thus dense in~$X_{TM}$.

Suppose now that the action of $\mathbb{Z}$ generated by~$\sigma$ is tame.
By the structure theorem recalled above, the subshift $X_{TM}$ carries a unique $\sigma$-invariant Borel probability measure~$\mu$.
The eigenfunctions of the unitary operator~$V$ on $L^2(\mu)$ that sends each function to its composite with~$\sigma$ span, by the same theorem, a dense subspace of $L^2(\mu)$.
Let $\vartheta(x)=(-1)^{x_0}$ for $x\in X_{TM}$, and let $c(n)=\int \vartheta\cdot(\vartheta\circ\sigma^n)\,d\mu$ for $n\in\mathbb{N}$.
By minimality and compactness, the subshift $X_{TM}$ contains a point~$x$ with $x_j=t_j$ for every $j\ge0$.
Every weak$^*$ limit of the measures $\frac1K\sum_{k<K}\delta_{\sigma^kx}$ is $\sigma$-invariant and hence equal to~$\mu$, whence
\[
c(n)=\lim_{K\to\infty}\frac1K\sum_{k<K}(-1)^{t_k+t_{k+n}} .
\]
Splitting the sum according to the parity of~$k$ and using $t_{2k}=t_k$ and $t_{2k+1}=1-t_k$, we obtain $c(2n)=c(n)$ and $c(2n+1)=-\frac12\bigl(c(n)+c(n+1)\bigr)$ for every $n\in\mathbb{N}$, and in particular $c(1)=-\frac13$.
For $K\ge1$, let $A(K)=\sum_{n<K}c(n)^2$ and $B(K)=\sum_{n<K}c(n)c(n+1)$.
For every $K\ge1$, the recursions yield the identities and inequalities
\begin{align*}
B(2K)&=-\tfrac12\sum_{n<K}\bigl(c(n)+c(n+1)\bigr)^2\le0,\\
A(2K)&=A(K)+\tfrac14\sum_{n<K}\bigl(c(n)+c(n+1)\bigr)^2\le\tfrac32A(K)+\tfrac12B(K),
\end{align*}
where the last inequality follows by expanding the squares, because $c(0)=1\ge c(K)^2$.
Hence we have $A(4K)\le\frac32A(2K)$ for every $K\ge1$, and, since $A(2)=1+c(1)^2\le2$, it follows that $A(2^k)\le2\bigl(\frac32\bigr)^{k-1}$ for every $k\ge1$.
Since the function~$A$ is nondecreasing, we conclude that $A(K)/K\to0$ as $K\to\infty$ (see Figure~\ref{fig:thue-morse}).
Now let $\lambda$ be a complex number of modulus~$1$, and let $\Pi_\lambda$ be the orthogonal projection of $L^2(\mu)$ onto the kernel of $V-\lambda I$.
This kernel consists of the eigenfunctions of~$V$ with eigenvalue~$\lambda$ and the function~$0$.
By von Neumann's mean ergodic theorem, applied to the unitary operator $\bar\lambda V$, the averages $\frac1K\sum_{n<K}\bar\lambda^n\,\vartheta\circ\sigma^n$ converge to $\Pi_\lambda\vartheta$ in $L^2(\mu)$.
Consequently, since $c(n)$ is the inner product of $\vartheta\circ\sigma^n$ and~$\vartheta$, the Cauchy--Schwarz inequality yields
\[
\|\Pi_\lambda\vartheta\|^2=\lim_{K\to\infty}\Bigl|\frac1K\sum_{n<K}\bar\lambda^n c(n)\Bigr|\le\lim_{K\to\infty}\Bigl(\frac{A(K)}{K}\Bigr)^{1/2}=0 .
\]
Thus the function~$\vartheta$ is orthogonal to every eigenfunction of~$V$, and therefore to every function in $L^2(\mu)$, since these eigenfunctions span a dense subspace of $L^2(\mu)$.
This conclusion is absurd, since $|\vartheta|=1$.
The action is therefore not tame, and, by Proposition~\ref{prop:cascade-tame}, the CCS with the group $\{\sigma^n:n\in\mathbb{Z}\}$ is not $\NIP$.

\begin{figure}[htbp]
  \centering
  \includegraphics[width=\textwidth]{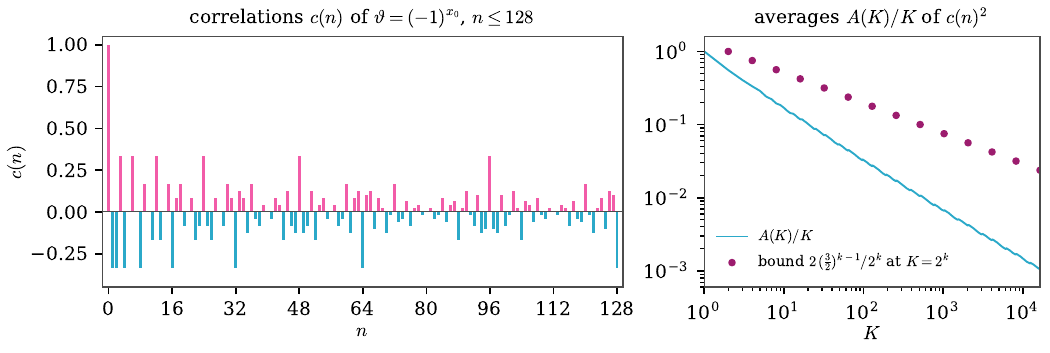}
  \caption{The correlations $c(n)$ of the function $\vartheta(x)=(-1)^{x_0}$ on the Thue--Morse subshift.
    They are computed exactly from the recursions $c(0)=1$, $c(2n)=c(n)$ and $c(2n+1)=-\frac12\bigl(c(n)+c(n+1)\bigr)$, which give $c(1)=-\frac13$.
    The left panel shows the values $c(n)$ for $0\le n\le128$.
    The right panel shows, on logarithmic scales, the averages $A(K)/K$ for $K\le2^{14}$, where $A(K)=\sum_{n<K}c(n)^2$.
    The dots mark the bound $2\bigl(\frac32\bigr)^{k-1}/2^k$ on $A(2^k)/2^k$.
    The decay of these averages is the key step in the proof of Proposition~\ref{prop:thue-morse}.}
  \label{fig:thue-morse}
\end{figure}

For the cascade CCS of~$\sigma$, we use a symmetry of the Thue--Morse sequence.
The words $W_{2m}$ are palindromes, and every finite block of $(t_n)$ occurs in one of them.
Consequently, the reversal $\iota$ of two-sided sequences, defined by $\iota(x)_n=x_{-n}$, maps $X_{TM}$ onto itself and satisfies $\iota\circ\sigma=\sigma^{-1}\circ\iota$.
Since the action of~$\mathbb{Z}$ is not tame, there are a continuous function $f\colon X_{TM}\to\mathbb{R}$ and pairwise distinct integers~$n_k$ such that the sequence $(f\circ\sigma^{n_k})_k$ is independent.
Moreover, infinitely many of the indices $n_k$ have the same sign.
If these indices are positive, the corresponding terms already form an independent sequence in the cascade.
If these indices are negative, the functions $f\circ\sigma^{n_k}\circ\iota=(f\circ\iota)\circ\sigma^{-n_k}$ form an independent sequence with positive exponents, because the reversal~$\iota$ is a bijection of $X_{TM}$.
In either case the cascade of~$\sigma$ is not tame, and, by Proposition~\ref{prop:cascade-tame}, its CCS is not $\NIP$.
\end{proof}

\bibliographystyle{amsalpha}
\bibliography{bib}

\end{document}